# AN EXPLANATORY NOTE

  A brief word about the plan of the proof. The basic idea is to show that if there is a minimal flow on a compact 3-manifold M, then there is a co-dimension one foliation on M transverse to the flow. In order to show this, I have partitioned M into subsets called tubes, each of which is homeomorphic to a cylinder D×I, where D is a 2-dimensional disc of points at distance ≤1 from the origin, on the plane and I = [0,1].

There is an inductive procedure described for foliating each tube with discs. Further, this foliation has to be done, for each disc, in a way that respects certain stringent boundary requirements − viz., that the procedure should not create an obstruction which we have termed a singular cycle, to the inductive foliation procedure at the succeeding stage. The core of the paper is devoted to showing this. A procedure is described, developed through a series of technical lemmas in Section 3, showing how this may be done without violating these boundary restrictions. Prior to this, in Section 2, we establish that there is a method of selecting a space $S_0$ on which to begin constructing these tubes, such that before beginning the inductive procedure, there are no singular cycles on $S_0$, thereby preparing the way for the start of the process.

However, I must also own to a certain tactical error I have committed. If I had attempted to show the non-existence merely of $C^1$ minimal flows, the whole proof could have been made materially simpler. In particular, the whole of Section 1, together with some complicated manoeuvres later on, would have been unnecessary. However, I overreached myself by establishing the result for $C^0$ minimal flows. In essence, therefore, I have been obliged to approximate the $C^0$ flow sufficiently closely by a smooth flow. This, as I have remarked, has complicated the proof somewhat. In view of this, you may, if you choose, begin reading the paper from Section 2, after a cursory glance at Section 1, or even skip it altogether.

Also, the references I have provided are, I am afraid, incomplete. It was my intention to provide a complete list, in the event of the paper being accepted for publication. But I thought that it would not be a good idea to delay sending the paper on this account. I have therefore restricted the references to the two results I have found it necessary to use.

## A plan of the proof

**SECTION 1:**  Pages 1 − 7   Here we show that a $C^0$ minimal flow can be approximated by a suitable smooth minimal flow.

**SECTION 2:**  Pages 8 − 24   Beginning with a disc transverse to the flow, we describe a construction on it to get a 2-manifold with boundary which is transverse to the flow, which is partitioned by its first return function into a finite number of regions on each of which the first return function is a homeomorphism. Using this space as a base, we are now ready to start the process of foliating $S^3$ by an inductive procedure.

**SECTION 3:** This section is devoted to the description of the foliating procedure, an overview of which is given in page 26. In a series of technical lemmas (Lemma 3.1 to 3.9), we establish the conditions necessary for the foliation procedure. In Theorem 3.1 we demonstrate the procedure by which the foliation is carried out. The second part of the theorem shows that the foliation thus constructed has no singular cycles, which form an obstruction to the foliation procedure at the next stage. Novikov's Theorem on the existence of compact leaves of co-dimension 1 foliations on $S^3$ now applies to give a contradiction, which proves the non-existence of minimal flows.

# The Gottschalk Conjecture

<u>Introduction.</u> The Gottschalk conjecture is one of the oldest problems relating to the geometry of vector fields on $S^3$. In some sense, the Gottschalk and the Seifert conjectures, between them, span the entire possible geometry of the behaviour of vector fields on $S^3$. Though the latter conjecture was shown to be false, the Gottschalk Conjecture, proposed in 1948, has resisted all attempts at solution. With the exception of Katok's result that there exists a smooth vector field on $S^3$ almost every integral curve of which is dense, there have been very few partial results. In what follows we shall establish the Gottschalk Conjecture by showing that there can be no $C^0$ vector field on $S^3$ with every integral curve being dense in $S^3$.

## SECTION 1

### Two Fundamental Principles

The two following principles are fundamental consequences of a flow being minimal and follow from the definition. We shall be using them repeatedly in Sections 1 and 2.

**FP I :** Let v be a minimal flow on a 3-manifold M and X a 2-manifold transverse to v in M. Then X intersects every flow line of v.

**FP II** : If Q is a point of X, there is a point P of X such that the flow line through P, when produced in the forward direction of the flow, reaches Q.

**Convention**: All vector fields in the sequel will be non-singular.

In this section, we begin with a $C^0$ minimal flow on $S^3$. In the discussion that follows, we will be using freely geometric terms like parallelism and orthogonality. In defence of appropriating such Euclidean terms to $S^3$, it can be observed that stereographic projection allows us to keep moving from $R^3$ to $S^3$-{a distinguished point}.

Let v be a minimal $C^0$ flow on $S^3$. The basic idea of the construction is to keep smoothing out the flow, by means of isotopies, inside certain tubular sets which we shall be describing presently. We begin by describing the first-hit function. Let S be a compact two manifold in $S^3$. We say that S is transverse to the flow v if there is an open set N containing S such that inside N, v and S are conjugate by a homeomorphism H to a topological situation in which H(S) is transverse H(v), where the former is a smooth imbedded surface and the latter is a smooth flow.

<u>Definition 1.1:</u>  <u>The first-hit function</u> Let S and $S_1$ be two compact 2-manifolds (generally with boundary), transverse to the flow. The first-hit function h : S $\longrightarrow$ $S_1$ (denoted by h(S,$S_1$)) is defined as follows. For x ε S, we consider the flow line through x in the positive direction of the flow v. If y is the point in which this flow line first intersects $S_1$, then h(x) = y. Since we will be concerned with a minimal flow in the sequel, we can be sure that a point like y will always exist. h will always stand for such a first-hit function. We may also refer to h as a first return function when S = $S_1$ and in this case denote it by h(S). h[S,$S_1$] will denote the image h(S) in $S_1$.

<u>Definition 1.2: Tube</u>. Let $D_1$, $D_2$ be two discs transverse to the flow ( where disc will always mean a closed two dimensional disc transverse to the flow) such that h($D_1$, $D_2$) is a homeomorphism onto  $D_2$. ( In this sense it will sometimes be useful to look upon $D_2$ as being obtained from $D_1$ when the latter is moved along  the flow to the position occupied by $D_2$). In such a case we call the union of the closed line segments between $D_1$ and $D_2$ a tube and denote it by T($D_1$,$D_2$). ( If $D_1$ and $D_2$ are two manifolds transverse to the flow, and h($D_1$,$D_2$) is a homeomorphism not necessarily onto $D_2$, then the same symbol $T(D_1, D_2)$ will denote the tube $T(D_1, h(D_1))$.

<u>Definition 1.3: Initial and end faces.</u> For a tube $T(D_1, D_2)$ as above with $h(D_1, D_2)$ a homeomorphism onto, $D_1$ and $D_2$ will be called the initial and end faces respectively, of the tube. This symbol T may also be used to indicate a tube where there is no possibility of confusion.

<div align="center">Fig 1</div>

Before we begin the smoothing construction referred to earlier, we assume that the flow has been "straightened out", by means of suitable isotopies, in a small open set, so that the segments of the flow lines in this open set are all parallel straight line segments which are orthogonal to a circular disc D contained in an open disc W. For $x \, \varepsilon \, D$ we consider $h(x)$, where h is the first return function on D. There are two possibilities.

<u>Case 1:</u> $h(x) \, \varepsilon \, \text{int} \, D$. In this case, let $V_x$ be a closed disc neighbourhood of x in W, such that h carries $V_x$ homeomorphically into int D, with $V_x \cap h(V_x) = \phi$. Let $T_x$ denote the tube $T(V_x, D)$. We now isotope the flow lines in $T_x$ within an $\varepsilon$-neighbourhood N of $T_x$, to a flow $v_1$ which is smooth in a neighbourhood of $T_x$, while the isotopy does not move any point in $V_x$ or outside N.

<div align="center">Fig 2</div>

<u>Case 2:</u> $h(x) \, \varepsilon \, \partial D$. Then because v is minimal, there is a positive integer n such that the forward or positive orbit through x intersects $\partial D$ n times before returning to int $D$ – i.e $h^k(x) \, \varepsilon \, \partial D$ for $1 \leq k \leq n$ and $h^{n+1}(x) \, \varepsilon \, \text{int} \, D$. Let $h^i(x) = x_i$, $1 \leq i \leq n+1$ . As before, let $V_x$ be a small closed disc

<div align="center">5</div>

neighbourhood of x in W and $T_x$ the tube with $V_x$ and $h^{n+1}(V_x)$ (h is now the first return function on W) , as initial and end faces. $V_x$ is chosen small enough so that $T_x \cap W$ consists of  n+2 pairwise disjoint discs $D_i$, i = 1,2,…. n+2 in the interior of W. Again, we isotope the flow lines in a neighbourhood of $T_x$ to a flow $v_1$ such that if $h_1$ is the new first return function on W with respect to $v_1$, then :

(i)$h_1$/ $D_i$, $1 \leq i \leq n +1$ is a smooth diffeomorphism.

(ii) $h_1$ / ($D_i \cap \partial D$) is transversal to $\partial D$, for $1 \leq i \leq n+1$

(iii) $v_1$ is smooth in a neighbourhood of $T_x$.

It is to be noted that to achieve the second requirement it may be necessary to redraw pieces of $\partial D$ within ( $T_x \cap W$) – i.e, in effect choosing a new disc D. A more detailed description of this procedure can be found in [1].

[ Convention :  As a consequence of the isotopy, the flow v has been deformed into another flow $v_1$. We shall, however, from this point on, continue to denote it by the same symbol v, since the essential geometric feature – viz, minimality – is preserved. This understanding will also govern our use of other symbols regarding regions, sets, etc. Even though these may be deformed into new positions, in order to avoid a proliferation of symbols, we will continue to denote them by the same symbols. In case we need to distinguish between, say, the first-hit function before and after a particular deformation is effected, we shall denote the functions associated with these two situations by different symbols when required to do so].

Next, for x ε D and h(x) ∈ int D, let us assume that we have isotoped the flow in a neighbourhood of a tube whose initial face U contains x in its interior, so that h : U ----→ D is a smooth diffeomorphism into the interior of D. Let $U_0$ be a disc in D or W (respectively represented by (a) or (b) in the diagram below).

Fig 3

It can be seen that the flow can be isotoped in a neighbourhood of $T(U_0)$, (without moving points of h(U)), so that h : U $\cup$ $U_0$ --------$\rightarrow$ int D is smooth. It is not difficult to see how the procedure should be modified (as suggested in case 2 above), when h(x) $\varepsilon$ $\partial$D or h($U_0$) intersects $\partial$D. The fundamental idea behind the procedure here is to keep enlarging the domain of smoothness of h over successively larger regions obtained by considering adjacent neighbourhoods on D.

Definition1.4: The procedure by which the flow lines in the tube T(U) were isotoped to make h smooth on U, will be called a tube modification.

But one point of difficulty needs to be cleared up. There may be a neighbourhood $U_m$ such that for h : D --------$\rightarrow$ D, h($U_m$) may intersect one of the $U_i$'s for i < m. For instance, let h($U_m$) $\cap$ U $\neq$ $\phi$. Thus, the process of smoothing out h in a neighbourhood of $U_m$ may now require some points of U to be moved by this isotopy. In general, this will mean that the new first return function will no longer be smooth at the displaced points of U. But the differentiability of h at these points can be restored by a fresh isotopy of U at these points. For example, let h($U_m$) = V.

The tube T($U_m$, V) is deformed in an $\varepsilon$- neighbourhood of itself to make h smooth on $U_m$. But in this process, points in the neighbourhood of U$\cap$V would have been moved. Obviously, a tube modification needs to be applied on U$\cap$V to restore the smoothness of h on U$\cap$V. In fact, this will entail a sequence of tube modifications on each of the $U_i$'s in turn, till we again reach V. On reaching V, there are two possibilities :

<u>Case 1 :</u> On finally reaching V we may find that the region R of V that now has to be isotoped does not intersect U∩V. In this case, the final isotopy will move points in a neighbourhood N of R such that N∩U = φ, while no points outside N will be moved.

Figure 4

The sequence of tube modifications thus terminates at this step.

<u>Case 2 :</u> ( U ∩V) ∩R = $R_1 \neq \varphi$.

Fig 5

In this case, after a tube modification on $U_m$ to restore smoothness there, we have to carry out a fresh sequence of tube modifications beginning with $R_1$ as before. As a consequence of the minimality of the flow, it can be seen that this procedure will terminate after a finite number of such repetitions. Hence, after a finite number of such tube modifications, it is clear that the smoothness of h can be restored to the regions which lost it at points which had to be moved in the process of an earlier tube modification. The technique now consists in employing the smoothing procedure described above, over progressively larger regions obtained by adjoining new regions to old ones as described.

<u>Remark 1.1 :</u> In fact, it is obvious that by choosing the regions adjoined suitably, we can get a sequence of regions over which h is smooth and whose union is $D - \{D_1\}$, where $D_1$ is a disc which is as small as we please.

<u>Definition 1.5 :</u> We will refer to the disc $D_0$ as the <u>non-minimal region.</u>

<u>Definition 1.6 :</u> Let v be a flow and D an open disc in general position transverse to v such that :

(i) v is a minimal flow

(ii) There is a disc $D_0$ in the interior of D such that h : $D_0$ --------$\rightarrow$ D is a homeomorphism. Let $T_1$ be the tube with initial face $D_0$ and end face h($D_0$).

(iii) If T is a tube with initial face U and end face V for two discs U and V in D, with T $\cap$ $T_1$ = $\varnothing$, h : U --------$\rightarrow$ V is a smooth diffeomorphism.

The flow $v_1$ arrived at by first carrying out a tube modification n times on v and then "filling in" the tube $T_1$ with line segments whose end points are on the initial and end faces are on the tube $T_1$, to smooth the flow, will be termed a n-derived minimal flow.

It is to be observed that every time we carry out a tube modification, we progressively shrink the region over which the flow is non-minimal, by choosing an arbitrary region inside the non-minimal region and smoothing out the flow inside this arbitrary region. Let $U_n$ be the region we choose in the non-minimal region obtained after applying the tube modification procedure (n-1) times. At the nth application, we smooth out the flow over $U_n$. Evidently, the flow obtained at the end of this application will depend on the regions $U_1$, $U_2$, ......, $U_n$ we have arbitrarily selected for the n modifications of the tube modification procedure. We record the following fact:

**<u>Remark 1.1:</u>**

<u>Another feature of the smoothing procedure is that any minimal set of the flow v must intersect the non-minimal region.</u>

**SECTION 2**

2.1  At the end of the last section we have shown the existence of a disc D which is :

(i) transverse to the flow

(ii) for which the first hit function h : $(D - D_0)$ ------$\rightarrow$ D is smooth, where $D_0$ is a closed disc which represents the non-minimal region. In fact, we will require h to be smooth in a neighbourhood of $\partial D_0$

(iii) the first-hit function h : $D_0$ ----------$\rightarrow$ D is a homeomorphism into the interior of D.

Since, in what follows, we shall be dealing with spaces more general than a disc, we will have to recast the concepts of Section 1 to suit a more general context. To begin with, when we describe a compact 2-manifold X as being in "general position", we will mean :

(i) That it is transverse to the flow

(ii) For any point x ε $\partial X$ and h the first return function h : X → X, one of the following conditions holds :

(a) h(x) is in the interior of X, or

(b) h(x) = y ε $\partial X$ and for two small closed disc neighbourhoods $V_1$, $V_2$ of x and y respectively, with X $\cup$ $V_1$ $\cup$ $V_2$ a smooth manifold transverse to the flow, the arc h($\partial X$ $\cap$ $V_1$) intersects ($\partial X$ $\cap$ $V_2$) in a single point transversally in $V_2$ for h the first return function on [X $\cup$ $V_1$ $\cup$ $V_2$]. This implies, in particular, that there are only a finite number of points $x_i$ on $\partial X$ such that h($x_i$) ε $\partial X$, for h the first return function on X.

Definition 2.1.1 : Let X be a compact orientable 2-manifold in general position and $A_i$, $B_i$ be a finite number of arcs such that :

(i) The $A_i$ are closed arcs lying on the boundary $\partial X$ of X, all of them oriented in the same direction – that is, the orientation being one they inherit from a fixed orientation of $\partial X$ and X. In particular, this means that if arcs $A_p$ and $A_q$ belong to the same component of $\partial X$, then the orientations they inherit from C and D (defined below) make them <u>both</u> positively or negatively oriented with respect to the orientations they inherit from X and $\partial X$.

(ii) The $B_i$ are closed line segments of flow lines, all having the same positive orientation inherited from the flow.

(iii) $C = A_1 \cup B_1 \cup A_2 \cup B_2 \ldots \ldots \cup A_k \cup B_k$ is a simple closed curve.

(iv) C is the boundary of a disc D such that :

(a) The interior of D is a smooth imbedding of an open disc in $S^3$, such that $D \cap X$ is the union of the $A_i$'s.

(b) Whenever the interior of D intersects a flow line, it does so transversally.

We will then call C a singular cycle of X and D a singular disc spanning C.

(Fig. 6 illustrates the case when k = 3).

  Before we embark upon the construction of the foliation in the next section, we will have to construct a space X in general position that does not have a singular cycle. This will be seen to be a necessary condition for the construction of the foliation which we will describe in Section 3. We will develop in this section a method of achieving this object.

<u>Definition 2.1.2 :</u>  The $A_i$ will be called the boundary segments and the $B_i$ the flow segments.

Figure 6.

<u>Definition 2.1.3</u> :  For a 2-manifold X in general position, let P be a point on $\partial X$ such that h(P) = Q $\varepsilon$ $\partial X$, where h is the first return function on X with respect to the flow v. We will refer to P and Q as nodal points of X and the flow line segment with P and Q as end points as a nodal segment of X (Fig. 6).

<u>Definition 2.1.4</u> :  If the flow line segment in the above definition is part of a singular cycle, then it will be called a singular (flow) segment and its end points P and Q, singular nodal points.

<u>Note 2.1.1</u> :  Since X is in general position, there can be only a finite number of nodal points on it.

For a disc D in general position, we now consider a closed disc $D_1$ such that $D_1$ is contained in int D. We will denote $\partial D$ by $b_1$ and $\partial D_1$ by $b_2$.

Figure 7

We now look at [D − (int $D_1$)]. This will be an annulus which we shall refer to as A in the sequel. If A is sufficiently narrow, the following condition will hold :

**<u>Property $P_0$ :</u>**

(i) Let P and Q be on $\partial D_1$ ∪ $\partial D$, with h(P) = Q, where h : D → D is the first return function. In other words, the flow segments with end points P and Q does not intersect interior $D_1$. Then there is a strip S of A containing P and a disc neighbourhood N of Q, with D ∪ N a smooth manifold transverse to the flow, such that h(S) is a homeomorphic image in (D ∪ N) for h : D → (D ∪ N) a first hit function, with A ∩ h(S) in N being as shown in Fig. 8 below.

Figure 8

In other words, the two strips A and h(S) intersect transversally. This is seen to follow from the fact that since D is in general position, there is an arc $\alpha$ with P in its interior such that h($\alpha$) (for h : D $\rightarrow$ D $\cup$ N as before), is a homeomorphic image of $\alpha$ which cuts $\partial$D transversally in a single point in N. We shall refer to this property of the annulus A as <u>Property $P_0$.</u>

(ii) A does not intersect $D_0$ or h($D_0$), where $D_0$ is the non-minimal region defined in Section 1 and h is the first return function on D.

In what follows, we shall assume that $D_1$ has been chosen so that A satisfies this condition. Next, we turn our attention to the flow v. Evidently, the flow v is smooth except inside the tube T with initial and end faces $D_0$ and h($D_0$) respectively. Hence, it can be made smooth by smoothing out the flow lines outside T, to a flow $v_1$. Further, since v is a minimal flow, by choosing $v_1$ to be sufficiently close to v, we can ensure the following properties for $v_1$:

(i) A intersects transversally all flow lines of $v_1$

(ii) Every flow line of $v_1$ which starts at a point of A returns to A

(iii) A is in general position with respect to $v_1$

It has to be remarked that, in general, $v_1$ will not be minimal. In what follows, we shall continue to designate the modified flow $v_1$ by the same symbol v. The reason for this is that we do not have any further necessity for the minimality of v. The three properties of $v_1$ listed above, together with Property $P_0$, will suffice to validate the arguments we shall be using.

<u>Definition 2.1.5:</u>  We will now consider the annulus A as a 2-manifold transverse to the flow. The boundary $\partial$A of A will consist of two copies of $S^1$ which we will denote by $b_1$ and $b_2$, with $b_1$ being the boundary of $\partial$D also. (Figure 7(b)). As in the case of D, we will define nodes and nodal segments of A in a similar fashion. The nodal segments of A can be partitioned into two classes as follows :

(i) Those nodal segments which are also nodal segments of D.

(ii) The nodal segments of A which are not nodal segments of D.

Note 2.1.2 :  It is to be observed that the nodal segments of the first class will have both their end points on $b_1$ and will not intersect $D_1$. The members in the second class will either have one or both end points  on $b_2$ or both end points on $b_1$, in which last case they will necessarily intersect the interior of the disc $D_1$. (Figs. 9 and 10).

Figure 9

The terms like singular node and singular nodal segment for the annulus A  will again be defined in similar fashion as earlier for the disc D.

Definition 2.1.6.  If a singular node or singular nodal segment of A as defined in Definition 2.1.4 is not already a singular node or singular segment of D, we will refer to it as a singular A-node and singular A-segment respectively.

Definition 2.1.7 :  If a singular A-segment of A intersects the underline{interior} of the disc $D_1$, it will be called an irregular singular A-segment and if it does not intersect int $D_1$, it will be called a regular A-segment.  Fig. 10 shows an irregular singular A-segment.

Figure 10

<u>Definition 2.1.8:</u> Let P and Q be nodal points on A with a flow segment oriented in the positive direction with respect to the flow, going from P to Q. We denote this by Q > P. If $s_1$, $s_2$ on $\partial A$ are segments with end points which are nodal points P, Q, with P $\varepsilon$ $s_1$, Q $\varepsilon$ $s_2$ and Q > P, we denote the relationship by $s_2 > s_1$. (A priori, we do not rule out the possibility $s_1 > s_2$ and $s_2 > s_1$). Fig.6 shows a typical singular cycle of A, with three singular boundary segments $A_1$, $A_2$ and $A_3$, with $A_1 > A_2 > A_3 > A_1$.

**2.2 The Construction** : We will now outline a procedure for killing regular singular cycles of A. We will see in Lemma 2.2.3 that by <u>Property $P_0$,</u> every regular singular node and cycle is associated with a singular node and cycle respectively of D; hence there can be only a finite number of such regular singular cycles for A. (Note 2.1.1).

Let $C_1$ be a regular singular cycle of A with three singular boundary and flow segments, the former being $A_1$, $A_2$ and $A_3$ and six singular nodes. Let $\alpha_1$ be a singular flow segment of $C_1$ with end points $p_1$ and $p_2$ and $h(p_2) = p_1$, with h being the first return function on D (Fig. 11).

Figure 11

Let $N_1$ be a tube (Definition 1.2) which is a neighbourhood of $\alpha_1$ in $S^3$ and let $S_1$ and $S_2$ contained in $N_1$, be the two components of A $\cap$ $N_1$, with $p_i$ $\varepsilon$ $S_i$, i = 1,2 and $S_1 > S_2$. We now isotope $S_2$ by an isotopy which does not move any point outside $\partial N_1$, so that in a small smooth disc neighbourhood $N_0$ of $p_2$, transverse to the flow, $S_2$, in its isotoped position $S_{21}$, intersects $S_1$ in a region $R_{12}$ whose boundary is contained in [$\partial S_1$ $\cup$ $\partial S_{21}$]. (In its new position, $A_2$ will be denoted by $A_{21}$). In other words $R_{12} = h(S_2)$, where h, in the tube $N_1$, is the first hit function from $S_2$ to $S_1$ (Fig. 12). This will imply, of course, that the two pieces $S_{21}$ and $S_1$ cross transversely. This is a consequence of  <u>Property P$_0$</u>  formulated earlier.

Figure 12

Figure 12 shows the isotopy being carried out

  We now carry out a similar isotopy inside each of the other tubes $N_2$ and $N_3$ (Fig 11) – that is, in each case, the lower strip is raised along flow lines till it intersects the upper strip. As earlier in the case of $A_{21}$, we shall denote the new segments by $A_{32}$ and $A_{13}$. Upon the completion of these isotopies, we get a copy $C_2$ of $S^1$ in M given by $C_2 = (A_{21} \cup A_{32} \cup A_{13})$ (Fig. 12.1). This can be seen to be the boundary of a singular disc $D_{123}$ in M. In fact, if $X_1$ is the 2-manifold of which $C_2$ is a boundary component, $D_{123}$ can be chosen so that $X_{11} = X_1 \cup D_{123}$ is a smooth 2-manifold which will be transversal to the flow and in general position. We will term $C_2$ <u>a first level quotient cycle</u> and denote it by $\beta_1$.

<u>Note 2. 2. 1:</u>

(i) By choosing the tubes $N_1$, $N_2$ and $N_3$ suitably small, it is obvious that the construction can be carried out in such a way that the isotopy does not move any point of any nodal segment of A except those that are collapsed to a point by it inside $N_i$, i = 1,2,3. A second consequence of such a choice is that there will be no nodes of A inside the tubes, apart from the ones that are being identified by the isotopies.

(ii) It is to be observed that the spaces $X_1$ and $X_{11}$ are both orientable. This follows from the fact that the two operations – identification through isotopy and sealing with a disc – both preserve orientation.

(iii) Because of (i), all regular singular cycles of A and $X_1$ (with the exception of $C_1$ which has been collapsed by the isotopies), will be the same.

We can now consider a regular singular cycle of $X_1$ (if it exists), and by a sequence of collapsing isotopies like the above, replace it by a quotient cycle $\beta_2$ in a space $X_2$.

Thus, by performing this operation successively on each regular singular cycle of the spaces $X_1$, $X_2$, ……, we have described a method by which we get a finite sequence of spaces $X_1$, $X_2$, …., $X_k$, where each $X_{i+1}$ is got as a quotient of the space $X_i$.

To summarise this procedure, we begin with the annulus A = $X_0$ and with a sequence of three isotopies, (as illustrated in Figs. 11 and 12), we get a first level quotient cycle $\beta_1$ and a corresponding quotient space $X_1$. We now choose a regular singular cycle of $X_1$, carry out a sequence of isotopies on it as before, to get a second level quotient cycle $\beta_2$ contained in the space $X_2$. The sequence of spaces $X_1$, $X_2$,….., $X_k$ is thus generated, where $X_i$ is got from $X_{i-1}$ by the formation of an ith level quotient cycle $\beta_i$. In the end space $X_k$, there will be no regular singular cycles which were originally present.

But it will be observed that the procedure just outlined for killing the regular singular cycles may create new singular cycles. An example is given in Fig 12 below, where two sections $B_1$ and $B_2$ of A are as shown below in Fig 12 (a).

Figure 12 (a)

Let $P_1$, $Q_1$ be points on $B_1$, $B_2$ respectively, such that:

(i) There is a flow line segment $L_1$ running from running from $P_1$ to $Q_1$ in the positive direction of the flow. (Fig. 12(b)).

Figure 12(b)

(ii) Under the procedure outlined above, it becomes necessary to raise a small neighbourhood of $P_1$ and bring it into coincidence with a neighbourhood of $Q_1$ (Fig. 12(c)), as has been described above, with the branch $B_1$ being shown in the new position by $B_{11}$.

Figure 12(c)

But this step might give rise to a new singular cycle $\alpha$, shown below in Fig. 12(d), which consists of the union of a boundary segment C of $B_{11}$ and a flow line segment $L_2$ going from $B_{11}$ to $B_2$, denoted in Fig. 12(d) by a bold line and a broken line respectively.

Figure 12 (d)

We shall assume that $\alpha$ is a regular singular cycle. (In fact, from Lemma 2.2.1, we will see there are no other kinds of singular cycles). Accordingly, the limbs $B_{11}$ and $B_2$ will cross each other as in Fig. 12.(a), rather than in the position shown in Fig. 18(d).

Finally, since $\alpha$ is regular singular cycle, the method described above can be used to eliminate it − i.e., the lower limb $B_{11}$ will be moved up by an isotopy inside a tube and attached to the upper limb $B_2$. (Fig. 12(e)).

Figure 12(e)

After a finite number of applications of this procedure, we will reach a 2-manifold $A_1$ with boundary, *without any regular singular cycles.* Some of the components of this boundary will be quotient cycles. Each of these quotient cycles can now be sealed by attaching discs along each of these cycles in such a way that after these attachments, we get a smooth 2-manifold with boundary which is in general position with respect to the flow, without any regular singular or quotient cycles. However, in order to assure ourselves that $A_1$ does not contain any *singular cycles,* we need to show :

(i) We do not introduce any new <u>irregular</u> singular cycles at any stage during the application of the isotopy technique (since the only remaining singular cycles would have to be irregular, after the regular ones have been sealed off).

(ii) There were no irregular singular A-cycles to begin with

Hence, in the concluding part of this section, we shall establish a key lemma which shows that the above construction can be carried out without the introduction of any such cycles. In fact, we will show the impossibility of irregular singular A-cycles in the initial or any subsequent stage of the construction.``

Before we take up the next lemma which establishes this, an explanation of its significance becomes necessary. In an attempt to produce a 2-manifold in general position without singular cycles, we have outlined in the Construction a technique for killing regular singular cycles. This involved the use of isotopies which effected identifications on the annulus A. In the final stage, the holes obtained thus were sealed off by discs transversal to the flow.

An important fact which made the construction possible was that ensured by property $P_0$ and the fact that the nodal segments were regular; this was, that within the cylinders in which the isotopies were carried out, two branches of A crossed in a simple fashion illustrated by Figs. 11 and 12(a). This need not always be the case and, in general, the crossing scheme in the neighbourhood of an *irregular* nodal segment can be extremely complicated. But the following lemma, by showing that irregular nodal segments cannot be parts of singular cycles, obviates the necessity of considering these geometric complications while removing singular cycles.

We set out here a few preliminaries which would be necessary in the proof of lemma 2.2.1.

Definition 2.2.1:  . Let $X_n$ be obtained by carrying out a number of identifications on A effected by a finite sequence of isotopies $I_1$, $I_2$ …………, $I_p$. Let F be the composition of all these isotopies – that is, F = $I_p \circ I_{p-1} \circ \ldots \ldots \circ I_2 \circ I_1$.

Note 2.2.2: We look at a nodal segment $\beta$ of $X_n$. Clearly, there is a nodal segment $\alpha$ of A such that F($\alpha$) = $\beta$, with F as in the above definition. (It is to be observed that from Note 2.2.1, $\alpha$ and $\beta$ are in fact, identical as subsets of $S^3$).

Definition 2.2.2: Let C be a singular cycle of $X_n$, having $\beta$ as a singular flow segment. If $\beta$ = F($\alpha$) (for F as above) where $\alpha$ intersects the disc $D_1$ in its interior, then $\beta$ will be referred to as an irregular singular $X_n$-segment.

Note 2.2.3: From Note 2.2.1, it will be clear that $\beta$ cannot intersect any of the tubes used in the isotopy. Hence, $\alpha$ and $\beta$ are identical as subsets of $R^3$. Therefore, the irregular singular cycle $\beta$ will also intersect the disc $D_1$.

Lemma 2.2.1:  The space $X_n$ will have no singular cycles among whose flow segments there is an irregular singular $X_n$-segment.

Proof:  If not, let there be a singular cycle $\alpha$, with irregular singular $X_n$-segments. Let $D_\alpha$ be a singular disc spanning $\alpha$. (Fig.13(a)).

We begin by considering the space $X_n$. It is to be recalled that $X_n$ is obtained by the repeated application of two procedures on the annulus A :

(i) First, certain identifications are effected on it by isotopies carried out within tubes, a typical one shown in Fig. 13(b) below.

Figure 13

The effect of the isotopy that is carried out is shown in Fig.14, where R is the region at which the two branches of A are identified as already described in the Construction.

Figure 14

(ii) Second, the resulting cycles are sealed off with transversal discs.

We now consider a simple closed curve $\beta$ in the interior of A and homotopic to both bounding circles of A. (Fig.15).

Figure 15

We now claim that $D_\alpha \cap \beta = \Phi$.

In fact, we will show that $D_\alpha$ does not intersect the interior of A. To show this, it is enough to show that $D_\alpha$ does not intersect (A $\cap \cup T_i$), where $T_i$ are the tubes inside which the isotopies are carried out, since outside these tubes $D_\alpha \cap X_n$ must be contained in $\partial A$ by definition. Accordingly, we consider a typical tube $T_1$ which intersects $D_\alpha$ and isotope it so that it is a vertical cylinder (which we will continue to denote by $T_1$) in $R^3$ whose initial and end faces are horizontal and the flow lines are vertical (Fig 14). Fig.14 also shows how $(D_\alpha \cap T_1)$ is located inside $T_1$ and why it cannot, therefore, intersect $\beta$ or the interior of A.

Hence, we can assume $D_\alpha \cap \beta = \Phi$. ……………………………………………………. (A).

We now return to a consideration of the disc D and the annulus A. Let $A_1$ denote the set h(A,D), where h is the first hit function from A to D.

Figure 16

A and $A_1$ partition D into a number of disjoint compact regions $R_i$ which intersect $A \cup A_1$ along their boundaries. (Fig.16 shows D partitioned into six such regions). Let $\alpha_1$ be an irregular flow segment of $\alpha$. Since $\alpha_1$ intersects D but cannot intersect the interior of $A \cup A_1$, $\alpha_1$ must intersect one or more of the $R_i$. Let it intersect the region $R_1$ of Fig.16.1. The point of intersection would be either in the boundary or in the interior (or both) of $R_1$. But it is clear that $\alpha$ can be slightly deformed to $\alpha_2$, so that the following relations hold:

(i) $\alpha_2$ intersects int $R_1$ at a finite number of points

(ii) At each of these points the intersection is transversal

(iii) Orientations in $R_1$ and $\alpha$ (and hence $\alpha_2$), can be chosen such that the intersection number of $\alpha_2$ with $R_1$ at each point of intersection would be +1. ……………………………………. (B)

Figure 16.1

Let $\beta_1$ (Fig.16.2), be a curve in interior $A_1$ such that:

(i) $\beta_1$ is homotopic to the boundary of the region $R_1$

(ii) $\beta_1$ is the union of pieces of $\beta$ which have been carried to $A_1$ under $h: A \rightarrow D$

Hence there is a simple closed curve $\beta_2$ such that $\beta_2$ is the union of:

(i) a finite number of pieces $\beta_{11}, \beta_{12}, \ldots\ldots\beta_{1k}$ of $\beta$ with $h(\beta_{11} \cup \beta_{12} \ldots\ldots\cup \beta_{1k}) = \beta_1$, where $h : A \rightarrow A_1$ is the first hit function.

(ii) segments of flow lines $\gamma_1, \gamma_2, \ldots\ldots\ldots, \gamma_k$ with one end point on $A$ and the other on $A_1$. Since $h: A \rightarrow A_1$ moves each point along flow lines, we can imagine $h$ homotoping $\beta_2$ to $\beta_1$ by collapsing each such flow line segment to a "vertex" like $P_1, P_2, P_3, P_4$ or $P_5$ and carrying each $\beta_{1i}$ into $\beta_1$. ( Fig.16.2) of $\beta_1$. Hence, we can extend $h$ to a map $h_1$ defined on $(A \cup \gamma_1 \cup \gamma_2 \ldots\ldots \cup \gamma_k) \rightarrow A_1$, which is such that $h_1 = h$ on $A$, while $h_1$ additionally collapses each flow line segment like $\gamma_i$ to a vertex.

.

Figure 16.2

It follows that $\beta_2$ is homotoped to $\beta_1$ by $h_1$. From this we get $\beta_1$ and $\beta_2$ are homotopic.

Further, the curve $\beta_2$ can evidently be homotoped to $\beta_1$(by moving each point along the flow line through it) in such a way that in any intermediate position during the homotopy, $\beta_2$ does not intersect the disc $D_\alpha$. Hence, the linking number $Lk(\beta_2, \alpha) = Lk(\beta_1, \alpha)$ ……………… (C)

But from (A), $D_\alpha$ does not intersect $\beta$ and hence $\beta_{11}, \beta_{12}, \ldots, \beta_{1k}$. Nor is it possible for $D_\alpha$ to intersect the flow segments of $\beta_2$. Hence, $D_\alpha \cap \beta_2 = \Phi$. It follows that $Lk(\beta_2,\alpha) = 0$ …….(D)

$$= Lk(\beta_1, \alpha) \text{ from (C)} \ldots\ldots\ldots(E)$$

However, from (B) it is clear that $Lk(\beta_1, \alpha) = Lk(\beta_1, \alpha_2)$ is a positive integer …………..(F)

Since (E) contradicts (F), the impossibility of an irregular singular cycle like α is proved. This establishes the lemma.

Figure 17 below shows the sequence of operations that are carried out, beginning with the disc D, till we reach $D_{12}$.

Figure 17

Note 2. 2. 3: The lemma also shows that there are no irregular singular A-segments.

Lemma 2. 2. 2 : The Construction described in 2.2 will terminate after a finite number of stages.

Proof : In the construction, certain identifications are made on the annulus A during the course of the isotopies described. As a result of these identifications, some singular cycles may be introduced. But by Lemma 2.2.1, none of these can contain an irregular singular $X_n$-segment.

Let us, therefore, look at the kind of singular cycles α that may be formed – in particular, what happens at the flow segments of such singular cycles. These are indicated in Fig. 18 below. The diagram shows that in a small neighbourhood of a nodal segment $α_1$ of $X_n$, <u>two branches of A cross each other.</u> They cannot, for instance, be as in Fig. 18(d).

Figure 18

In every other case above, in Fig. 18, it will be seen that, in a suitable neighbourhood of α, A will be seen to cross itself transversally as already displayed in Fig. 8. This is a consequence of <u>Property $P_0$</u> enunciated earlier. This will imply that in a neighbourhood of a (regular) singular $X_n$-flow segment, there will be a corresponding nodal segment $α_2$ of D which will connect two points of $b_1$. Underline{Hence, with every singular $X_n$-flow segment, there will be associated a nodal segment of D.} Also, since the disc D had only a finite number of nodes to begin with, the number of crossings like the ones in Fig. 11, and hence the number of cycles, cannot be greater than the number of nodes and therefore must be finite. Further, with each step of the construction, the number of nodal segments keeps getting reduced. This, in turn, implies the lemma.

Note 2.2.4 : It is possible that $\underline{X_n \text{ does not have any boundary.}}$ In other words, we find that when the process of sealing off the holes is completed, there are no free edges and hence what is obtained at the end of this process is a 2-manifold without boundary. This, of course, cannot happen when the space M on which the minimal flow is defined is $S^3$, since, in this case, $X_n$ would separate $S^3$. But it is certainly possible in a more general M. The process of constructing a foliation transverse to the flow, which will be demonstrated in the next Section, will become particularly simple in this special case. (In fact, the presence of such a foliation can also be established independently by a more direct method in this case). We shall, therefore, assume in the sequel, that $X_n$ has a boundary.

**The topological nature of $X_n$**

Let $X_n$ be the 2-manifold resulting from successive applications of the Construction described in 2.2. The following geometric features of $X_n$ are set down for the record :

(i) $X_n$ does not have any singular cycles.

(ii) $X_n$ is in general position with respect to a flow $v_0$ which is such that the flow line of $v_0$ through any point on $X_n$ returns to $X_n$.

We shall sum up the above observations in the following theorem.

<u>Theorem A :</u> $X_n$ is a compact orientable 2-manifold with boundary, which is in general position with respect to the flow and which is without any singular cycles.

<br>

## <u>SECTION 3</u>

## <u>The Foliation</u>

<br>

In this section, we use the flow v to construct a co-dimension one foliation on $S^3$ transverse to the flow.

We shall begin this Section by considering the space $X_n$ of Theorem A in Section 2 and shall refer to it as X in the sequel.

As in the earlier instance with the disc D in Section 1, we recall that X can, through a finite number of modifications carried out in an arbitrarily small neighbourhood of its boundary, be made to satisfy the following general position property :

We will consider a partition of X obtained as follows. Let Y be the system of arcs given by Y = h[$\partial$X, X], where h is the first hit function with respect to $-v$. By the process used in Section 1, X can be perturbed in an arbitrarily small neighbourhood of its boundary so that :

(a) Y consists of a finite number of smooth arcs such that whenever two arcs intersect each other, they do so transversally

(b) No interior point of X has more than two arcs meeting there and no point of $\partial$X has more than one arc of Y intersecting $\partial$X there

(c) Whenever an arc of Y intersects $\partial$X it does so transversally

(d) By introducing a finite number of arcs, the system of arcs Y can be augmented to a system $Y_1$ such that $Y \subset Y_1$ and $(X - Y_1)$ consists of a finite number of regions the closure of each of which is topologically a disc.

In what follows, we will therefore assume that Y itself has the properties (a) – (d) above.

<u>Definition 3.1:</u> We shall call the space X together with the system of arcs Y and the regions (X-Y) the **Partition Space** corresponding to the partition $\pi$ or resulting from the partition $\pi$ and denote the connected regions of (X - Y) by $P_1$, $P_2$, ……. $P_n$. The $P_i$ will be referred to as **pieces**. $\pi$ itself will be used to denote either the partition or the Partition Space obtained therefrom, where this will not lead to any confusion.

It has to be remembered that, in general, X need not necessarily be homeomorphic to a subset of $R^2$, as it is represented in Fig.19.

Figure 19

<u>Note 3.1:</u> Let $R_i$ be the closure of the pieces $P_i$. We shall refer to the $R_i$ as the regions of X. It is to be observed that, the first hit function $h(R_i, D)$ with respect to v will not, in general, be a homeomorphism, unlike $h(P_i, D)$. But by convention, whenever we refer to $h(R_i, D)$, we shall assume that it is a homeomorphism, by extending h in an obvious fashion to the closures $R_i$ of the $P_i$.

At the heart of the proof is a technique for constructing a co-dimension one foliation on $S^3$ which is transversal to v. We do this by progressively foliating certain closed subspaces of $S^3$ which we shall call tubes. Each tube will be :

(i) Topologically $D^2 \times I$, where D is the closed two dimensional disc and I is the interval [0,1].

(ii) Made up of a union of closed flow line segments whose end points are on X.

(iii) No point in X is in the interior of any tube.

**Convention** :  For two 2-manifolds $S_1$, $S_2$ transverse to the flow v, $h^1(S_1, S_2)$ will denote the first-hit function from $S_2$ to $S_1$ with respect to $-v$.

The foliation that we will undertake progressively, on one tube after another, will be referred to as the **Foliation Procedure** in what follows and is described below. We will use the term **leaf** to denote any region or union of regions of $\pi$ or any leaf in the foliation at a later stage.

**The Foliation Procedure:**  We begin with the regions of $\pi$. Given a region U of $\pi$, we consider the tube with initial face U and end face V = h(U), for the first return function h on X. From note 3.1, we assume h/U is a homeomorphism. (In the sequel, a tube will always be understood to mean one with its initial and end faces on X). We now describe a procedure for foliating these tubes with a co-dimension one foliation, one after another in turn, subject to the following conditions :

(i) Each leaf of the foliation of a tube will be a two dimensional disc with its boundary on the boundary of the tube, with U and V themselves being leaves.

(ii) Whenever two adjacent tubes $T_1$ and $T_2$ – i.e, two tubes which intersect – are foliated, the union of the leaves gives a smooth foliation for $T_1 \cup T_2$. In particular, whenever a leaf L intersects a region R of $\pi$ , $L \cup R$ is a smooth 2-manifold

(iii) The leaves will be transverse to all the flow lines they intersect.

Figure 20

<u>Definition 3.2</u> : Let us suppose we have foliated successively n tubes $T_j$, j = 1,2,……,n. Then the space $X \cup [\cup T_j]$ is called the nth tube space and will be denoted by $S_n$. By convention, X, which is the union of the regions $R_i$, will be called the zeroth tube space and will be denoted by $S_0$. Hence we have the inclusion relations $S_0 \subset S_1 \subset \dots \dots \subset S_n$.

<u>Definition 3.3</u> : For a tube T with initial and end faces $F_1$ and $F_2$, closure $[\partial T - (F_1 \cup F_2)]$ will be called the lateral surface of T. The lateral surface of $S_n$ will be the union of those parts which are not in the interior of $S_n$, of the lateral surfaces of the tubes or the $R_i$ whose union is $S_n$. In the sequel, the lateral surface of a tube T or a tube space $S_n$ will be denoted by $\partial T$ or $\partial S_n$ respectively.

<u>Definition 3.4</u> : By a ruling on a tube space $S_n$, we will mean the one- dimensional foliation of any part of the lateral surface of $S_n$ by curves which are the union of the curves $\partial X$ and the boundaries of the leaves foliating $S_n$. Any connected subarc of one of the curves constituting a ruling R on a tube space $S_n$ will be called a **<u>curve of the ruling.</u>**

<u>Definition 3.5</u> : Let R be a ruling on a tube space $S_n$. Let $A_1$, $B_1$, $A_2$, $B_2$,……,$A_k$, $B_k$ be line segments such that :

(a) $A_1$, $A_2$, ……,$A_k$ are closed line segments of the curves of the ruling R which are coherently oriented. In other words, if $A_i$, $A_j$ are subarcs of the same curve α of the ruling, then, with respect to the orientation they inherit from the orientation of the closed curve C (defined in (c) below), both are either positively or negatively oriented with respect to a fixed orientation on α.

(b) $B_1$, $B_2$, ……..,$B_k$ are closed, directed (in the positive direction of the flow) line segments of flow lines.

(c) The union of $A_i$'s and $B_i$'s forms a simple closed curve C for which $A_1$, $B_1$, $A_2$, $B_2$, ..........$A_k$, $B_k$, $A_1$ are adjacent segments having an end point in common. ( It is to be noted that some of the $B_i$'s may degenerate into single points). C will be referred to as a **general cycle** of the ruling.

Definition 3.6 : Let C be a general cycle of a ruling R on a tube space $S_n$ such that at least one of the $B_i$'s does not degenerate into a point. Let D be a disc such that:

(i) $\partial D = C$.

(ii) (D - $\partial D$) is a smooth imbedding of an open disc into $S^3$.

(iii) (D - $\partial D$) is transverse to all flow lines it meets.

(iv) If a leaf L intersects int D, then $L \subset$ int D.

In this case, C will be called a singular cycle and D a singular disc and the ruling R will be said to contain a singular disc on $S_n$. D will then be said to **span** C.

Note 3.2:  (i) There is an equivalent definition of a singular disc which is sometimes easier to use. Here, instead of requiring C to be the boundary of D as in condition (i) above, we will require the existence of a closed curve $C^0$ which is the boundary of a disc D satisfying (ii), (iii) and (iv) and which can be carried along the flow lines – either backward and/or forward – to coincide with C. [Refer Principal [P] under Remark 3.16].

(ii) When we speak of the ruling containing a singular disc, the terminology may suggest that the singular cycle which is the boundary of the disc will be contained in $\partial S_n$. But it is to be observed that only the segments $A_i$ of the boundary of $\partial D$ are necessarily in $\partial S_n$. In particular, some of the $B_i$'s may be partly outside $\partial S_n$.

(iii) It will be observed that a singular cycle cannot intersect the initial or end face of a tube, since this would imply that the cycle has positive linking number with $\partial X$.

Remark 3.1 : All the singular discs like D above, that will be considered in the sequel, will be assumed to either not intersect int $S_n$, or if any of them do, the intersection will be assumed to satisfy condition (iv) of definition 3.6 above.

Definition 3.7 : If a ruling on a tube space $S_n$ has no singular cycles, the ruling as well as the space $S_n$ will be called **acyclic.** (Hence, the ruling on the tube space $S_0$ – i.e, $\partial X$ - will be acyclic by Theorem A of Section 2. This ruling on $S_0$ we shall denote by $\Omega$).

Definition 3.8 : Beginning with the space $S_0$, we "fill in" the tubes one after another with discs to get a foliation of $S_n$ by 2-manifolds. But this foliation should be carried out such that at each successive stage we get a tube space for which the ruling associated with the foliation is an acyclic ruling. This condition, in what follows, will be referred to as the **Constraint Condition**.

Definition 3.9 : Given a tube space $S_n$, a disc D with $\partial D$ made up of a union of sub-intervals of flow lines and ruling segments $B_i$ and $A_i$ respectively, is said to be transversal for $S_n$ – or just transversal - if D is as in 3.6 – i.e, it satisfies conditions (ii), (iii) and (iv). (We allow the possibility of all the $B_i$'s being points).

The following alternative formulation of the Constraint Condition is useful in clarifying the geometrical implications of the presence of a singular disc.

**The geometric meaning of the constraint condition**. Let $S_n$ be a tube space with a ruling $R_n$ on it. Let T be a tube not in $S_n$ and $S_{n+1}$ the space $S_n \cup T$. Let D be a disc such that :

(i) D is transversal to the flow

(ii) $\partial D = C \cup C_1$, where C consists of arcs $A_1$, $B_1$, $A_2$, $B_2$,.......$A_k$, $B_k$ with the $A_i$'s being segments of curves of the ruling $R_n$ and the $B_i$'s are directed (in the positive direction of the flow) segments of the flow lines not all of which degenerate to points and $C = A_1 \cup B_1 \cup A_2$ .......$\cup A_k \cup B_k$. $C_1$ is a curve on $\partial T$ whose interior is in a region which does not contain any rulings of $S_n$, with C and $C_1$ intersecting at a pair of common end points $p_1$ and $p_2$, which are one end point of respectively $A_1$ and $B_k$.

Then when T is being foliated to give a ruling $R_{n+1}$ on $S_{n+1}$, the constraint imposes the condition that if the curve L on $\partial T$ is the boundary of one of the discs foliating T, and L passes through $p_1$, there is a point p on L which is below $p_2$, in the sense that there is a segment of a flow line running from p to $p_2$ in the positive direction.

In the special case when C is a curve of the ruling, i.e, when all the $B_i$'s degenerate to points, the constraint implies that when T is foliated, some disc of the foliation on T should have a boundary with a sub-arc $C_2$ such that $C_2 \cup C$ is a closed curve on $\partial S_{n+1}$, with $C_2 \cup C$ being the boundary of the disc.

<u>Definition 3.10 :</u> A connected curve of $S_n$ without self-intersections which is made up of a union of a finite number of segments of flow lines and curves of the ruling will be known as a **<u>path.</u>** A simple closed path is similarly defined.

<u>Definition 3.11 :</u> Let C be a path of $S_n$ which is made up of $A_1$, $B_1$, ………. ,$A_k$, $B_k$, with the $A_i$'s being segments of the curves of the ruling and the $B_i$'s being segments of the flow lines not all of which degenerate into points, as before. If any orientation on C induces on the $B_i$'s orientations which are all positive or all negative with respect to the standard orientation on the flow lines, C will be called an **<u>ascending or descending path</u>** respectively.

**<u>The Construction of the Foliation</u> :** Let $S_n$ be an acyclic tube space obtained at the end of the nth stage of the foliation procedure, with a (co-dimension one) foliation $F_n$ defined on it. We now give a construction for extending $F_n$ to a foliation $F_{n+1}$ on $S_n \cup T$, where T is a tube not contained in $S_n$, consistent with the Constraint Condition. The task obviously becomes trivial in the case when $T_i \cap T = \varnothing$ (or even when the lateral surfaces of $T_i$ and T do not intersect at an interior point of $T_i \cup T$), for every tube $T_i \subset S_n$. Hence, let us assume that $T_i \cap T \neq \varnothing$ for some $T_i \subset S_n$ – i.e, their lateral surfaces do intersect at an interior point of $T \cup T_1$. In this case, the ruling $R_n$ on $\partial S_n$ associated with $F_n$ induces a ruling $R$ on $T$. From the fact that $S_n$ is acyclic, it follows that the curves of the ruling $R$ are either closed line segments or circles. There cannot, for example, be a curve spiralling down to a circle. Let the components of the regions of $\partial T$ which do not have any rulings be called the **Unfoliated Regions** . Each such region will therefore be an open set of $\partial T$ bounded by segments of flow lines and of the curves of the ruling $R$. The foliation of $T$ by discs is carried out by first foliating $\partial T$ by circles and then extending this to a foliation of $T$ by discs. The foliation of $\partial T$ is in turn carried out by first extending the ruling $R$ to a foliation of each unfoliated region by curves, one after another.

To sum up, the procedure we shall be using in the sequel is essentially an inductive one and will consist of two parts :

(i) An explicit assumption that the tube space $S_n$ obtained at the nth stage is acyclic

(ii) Under the above assumption, we shall describe a procedure for extending the foliation on $S_n$ to $S_{n+1} = S_n \cup T$, where $T$ is a tube not contained in $S_n$, such that $S_{n+1}$ is acyclic

<u>Definition 3.12 :</u>  Let $P_1$ be a curve of the ruling on $\partial S_n$ and $C_1$ an arc in a connected unfoliated region $U_1$ of $\partial T$ such that :

(i)  $C_1 \cap \partial U_1 =$ the end points $\{x_1, x_2\}$ of $C_1$.

(ii) $P_1 \cap \partial U_1 = \{x_1, x_2\}$, the end points of $P_1$.

(iii) $C = C_1 \cup P_1$ is a simple closed curve such that there exists a disc $D$ whose interior cuts all flow lines it meets transversally and $\partial D = C_1 \cup P_1$. (This implies that $C_1$ cuts transversally, all flow lines it meets).

(iv) If D intersects interior $S_n$, then every component of $D \cap S_n$ is a leaf of $S_n$ which is in the interior of D.

In this case, we refer to either $C_1$ or $P_1$ as a C-curve for $U_1$ and D as a C-disc for $C_1$ or $P_1$.

Figure 21

As a first step in constructing a foliation on a tube T, we choose an unfoliated region $U_1$ of $\partial T$ and extend the ruling R on $\partial T$ to a foliation $R_1$ of $U_1$. We begin this process by first connecting by an arc $C_i$ in $U_1$, the end points of a curve $\Lambda_i$, whenever $\Lambda_i$ is a C-curve. Such an extension is possible, of course, only when the curves $C_i$ can be drawn without cutting each other. It will be the object of lemma 3.1 to show this. Before that we introduce some definitions.

<u>Definition 3.13:</u> Let $S_n$ be a tube space obtained at the nth stage and T a tube not contained in $S_n$. As already observed, if the lateral surfaces of $\partial S_n$ and $\partial T$ intersect, a subset of $\partial T$ will be foliated by the rulings on $S_n$. This part of $\partial T$ will be called the foliated region of $\partial T$ and each connected subset of the closure of the complement in $\partial T$ of the foliated regions will be called an unfoliated region (of $\partial T$). It is to be noted that the unfoliated regions may contain the rulings contributed by the boundaries of the regions $R_i$ – i.e. the zeroth tube space $S_0$ – with which we began the foliation procedure. These may not be in the boundary of $\partial T$.

<u>Definition 3.13.1 :</u>  The boundary of an unfoliated region $U_m$ of $\partial T$ is made up of a number of closed line segments which are sub-intervals of flow lines and curves of the ruling arranged alternately, with two adjacent sub-intervals sharing a common end point. These will be referred to as **<u>flow segments</u>** and **<u>ruling segments</u>** respectively. (Figure 21.1).

Figure 21.1

<u>Definition 3.13.2 :</u>  Let I be a ruling segment of $\partial U_m$ with end points $x_1$, $y_1$. We term I a **free segment** if both $x_1$ and $y_1$ are such that any sub-intervals $I_1$, $I_2$ of flow lines in whose interiors $x_1$ and $y_1$ respectively lie, intersect the interior of $U_m$. (I and $I_0$ are free segments in Fig.21.2(a)).

<u>Definition 3.13.3 :</u>  Let I be a free segment of $\partial U_m$, where $U_m$ is an unfoliated region of $\partial T$. By a subdivision of $U_m$, we mean the process of joining a point $x_1$ of I to a point $x_2$ of $\partial U_m$ by a segment of a flow line. (Fig. 21.2 (b)).

<u>Note 3.13.3 :</u>  (i) In the above figure, we see $U_m$ being subdivided into three regions $U_{m1}$, $U_{m2}$ and $U_{m3}$ by two segments of flow lines with end points $x_1$, $x_2$ and $y_1$, $y_2$.

(ii) For convenience, we have represented $U_m$ as a homeomorphic copy of a disc. But it could possibly be an annulus, or an annulus with holes, without affecting any of the arguments used here.

<u>Definition 3.13.4 :</u> An unfoliated region $U_m$ of $\partial T$ will be called **Ruling Convex** if $\partial U_m$ does not have any free segments. Figures 21.2 (a) and (b) give examples of unfoliated regions which are not Ruling Convex.

Figure  21.2

It is to be observed that all the three regions obtained by the subdivision in Fig.21.2 (b) (where the subdivisions are effected by the introduction of two flow line segments), are Ruling Convex.

In fact, it is not difficult to see that by a suitable subdivision by flow line segments, <u>any</u> unfoliated region of $\partial T$ can be partitioned into regions each of which is Ruling Convex.

<u>Note 3 .13 .4 :</u>  It will be our objective to extend the curves of the ruling on $\partial T$ to a foliation on each unfoliated region of $\partial T$, one after another, so as to get a foliation on $\partial T$ finally. But the technique that we will be employing will require that the region over which the foliation is sought to be extended, should be Ruling Convex. We will, therefore, in the light of the above observation, assume that the unfoliated region under consideration is Ruling Convex.

<u>Definition 3. 14 :</u>  Let $U_1$ be an unfoliated region of $\partial T$. $\partial U_1$ will be made up of a finite number of ruling segments and flow segments. A curve in the former class will be called a **top segment** if no point of it can be connected to some point in $\partial U_1$ by an ascending curve $\partial T$. A **bottom segment** is similarly defined.

<u>Note 3.14.1 :</u>  If an unfoliated region $U_1$ is Ruling Convex, it has the following properties :

(i) It has a unique top segment and a bottom segment.

(ii) There are just two ascending curves on $\partial U_1$ connecting the top and bottom segments. These will be called the **lateral curves**.

(iii) If two points are chosen on one of these lateral curves such that the points are not on the same ruling segment of $\partial U_1$, they can be connected by an ascending path. In fact, this ascending path can be chosen so that every flow line on $\partial T$ intersects it in a connected set.

<u>Definition 3.15 :</u> Let A, B be two discs transverse to the flow such that $\partial A = C_1 \cup P_1$ and $\partial B = C_2 \cup P_2$, where the $P_i$'s belong to the boundaries of two leaves $L_i$ of $S_n$, while the $C_i$'s are two curves of an unfoliated region $U_n$ of $\partial T$. (It has to be observed that the $P_i$'s may intersect $\partial T$ only at their end points on $U_n$). Then A and B will be said to be **boundary intersecting discs** if the following condition holds : any arc whose end points are those of $C_1$ and whose interior lies in int $U_n$, will intersect any arc whose end points are those of $C_2$ and whose interior is in the interior of $U_n$. In such a case, it can be shown that one end point of $C_1$, will be above that of $C_2$, while the other end point will be below the other end point of $C_2$. (Here, a point P will be above or below another point Q if there is an ascending or descending path respectively, running from Q to P). The reason for terming A and B to be boundary intersecting discs becomes apparent when we note that the above condition implies that $C_1$ and $C_2$ have to intersect somewhere on $U_n$. (Figs. 22 and 23).

Figure 22

Definition 3.16 : If A and B are two boundary intersecting discs which are on opposite sides of an unfoliated region $U_n$, as in Fig 22, they are called **boundary intersecting discs of the first kind** and if they are on the same side of $U_n$ they will be called **boundary intersecting discs of the second kind.**

Remark 3.16 : Two boundary intersecting discs A and B obviously cannot both be contained in T.

We now enunciate, without proof, a principle which we shall have occasion to use repeatedly.

**Principle [P] :**  Let $D_1$ be a transversal disc in $\partial S_n$. Let $C_1$, $C_2$ be closed piecewise smooth curves which intersect every flow line they meet in a connected set – ie., a point or a closed line segment - such that:

(i) $C_1 = \partial D_1$, and

(ii) Every sub-interval of a flow line through a point of $C_1$ which meets it in a connected set also meets $C_2$ in a connected set and $C_2$ can be isotoped to $C_1$ by moving every point of it along the flow line through that point.

Then there is a transversal disc $D_2$ such that $\partial D_2 = C_2$.

<u>Note :</u> In applications of this principle, $C_1$, $C_2$ would be made up of segments of curves of the ruling and flow lines.

The following lemma shows a fundamental property of C-curves.

<u>Lemma 3.1 :</u> Let $U_1$ be an unfoliated region on $\partial T$ of a tube T. Let $C_1$ and $C_2$ be two curves in $U_1$ such that :

(i) The interiors of $C_1$ and $C_2$ are in int $U_1$ while their end points are in the boundary of $U_1$

(ii) There exist two discs $D_1$, $D_2$ transverse to the flow such that $\partial D_1 = C_1 \cup P_1$ and $\partial D_2 = C_2 \cup P_2$ where $P_1$ and $P_2$ are curves on $\partial S_n$ lying on the boundary of two leaves $L_1$ and $L_2$ on $S_n$. (In other words, $C_i$ and $P_i$, i = 1,2 are C – curves).

Then there exist two curves $Q_1$ and $Q_2$ on $U_1$, with their interiors in the interior of $U_1$ and their end points on $\partial U_1$, such that $Q_1 \cap Q_2 = \varnothing$ and $P_1 \cup Q_1$ and $P_2 \cup Q_2$ are the boundaries of two transversal discs on $\partial T \cup \partial S_n$. (Hence, $Q_1$, $Q_2$ must have the same end points as $C_1$,$C_2$).

<u>Proof :</u> If two such curves as $Q_1$ and $Q_2$ cannot be drawn without intersecting each other, it must happen that $Q_1$ and $Q_2$ can be drawn such that they intersect each other at a single point, while $Q_1$ and $Q_2$ cut all flow lines they meet, transversally. This implies that every point of $C_1$ can be isotoped along flow lines to coincide with $Q_1$. A similar remark holds for $C_2$ and $Q_2$.

By applying Principle P to the curves $[C_1 \cup P_1$ and $Q_1 \cup P_1$ ] and $[C_2 \cup P_2$ and $Q_2 \cup P_2$ ], we see that $(Q_1 \cup P_1)$ and $(Q_2 \cup P_2)$ are the boundaries of transversal discs $D_{11}$ and $D_{22}$ respectively.

Figure 23

Since $D_{11}$ and $D_{22}$ are the transversal discs on whose boundaries $Q_1$ and $Q_2$ lie, $D_{11}$ and $D_{22}$ are obviously boundary intersecting discs. Two cases now arise with respect to the discs $D_{11}$ and $D_{22}$.

<u>Case (i)</u> : $D_{11}$ and $D_{22}$ lie on opposite sides of $U_1$ – i.e, $D_{11}$ and $D_{22}$ are boundary intersecting discs of the first kind. From the fact that $U_1$ can be assumed to be Ruling Convex by Note 3.13.4, it follows that there are paths $s_1$ and $s_2$ connecting the end points of $C_1$ and $C_2$ such that one of them, say $s_1$, is a descending curve and the other, $s_2$, is an ascending curve (Note 3.14.1, (iii)).

This implies that the curve $C_0$ = ( $P_1 \cup s_1 \cup$ - $P_2 \cup s_2$ ), where – $P_2$ is oppositely oriented to $P_2$, can be oriented so that so that it is now an ascending or descending singular path. (In the analogous geometric situation in Fig 22, this is shown to be a descending path). <u>Obviously, $s_1$ and $s_2$ can be chosen so that any flow line segment on $\partial T$ which either of them intersects, meets them in a connected set.</u>

We have to show that $C_0$ is the boundary of a singular disc. To do this, we first confine our attention to the part of $D_{11}$ and $D_{22}$ which lies in a small neighbourhood U of $U_1$ in $S^3$. (Fig 23.1).

Fig. 23.1

The four curves $s_1$, $Q_1$, $s_2$, $Q_2$ lie on $\partial T$ and $s_1 \cup Q_2 \cup s_2 \cup Q_1$ forms a closed curve which intersects itself at the point O. We now carry out the following transformations in U. We first isotope U such that :

(i) The flow lines are now vertical

(ii) $U_1$ in its new position is in a vertical plane.

We shall, however, continue to use the same symbols $Q_1$, $s_1$, $Q_2$, $s_2$ for the isotoped positions of these curves. We shall now imagine that $Q_1$ and $Q_2$ are pulled away from each other slightly, so that in their new positions they are respectively behind and in front of the vertical plane containing $U_1$, while being contained in U. Their relevant geometric properties are to be preserved under such a deformation : viz., the curves are smooth and intersect all flow lines they meet transversally. Similar transformations are also carried out on the discs $D_{11}$ and $D_{22}$, so that the parts of the disc $D_{11}$ and $D_{22}$ in U are now deformed slightly while remaining transversal to the flow and the curves $Q_1$ and $Q_2$ are in U, but on different sides of $U_1$ – in other words, in their isotoped positions, the interiors of $D_{11}$ and $D_{22}$ will not intersect int $U_1$. Let these curves in their new positions be denoted by $Q_{10}$ and $Q_{20}$. (Fig. 23.2)

Figure 23.2

If we look at the curve $\alpha = s_1 \cup Q_{20} \cup s_2 \cup Q_{10}$, we see that:

(i) $\alpha$ is a simple closed curve

(ii) $\alpha$ intersects $\partial T$ in $s_1$ and $s_2$, while the interiors of $Q_{10}$ and $Q_{20}$ are on opposite sides of $\partial T$.

(iii) Further, it can be seen that the deformation can be carried out so that every flow line that meets $\alpha$, intersects it in a connected set.

From the above observation, it follows that there exists a disc $D_3$ such that:

(i) $\alpha = \partial D_3$.

(ii) The interior of $D_3$ intersects flow lines transversally

(iii) If $D_{10}$ and $D_{20}$ are the discs into which $D_{11}$ and $D_{22}$ have been isotoped, then $D_{10} \cup D_3 \cup D_{20}$ is a smoothly imbedded transversal disc in $\partial S_n$ with boundary $C_0$, which is also singular there, since its boundary $s_1 \cup$ - $P_2 \cup s_2 \cup P_1$ can be oriented so that it is an ascending or descending curve.

This implies that $\partial S_n$ is not acyclic and this contradiction establishes the impossibility of this case.

<u>Case (ii)</u> : Let $D_1$ and $D_2$ be boundary intersecting C-discs of the second kind for $Q_1$ and $Q_2$ respectively (Fig.23).

Let $Q_1$, $Q_2$ intersect at X. Since $D_1$, $D_2$ intersect all flow lines they meet transversally and are on the same side of $U_1$, in a neighbourhood of X $D_1$ and $D_2$ must intersect transversally. Since $P_1$ and $P_2$ cannot intersect without coinciding, this implies that either $\partial D_1$ intersects the interior of $D_2$ or $\partial D_2$ intersects the interior of $D_1$. Let the first possibility hold – i.e, let $\partial D_1$ intersect int $D_2$ at the point P. Evidently, P can be chosen so that $\partial D_1$ intersects int $D_2$ transversally. This will, in turn, imply that $D_2$ intersects int $S_n$. Since $P_1$ and $P_2$ are C – curves, by (iv) of Definition 3.12, there is a neighbourhood N of P in $D_2$ such that N $\cap$ $\partial S_n$ is an arc $c_1$ of the curve of the ruling $R_n$. Let $c_2$ be a small arc neighbourhood of P in $\partial D_1$. Since $c_1$ and $c_2$, being curves of the ruling, cannot intersect transversally at P, they have to coincide in an arc neighbourhood of P. This can be seen to be clearly impossible.

Hence the Lemma.

<u>Definition 3. 17 :</u>  If $U_1$ is an unfoliated region of $\partial T$, we consider the quotient $X_1$ of $\partial U_1$ obtained by collapsing each ruling segment on $\partial U_1$ to a point. We denote by q the quotient map from $\partial U_1$ to $X_1$.

Since q is obviously a homeomorphism when restricted to any flow segment in $\partial U_1$, it is possible to carry the orientation on the flow lines to an orientation on the line segments of $X_1$. This orientation is displayed in Fig 24.

Figure 24

<u>Definition 3. 18 :</u>  In Fig 24, a point like H which cannot be connected to any other point in $X_1$ by a positively oriented segment, will be called the **<u>high point</u>** of $X_1$. Similarly, the point L will be called the **<u>low point</u>**. It is to be observed that the high and low points are the images under q of the top and bottom segments in $\partial U_1$.

<u>Definition 3.19 :</u>  Let x be a point in $\partial U_1$. x will be called a **<u>node</u>** if there is a sub-interval of a curve of the ruling without interior points in $\partial U_1$, joining x to a point y ($x \neq y$) which is either :

(i) On $\partial U_i$ where $U_i$ (i may be equal to 1) is an unfoliated region of $\partial T$ and y is a point at which a flow segment and a ruling segment of $\partial U_i$ meet, or

(ii) There is a sub-interval I of a flow line which is such that $y \in$ int I and $I \subset \partial S_n$, while y does not have a neighbourhood in $\partial S_n$ made up of a union of sub-intervals of flow lines. In this case, y will be on $\partial X$ and on the intersection of the boundaries of two initial faces or two end faces of two tubes. y will then be a node according to Definition 2.1.3. (See note below). $x_i$, $y_i$ denote the nodes in the Figure 25.

Figure 25

**Note:** The term "Node", as defined here, has nothing to do with terms like "Node", "Nodal segment", etc., occurring in Definitions 2.1.3, 2.1.5 and the discussion immediately following these definitions. In particular, "Node" will have, in the sequel, the meaning given to it in Definition 3.19.

<u>Definition 3.19.1</u> : If $x_1$ is a node of $U_1$, then $q(x_1)$ will be called a node of $X_1$.

<u>Definition 3. 20</u> : If, for two points x and y in $X_1$, there is a curve in $X_1$ running from x to y in the positive direction with respect to the orientation of $X_1$, we say that there is an ascending curve joining x and y. (In Fig 24, there is an ascending curve joining L and B, or C and H ). In this case, we say that L is lower than B, or H is higher than C.

<u>Definition 3. 21</u> : Let $x_m$ be a node of $X_1$. Let p be a node of $\partial U_1$ such that $q(p) = x_m$. $A(x_m)$ ( respectively $B(x_m)$ ) will be the set in $X_1$ defined as follows : For a node y in $X_1$ with $q(r) = y$, y ε $A(x_m)$ $(B(x_m))$ if there is an ascending (descending) path α in $\partial S_n$ running from p to r such that:

(i) There is a curve C in $U_1$ with end points (p,r) in $U_1$, with int C $\subset$ int $U_1$

(ii) α $\cup$ C is the boundary of a disc D which cuts transversally all flow lines it meets.

<u>Remark 3. 3</u> : It can be seen that there will be precisely two curves in $X_1$ (which are the images under q of the lateral curves defined earlier in Note 3.14.1) which have the high and low points of $X_1$ as end points and which we will, by analogy, call the **lateral curves** of $X_1$. If $x_i$ is a node on one of these curves, $A(x_i)$ will be on the other. Hence any two points of $A(x_i)$ can be connected by an ascending arc.

<u>Defininition 3. 22 :</u> For a node $x_i$, $A(x_i)$ will be a finite set. The lowest point of this set will be denoted by $a_i$. Similarly, $b_i$ will denote the highest point of $B(x_i)$. If x is a node which is one end of a C-curve the other end of which is y, we will, by convention, declare $A(x)$ and $B(x)$ both to be the point y.

<u>**Note 3.22 :**</u>

(i) For a node $x_i$, it can be seen that one or both of $a_i$ and $b_i$ can fail to exist. If both exist they can be connected by an ascending or descending path.

(ii) For a node $x_1$, $a_1$ cannot coincide with the lowest point of $X_1$, since this would imply that $\partial S_n$ is not acyclic. Similarly, $b_1$ cannot coincide with the highest point.

(iii) Let $x_1$, $y_1$ ε $\partial U_1$ be such that :

(a) $x_1$ is a node

(b) $q(x_1)$, $q(y_1)$ belong to two different lateral curves of $X_1$ (Fig. 26)

(c) $x_1$, $y_1$ are connected by a path in $\partial S_n$ which ascends in going from $x_1$ to $y_1$.

Then the $a_1$ corresponding to $x_1$ will either coincide with $y_1$ or be below it.

Figure 26

<u>Convention :</u> In the sequel, for two points $x_1$, $x_2$ on $\partial U_1$, $x_1 > x_2$ should be understood to mean that $x_1$ is either above $x_2$ (in the sense of an ascending arc on $\partial U_1$ going from $x_1$ to $x_2$) or coincides with it. The context will fix the intended sense in each case.

<u>Lemma 3.2 :</u> Let $C_1$, $C_2$, $C_3$, $C_4$ (Fig. 27), be curves such that :

(i) $C_1$, $C_3$ are ascending and descending curves respectively in $\partial S_n$ – i.e, $C_1$ is ascending in going from $x_1$ to $y_1$ and $C_3$ descends from $x_2$ to $y_2$ - and $C_2$, $C_4$ are in $U_1$.

(ii) $C_1 \cup C_2$ and $C_3 \cup C_4$ are the boundaries of transversal discs which intersect transversally all flow lines they meet.

(iii) $C_1 \cap C_2 = \{x_1, y_1\}$ and $C_3 \cap C_4 = \{x_2, y_2\}$ ( Fig 27), where $q(x_1)$, $q(y_1)$ are on two different lateral curves of $X_1$, as also are $q(x_2)$ and $q(y_2)$.

Then it is not possible for the relations $x_1 > x_2$ and $y_2 > y_1$ to hold simultaneously.

Figure 27

Proof : In the following discussion, we need to distinguish two possibilities. These correspond to the two positions relative to each other of the two discs whose boundaries are $C_1 \cup C_2$ and $C_3 \cup C_4$, - viz, we can assume that these discs are either on the same side of $U_1$ or on opposite sides of it. In other words, the discs will be boundary intersecting of the first or second kind.

Case (i) : Here, $D_1$ and $D_2$ are on opposite sides of $U_1$ (Fig 28 (a)). The reasoning here would be the same as in Case (i) of Lemma 3.1, for boundary intersecting discs of the first kind.

Figure 28

Let $\partial D_1 = C_1 \cup C_2$ and $\partial D_2 = C_3 \cup C_4$, where $D_1$ and $D_2$ are discs transversal to flow lines. Let, if possible, an ascending curve $\alpha$ run from $x_2$ to $x_1$ and another ascending curve $\beta$ from $y_1$ to $y_2$ on $\partial U_1$. It is to be observed that $\alpha$ and $\beta$ are oriented in the positive direction of the flow. Hence $C = \alpha \cup C_1 \cup \beta \cup (-C_3)$ is an ascending curve, where $-C_3$ is the curve $C_3$ with negative orientation, which, hence, a*scends* from $y_2$ to $x_2$. Under these conditions, it can be seen that there is a transversal disc D in $\partial S_n$ with $C = \partial D$, using the same argument as in Case (i), Lemma 3.1. This implies that $\partial S_n$ is not acyclic. This contradiction establishes the impossibility of this case.

<u>Case (ii)</u> : Here both $D_1$ and $D_2$ are on the same side of $U_1$ – i.e, they are boundary intersecting discs of the second kind. To begin with, we know that $C_1$ at its end point $x_1$, is above $C_3$ at its end point $x_2$. Beginning with the point $x_2$, we now move along $C_3$. In so doing, we will have to reach a point $p_2$ on $C_3$ such that : (Fig. 28 (b)).

(i) $p_1$ on $C_1$ is above the point $p_2$ on $C_3$ (i.e on the same flow line).

(ii) there is a point $q_1$ on $C_1$ such that $q_1$ is below a point $q_2$ on $C_3$

(iii) between the points $p_2$ and $q_2$, no part of the interior of the sub-interval of the curve $C_1$ with end points $p_1$, $q_1$ is above or below any point of the interior of the sub-interval of the curve $C_3$ with end points $p_2$, $q_2$.

This can be seen to imply that either : (i) between $p_1$ and $q_1$, the curve $C_1$ intersects the interior of the disc $D_2$ or (ii) between $p_2$ and $q_2$ the curve $C_3$ intersects the interior of the disc $D_1$. The third possibility that $\partial D_1$ and $\partial D_2$ can be shown, without difficulty, to be reducible to one of these cases. Since (i) and (ii) are really identical situations, to dispose of this case, it is enough to show that (i) is impossible.

The discussion that follows will be with reference to Fig. 28(b).

Let the flow lines through $C_1$ between $p_1$ and $q_1$ cut $D_2$ in a curve $C_0$ which has $p_2$ and $q_2$ as end points. Since $D_2$ intersects all flow lines it meets transversally, $C_0$ is a smooth image of a sub-interval of $C_1$ under this projection. We now introduce the following notation :

(i) $C_{11}$ is the part of the curve $C_1$ between $p_1$ and $q_1$

(ii) $C_{21}$ is the part of the curve $C_3$ between $p_2$ and $q_2$

(iii) $I_1$ and $I_2$ are segments of flow lines with end points $p_1$, $p_2$ and $q_1$, $q_2$ respectively.

Let C be the simple closed curve which is the union of $I_1$, $C_{11}$, $I_2$ and $C_{21}$. Properly oriented, C becomes either an ascending or descending curve. If $C_{00} = C_0 \cup C_{21}$, then it is to be observed that C and $C_{00}$ satisfy the hypothesis of Principle [P]. Hence C is actually seen to be a singular ascending curve. Since C is in $\partial S_n$, this implies that $\partial S_n$ is not acyclic and this contradiction establishes the impossibility of this case.

This concludes the proof of the lemma.

Corollary 3.2 :  (i) An obvious modification of the proof yields the same result even when $x_1$ and $x_2$ or $y_1$ and $y_2$ coincide.

(ii) Another modification of the proof of Lemma 3.2 can be used to derive also the following result :

Let p and q be points in $\partial U_1$ which are the end points of a C-curve in $U_1$. Let $p_1$, $q_1$ on $\partial U_1$ be the end points of another C-curve or of an ascending (descending) curve on $\partial S_n$ − i.e, the curve ascends (descends) in going from $p_1$ to $q_1$. Then it cannot happen that $p_1$ is above (below) p and $q_1$ is below (above) q.

Lemma 3. 3 :  Let $x_1$ be a node in $X_1$. Then $a_1 > b_1$, if both $a_1$, $b_1$ exist. ($a_1$, $b_1$ are related to $x_1$ as in Definition 3.22).

Proof :  Let $p_1$, $r_1$, $s_1$ be nodes in $\partial U_1$ such that $q(p_1) = x_1$ ; $q(r_1) = a_1$ ; $q(s_1) = b_1$. It is to be observed that the possible positions of $r_1$ and $s_1$ are restricted by the condition that by definition either of them must be connected to $p_1$ by a smooth arc whose interior lies in the interior of $U_1$ and is transversal to all flow lines it meets. For $p_1$ as in Fig.29 (a), this requirement rules out the possibility of $r_1$ and $s_1$ being located anywhere in the interior of the vertical segment AB.

Figure 29

If possible, let $s_1 > r_1$, with an ascending arc $C_0$ running from $r_1$ to $s_1$ (Fig.29(b)). From the definition of $p_1$ and $r_1$, it can be seen that there is a transversal disc $D_1$ such that $\partial D_1 = C_1 \cup C_2$, with $C_1$ an ascending curve in $\partial S_n$ and $C_2$ a curve in $U_1$, both curves having $p_1$ and $r_1$ as end points. Let $C_1$ and $C_2$ have the orientation shown in Fig 30. Let $D_2$ be another transversal disc with $\partial D_2 = C_4 \cup C_3$, where $C_4$ is a descending curve in $\partial S_n$ and $C_3$ is a curve in $U_1$ with $C_4$, $C_3$ having the orientations shown. Also, $C_2$ and $C_3$ can be assumed to intersect only at the end point $p_1$.

Figure 30

The positions of the curves $C_1$ and $C_4$ is ruled out by Corollary 3.2. This contradiction establishes the impossibility of this case, and hence the lemma.

<u>Note :</u> If one or both of $a_1$, $b_1$ does not exist, then the Lemma becomes trivially true in view of the Convention below, set out after Definition 3.23.

<u>Definition 3. 23 :</u> For a node $x_i$ of $X_1$, the interval $I_i$ with end points $a_i$, $b_i$ on $X_1$ which does not contain the high or low point of $X_1$ in its interior, will be called the **<u>choice interval</u>** of $x_i$.

The terminology is suggested by the fact that the process of foliating $U_1$ is initiated by a foliation on a disc A in which $X_1$ is imbedded as boundary. This in turn will require as a first step, $x_i$ to be connected by a smooth curve in A to a point chosen in the interior of the choice interval.

<u>Convention :</u> Let H, L be the highest and lowest points respectively of $X_1$. If $A(x_i)$, $B(x_i)$ are both empty, $I_i$ will, by convention, be the lateral curve on $X_1$ with end points H and L not containing $x_i$. If $A(x_i)$ is empty, $I_i$ will, by convention, be the interval with end points H and $b_i$ on $X_1$ not containing $x_i$, while if it is $B(x_i)$ which is empty, $I_i$ will be the interval with end points $a_i$ and L not containing $x_i$. Since, by Lemma 3.3, $a_i$ and $b_i$ cannot coincide, the choice interval is always non-degenerate except in the following exceptional case : Here, $x_i$ is a node and $x_i$ is also one end point of a C-curve. If the image under q of the other end point of the curve of the ruling is $y_i$, then the choice interval $I_i$ of $x_i$ will be considered to be the degenerate interval represented by the point $y_i$.

<u>Note :</u> With the above convention, Lemma 3.3 becomes obviously true when $a_i$ = H or $b_i$ = L. The same observation holds regarding the other propositions in the sequel involving the choice interval.

The proof of the following Lemma is evident.

<u>Lemma 3.4 :</u> Let $z_0$, $x_0 \in \partial U_1$ be points such that $x_0$ and $z_0$ are connected by an ascending curve with $x_0$ a node. Let $q(x_0) = x_1$ and $q(z_0) = y_0$, where $x_1$ is a node. If $[a_1, b_1]$ is the choice interval of $x_1$, it can be seen that $a_1$ (and hence, the interval $[a_1, b_1]$) is lower than $y_0$ if it does not coincide with it.

<u>Definition 3. 24 :</u> In what follows, A will denote a fixed topological disc whose boundary is $X_1$ (Fig 31 (a)).

<u>Definition 3. 25 :</u> Let $x_1$ and $x_2$ be two nodes of $X_1$ and $I_1 = [a_1, b_1]$, $I_2 = [a_2, b_2]$ be their choice intervals. We say that $I_2$ is higher than $I_1$ and denote it by $I_2 > I_1$, if either $b_2 > a_1$ or if $b_2$ and $a_1$ coincide. (Fig 31 (b), (c)).

Figure 31

<u>Lemma 3. 5 :</u> Let $\alpha$, $\beta$ be two nodes of $X_1$ and $I_\alpha = [\alpha_1, \alpha_2]$, $I_\beta = [\beta_1, \beta_2]$ their respective choice intervals. Then neither of the following possibilities can occur. (i) $\alpha > \beta$ and $I_\beta > I_\alpha$ (Fig 32 (a)) (ii) $\alpha_2 > \beta$ and $\beta_2 > \alpha$. In other words, $I_\alpha > \beta$ and $I_\beta > \alpha$ (Fig 32 (b).

Figure 32

<u>Case (i) :</u> Let $q(x_1) = \alpha$, $q(x_2) = \beta$, $q(I_1) = I_\alpha$, $q(I_2) = I_\beta$, where $I_1$, $I_2$ are paths in $\partial U_1$ with end points $a_1$, $b_1$ and $a_2$, $b_2$ respectively. The first possibility implies that $x_1 > x_2$ and $I_2 > I_1$.

Let $C_{11}$ be the ascending curve running from $x_1$ to $a_1$.

Let $C_{12}$ be the descending curve running from $x_1$ to $b_1$.

Let $C_{21}$ be the ascending curve running from $x_2$ to $a_2$.

Let $C_{22}$ be the descending curve running from $x_2$ to $b_2$. (Fig 33 (a)).

The relative positions of the curve $C_{11}$ and $C_{22}$ contradict Lemma 3.2, showing this case to be impossible.

Figure 33

.

<u>Case (ii) :</u> The possibility $b_1 > x_2$ and $b_2 > x_1$ (Fig 33 (b) ).

Let $C_{11}$, $C_{12}$, $C_{21}$, $C_{22}$ be the curves described in Case (i). As before, the relative positions of the curves $C_{12}$ and $-C_{22}$ (which *ascends* from $b_2$ to $x_2$ ) contradict Lemma 3.2, ruling out this case.

Hence the lemma.

<u>Note :</u> The lemma can also be seen to hold when the lower end of $I_2$ coincides with the upper end of $I_1$ in Case (i), or when $b_1$ or $b_2$ coincides with $x_2$ or $x_1$, respectively, in Case (ii).

<u>Definition 3. 26 :</u> Let A be a disc having $X_1$ as boundary. Let $x_1$, $x_2$ be nodes in $X_1$ and $I_1$, $I_2$ the choice intervals of $x_1$ and $x_2$. If, for a point x in the interior of $I_1$, any curve $c_1$ with its interior in the interior of A and having $x_1$ and x as end points is such that $x_2$ cannot be connected to any point in int $I_2$ by a curve whose interior lies in int A and which does not intersect $c_1$, then $x_1$ and **<u>x will be said to obstruct $x_2$</u>**. We will also say that $(x_1, x)$ obstructs $x_2$ (Fig 34 (a), (b)).

Figure 34

<u>Lemma 3. 6 :</u> Let $x_1$ be a node in $X_1$ and $I_1$ its choice interval. Then $\exists$ a point x $\varepsilon$ int $I_1$ (whenever $I_1$ is non-degenerate) such that $(x_1, x)$ does not obstruct any other node in $X_1$.

<u>Proof :</u> The boundary of A can be considered to be the union of the two curves $C_1$ and $C_2$ which are the lateral curves of $X_1$ in Fig 35. For reasons explained earlier, it is clear that if a node $x_p$ is on $C_1$, its choice interval $I_p = [a_p, b_p]$ will be on $C_2$ and vice versa. To begin with, we shall assume $x_p$ to be on $C_{1,}$ say.

Figure 35

Let $S_1 = \{$ $a_i$ / $[a_i, b_i] = I_i$ is the choice interval of the node $x_i$, where $x_i > x_p \}$

And $S_2 = \{b_j$ / $[a_j, b_j] = I_j$ is the choice interval of the node $x_j$, where $x_j < x_p \}$

By Lemma 3.5 : (i) every point of $S_1$ is above $b_p$.

(ii) every point of $S_2$ is below $a_p$

(iii) every point of of $S_2$ is below every point of $S_1$ (Fig 36 ).

Figure 36

These conditions imply that there is a non-degenerate interval $I_{11}$ in int $I_p$ between $S_1$ and $S_2$ – i.e, every point of $S_1$ is above $I_{11}$ and every point of $S_2$ is below $I_{11}$. For x ε $I_{11}$, it can be seen that $(x_p, x)$ does not obstruct any node which is above or below $x_p$. Hence, $(x_p, x)$ does not obstruct any other node in $C_1$………………………………………………….(1)

We now define two other sets $S_3$ and $S_4$ as follows ( Fig 37 ).

$S_3$ = set of nodes in $C_2$ whose choice interval is below $x_p$

$S_4$ = set of nodes in $C_2$ whose choice interval is above $x_p$.

Figure 37

By Lemma 3.5, it follows that : (i) $S_3$ has to be below $S_1$ and $S_4$ has to be above $S_2$

(ii) $S_3$ is below $S_4$

(iii) $S_4$ is above $I_p$ and $S_3$ below $I_p$

(i) and (ii) imply the existence in $I_p$ of an interval $I_0$ which is below $S_1$ and $S_4$ and above $S_3$ and $S_2$ and by a suitable choice of $I_{11}$, we can have $I_0 \subset$ int $I_{11}$, with $I_0$ having the property that if x ε $I_0$, then $(x, x_p)$ does not obstruct any node in $C_2$ which is in $S_3 \cup S_4$………….(2)

Let $S_5$ be the set of nodes (in $C_2$) whose choice intervals contain $x_1$ as an interior point. If we assume that x has been chosen such that x $\notin S_5$, then $S_5$ is the disjoint union of $S_{15}$ and $S_{25}$, where $S_{15}$ are points of $S_5$ which are above x and $S_{25}$ those which are below x. For $x_k$ in $S_5$, there is obviously an open interval I having $x_p$ in its interior, such that I $\subset I_k$, the choice interval of $x_k$. It can be seen that for $x_k$ in $S_{15}$ ($S_{25}$), it is possible to connect $x_k$ to a point $y_k \in I$ above (below) $x_p$, by a curve which does not intersect a curve connecting $x_p$ and x.

Hence, $(x_p, x)$ does not obstruct any node in $S_5$ ……………………………………..……(3)

  From (1), (2) and (3) the Lemma follows.

<u>Lemma 3.7</u> :  Let $x_1$, $x_2$,……..,$x_k$ be the nodes in $X_1$ and $I_1$, $I_2$, …….., $I_k$ be their choice intervals. Then $y_j$ can be chosen in $I_j$ such that $(x_j, y_j)$ does not obstruct any of the other nodes – i.e, there are curves $c_j$ joining $x_j$, $y_j$ are such that $c_i \cap c_j = \emptyset$, for $i \neq j$.

<u>Proof</u> :  The proof is essentially a repeated application of the procedure in Lemma 3.6 to the $x_i$, taken one after another.

The nodes are first designated in the following fashion : the lowest node on $C_1$ will be labelled $x_1$. The next higher node on $C_1$ will be labelled $x_2$, etc, till we reach the node $x_m$ which will be the highest node in $C_1$. Then the lowest node in $C_2$ wil be called $x_{m+1}$ and the next higher node in $C_2$ will be $x_{m+2}$ and so on, till we reach the highest node $x_k$ on $C_2$.

For $x_1$, Lemma 3.6 shows how to choose $y_1$ in $I_1$ such that $(x_1, y_1)$ does not obstruct any of the $x_i$ for $i \neq 1$. The proof now proceeds by induction. Let us assume that for the nodes $x_1$, $x_2$,……,$x_{n-1}$, points $y_1$, $y_2$,……..,$y_{n-1}$ in their choice intervals $I_1$, $I_2$, …….,$I_{n-1}$ have been found such that there are arcs $c_i$ connecting $x_i$ and $y_i$, $1 \leq i \leq n-1$, with $c_i \cap c_j = \emptyset$ for $i \neq j$, $1 \leq i,j \leq n-1$ and the $c_i$ do not obstruct any of the $x_j$'s, $0 \leq j \leq k$. The following cases arise.

<u>Case (i)</u> :  $x_n$ is on $C_1$. By induction hypothesis, since $(x_1, y_1)$, $(x_2, y_2)$,…….., $(x_{n-1}, y_{n-1})$ do not obstruct $x_n$, $a_n$ will be above $y_1$, $y_2$,………, $y_{n-1}$, where $I_n = [a_n, b_n]$ is the choice interval of $x_n$. (It will be observed that $x_1$, $x_2$,……,$x_{n-1}$ are below $x_n$ in $C_1$ and $y_1$, $y_2$,……,$y_{n-1}$ will be in $C_2$, with $y_i > y_j$ for $i > j$) …………………………………………… (1)

The relationships between the relative positions of the sets defined below are displayed in Fig. 38.

Let $S_1 = \{a_i \ / \ [a_i, \ b_i]$ is the choice interval of $x_i$, where $x_i > x_n$ in $C_1\}$. Since $(x_1, y_1)$, $(x_2, y_2)$, ……. ,$(x_{n-1}, y_{n-1})$ do not obstruct any of the nodes $x_1$, ….,$x_k$, by induction hypothesis, $S_1$ (if non-empty) should be above $y_1, y_2, …..,y_{n-1}$…………………………………………(2)

Also, it can be seen from Lemma 3.5 that $b_n$ is below $S_1$…………….(3)

Let $S_2 = \{ b_i \ / \ [a_i, b_i] = I_i$ is the choice interval of $x_i$, and $x_i < x_n \}$. Again, by Lemma 3.5, $S_2$ is below $a_n$ and $S_1$ …………………..…………………..(4)

Let $S_{11} = \{ x_i \ / \ x_i$ on $C_2$ and $I_i$ the choice interval of $x_i$ is below $x_n \}$. From Lemma 3.5, $S_{11}$ is below $S_1$ and $a_n$ …………………………………………..(5)

Similarly, if $S_{21} = \{ x_j \ / \ x_j$ is on $C_2$ and $I_j$ is above $x_n \}$, $S_{21}$ is above $S_2$ and $b_n$.

…………………………..………..(6)

Also, it can be seen that $S_{11} < S_{21}$…………………………..……….(7)

From (1) – (7), it can be seen that it is possible to choose an interval $I_0$ in int $I_n$ such that $I_0$ is :

(i) below $S_1$ and above $S_2$ and $y_1, y_2, ……, y_{n-1}$

(ii) below $S_{21}$ and above $S_{11}$.

It can now be checked that $(x_n, y_n)$ does not obstruct any node, where $y_n \ \varepsilon \ I_0$ and is not a node.

It can be seen that the case of a point $x_i$ on $C_2$ which is such that $[a_i, b_i]$ contains $x_n$ in its interior, offers no difficulty over being obstructed by $(x_n, x)$.

Figure 38

<u>Case (ii) :</u> Let $I_i$ be the choice interval of $x_i$ for $i \leq (n-1)$. The induction hypothesis now ensures the following facts :

(i) There exists $y_i \, \varepsilon \, I_i$ and curves $C_i$ connecting $x_i$ and $y_i$ such that $C_i \cap C_j = \Phi$, for $i,j \leq (n-1)$, $i \neq j$.

(ii) If $S_3 = \{a_i \, / \, [a_i, \, b_i] = I_i$ is the choice interval of $x_i$ for $x_i > x_n\}$, then $S_3$ is above $b_n$, where $[a_n, b_n]$ is the choice interval of $x_n$. (Since $x_n$ is assumed to be on $C_2$, $x_i > x_n$ will be on $C_2$ and hence $I_i$ will be on $C_1$).

(i) With the $y_i$ as in (i), if $x_n$ is in the interior of the interval with end points $y_p$, $y_{p+1}$ on $C_2$, there will be a point $y_n$ in the choice interval $I_n$ (in $C_1$) of $x_n$ such that $y_n$ is in the interval in $C_1$ with end points $x_p$, $x_{p+1}$ and $y_n$ will also be below $S_3$.

It can be seen that a curve $C_n$ with end points $x_n$, $y_n$ can be drawn such that $C_n$ does not obstruct any of the other nodes.

Hence the lemma.

Figure 39

In what follows, we will be generalising the concept of node introduced in Defn. 3.19.

<u>Definition 3.27 :</u> Let $x \in \partial S_n$. Then $x$ is a **<u>vertex</u>** if :

(i) $x$ is the interior of a sub-interval $I$ of the flow line through $x$, with $I \subset \partial S_n$

(ii) At a point $y$ arbitrarily close to $x$ on the curve of the ruling through $x$, it is not possible to find a sub-interval $I_y$ of the flow line through $y$ such that $y \in$ int $I_y$ and $I_y \subset \partial S_n$.

By this definition, the points where a ruling segment and a flow segment in the boundary of an unfoliated region meet, would be vertices. (Fig. 40).

Let T be a tube not in $S_n$, to which we wish to extend the foliation on $S_n$ and for which $T \cap S_n$ $\neq \emptyset$. It will be noticed that a significant difference between the above definition of vertex and the earlier one of node is that we have relaxed the condition that x should be connected to a point on $\partial U_1$, where $U_1$ is an unfoliated region of $\partial T$, by a curve of the ruling. Apart from such points, other possibilities for a node may include $X \cap S_n$, where X is the space described at the beginning of Section 3 on which we begin constructing the tubes.

We need one other fact introduced by the following situation. Let $C_1$, $C_2$ be two curves of the ruling on $\partial S_n$. Let x be a point on $C_1$. x can now move in one of two possible directions, say, right and left, along the curve of the ruling on which it is situated. For a position $x_i$ of x on $C_1$, let the flow line through $x_i$ first intersect $C_2$ at a point $y_i$ on $C_2$ (below $x_i$, say). We assume that the part of the segment of the flow line between $C_1$ and $C_2$ as x moves to the right along $C_1$, lies on $\partial S_n$. It may now happen that x reaches a position $x_1$ at which the following happens : $x_1$ is the left end point of an interval $I_1$ in $C_1$ such that the flow line through $x_1$ meets $C_2$ at a point $y_1$ on it, but for any $x \in I_1$ the flow line through which does not pass through $x_1$, the flow line through x does not meet $C_2$ (Fig. 40).

In what follows, we shall be concerned chiefly with the situation in which x is a point on the boundary $\partial D_1$ (which is a union of curves of the ruling and segments of flow lines) of a transverse singular disc $D_1$ and x begins to move along $\partial D_1$. ($C_1$, $C_2$, will, therefore, initially lie on $\partial T$). In this context, we will also additionally require the flow lines through x which meets $C_2$, to lie on $\partial T$ at the point when x begins to move. In, fact, when we encounter this situation for the first time in Theorem 3.1, which follows, $C_1$ and $C_2$ intersect at $x = z_1$, when x begins its movement.

.

Under these conditions, for two positions α, β of x on $C_1$ which are very close to each other, if α₁, β₁ are the points at which the flow lines through α, β meet $C_2$, then the segments of flow lines with end points α, α₁ and β, β₁ should be very close to each other.

We now introduce the following lemma. The proof is evident.

Lemma 3.8 :  The above situation can happen only if in the interval I of the flow line through $x_1$, with end points $x_1$, $y_1$, with $y_1 \in C_2$, there is a point P which is a vertex. (It may be possible that this exterior node may be $x_1$ itself. Barring this exceptional case, P will be below $x_1$). Further, of the two segments which meet at P -  namely the flow segment and the ruling segment – a sub-interval $I_2$ of the latter will have P as an end point and will lie on the side of P on which it will intersect flow line segments which will have end points on $C_1$ and $C_2$. The figure below illustrates some of the possibilities for P and $I_2$.

Figure  40

Note 3.8 :  (i) A corresponding result is true when x moves along $C_1$ to the left or when $C_1$ is below $C_2$.

(i) It is possible to define "right" and "left" coherently on $\partial S_n$ with respect to the orientation of the flow lines due to the fact that $A_1$ and therefore $\partial S_n$, are orientable. Also, there can be more than one vertex in the interval I of Lemma 3.8. Fig 40(b) illustrates this possibility.

Definition 3.28 :  In Lemma 3.8, when P is the left end point of $I_2$, as in Fig.40, it will be called a **right vertex** and if it is the right end point it will be called a **left vertex**.

 Definition 3.29 :  Let us consider the case when x moves left on a curve $C_1$. In the course of such a movement, x may reach a position $x_1$ such that the flow line through x meets $C_2$ for x close to $x_1$ and to the right of it, while it fails to meet $C_2$ for x close to $x_1$ on the left. x is then said to encounter a **_left obstruction_** on account of a right or left exterior node, according as P is a right or left vertex. A right obstruction is similarly defined. Fig. 40 shows a number of cases in which x encounters a right or left obstruction shown as end point of an interval $I_2$.

Definition 3.30 : A **nodal line** is a curve of the ruling passing through a vertex

The following lemma will be used in the proof of Theorem 3. 1.

**Lemma 3.9 :** Let $C_1$, $C_2$ be the lateral curves of $X_1$ and $S_i \subset C_i$, i = 1,2, be the end points of C-curves. Let f: $S_1 \rightarrow S_2$ be defined by f(x) = the end point of the curve whose one end point is x $\varepsilon$ $S_1$. f is then a homeomorphism.

Proof: The proof is a consequence of the fact that there are only a finite number of nodes and that to get f(x), we just follow the curve of the ruling with the end point at x, around $\partial$T till we reach f(x).

**Theorem 3.1** : The unfoliated region $U_1$ of $\partial$T can be foliated by a system S of curves such that:

(i) S $\cup$ { the system $S_1$ of curves of the ruling $R_n$ gives a smooth foliation for [$\partial S_n \cup U_1$]

(ii) If $C_0$ is a C - curve in $U_1$, there is a member of S with the same end points as $C_0$.

(iii) the curves of S cut all flow lines they meet transversally.

(iv) there is no singular ascending or descending curve in S $\cup$ $S_1$.

Proof : By Note 3.13.4, we can assume $U_1$ is Ruling Convex and use the properties of a Ruling Convex region listed in Note 3.14.1. By Lemma 3.7, it becomes possible in A to join nodes $x_i$ to points $y_i$ in the interior of their choice intervals (whenever these are non-degenerate) by smooth curves $c_i$ such that the $c_i$ do not intersect. We now extend these curves to a foliation on the space A by a system of smooth curves $S_0$ that is transverse to the two boundary curves of $X_1$ that it meets. We also require that $S_0$ have the further important property that if $z_1$, $z_2$ are two points of $X_1$ which are the images under q (q is, of course, the quotient map from $\partial U_1$ to $X_1$ of Definition 3.17) of the end points in $\partial U_1$ of a C-curve in $U_1$, then $z_1$, $z_2$ are connected by a curve of $S_0$. A brief explanation would be in order here about why this is always possible.

We first connect points like $z_1$ and $z_2$ by a smooth family $S^0$ of curves (which would be C-curves) which are pairwise disjoint and then draw the curves $c_i$ so that:

(i) none of them intersects any member of $S^0$ and

(ii) The $C_i$ are pairwise disjoint.

Lemma 3.1 shows that the members of $S^0$ can be drawn without intersecting each other, while Corollary 3.2 guarantees the legitimacy of the latter requirement. (It is worth observing that the union of the members of $S^0$ can be shown to be a closed set. This follows from Lemma 3.9). In order to get $S_0$ now, it is only necessary to extend the $c_i$'s and the members of $S^0$ to a smooth foliation of A satisfying the condition that the curves of the foliation intersect transversally the two boundary curves of A, while the end points of each curve lie one on each lateral curve of $X_1$. (There are only a finite number of nodes and hence a finite number of $c_i$'s).

We now choose a map p of $U_1$ to A which is such that $p/\partial U_1 = q$ and p is a diffeomorphism in the interior of $U_1$. P and $S_0$ can evidently be chosen such that the flow lines of $U_1$ are carried by p to lines which cut members of $S_0$ transversally. It can be seen that $U_1$ can be foliated by a system of curves S with the following properties :

(i) p carries the curves of the system S to those of $S_0$.

(ii) The curves of S intersect transversally all flow lines they meet. It follows that the members of S cut transversally all boundary segments of $\partial U_1$ which are made up of flow segments.

(iii) The union of curves of S and those of the system $S_1$ of the curves of the ruling $R_n$ on $\partial S_n$ gives a smooth foliation by a system of curves for $\partial S_n \cup U_1$.

(iv) for every C- curve $C_0$ in $U_1$, there is a member of S with the same end points as $C_0$.

**To show that there are no singular curves among the curves of $S \cup S_1$**.

Let, if possible, C be an ascending singular curve in the system of curves $S_{11} = S \cup S_1$. Let C = $C_1 \cup C_2$ where $C_1$ is a curve in $\partial S_n$ and $C_2$ a curve of the system S  in $U_1$. ( If C does not intersect int $U_1$, or even if the intersection of C with $U_1$ happens to be merely a flow line segment of C, then C would be in $\partial S_n$ and this would imply that $\partial S_n$ is not acyclic). If $C_1$ is not an ascending curve, then $C_1$ has to be the boundary of a leaf of $S_n$. This implies that $C_1$ is a C-curve and the end points $z_1$, $z_2$ of $C_1$ on $\partial U_1$ are the end points of a C-curve on $U_1$. In other words, $C_2$ will have to be a C-curve and hence C = $C_1 \cup C_2$ cannot be an ascending curve – i.e, it will not contain any non-degenerate sub-interval of a flow line. This follows from the fact that $C_2$ will have to meet all flow lines transversally. Hence it must happen that $C_1$ is an ascending curve.

Let $C_{11}$ be the curve of the ruling in $\partial S_n$ through $z_1$. As x moves on $C_1$, beginning from the initial point x = $z_1$, we consider below the following possibilities regarding the intersection with $C_{11}$ of the flow line through x. (There is a small ambiguity here, as there will be two curves of the ruling through x – one of them lying on $\partial T$ in a neighbourhood of x and the other intersecting $\partial T$ transversally at x. But the ensuing discussion will apply with equal validity to either alternative).

Case (i)(a):  As x moves on $C_1$ towards the point $z_2$, the flow line through x continues to intersect $C_{11}$ at a point below it (or at x itself, when x is on $C_1 \cap C_{11}$), till x reaches a position $z_0$ ($z_0$ may coincide with $z_2$) at which the flow line through x = $z_0$ meets $C_{11}$ at a point $z_{11}$ where $C_{11}$ first meets $\partial U_1$. (Fig.41).

Figure 41

Let $I_1$ be the segment of the flow line with end points $z_0$ and $z_{11}$ and $I_{11}$ the sub-interval of $C_1$ with end points $z_0$ and $z_2$. It can then be seen that either $(C_{11} \cup I_1 \cup I_{11} \cup C_2)$ in Fig. 41(a) or $C_{11} \cup I_1 \cup C_2$ in Fig. 41(b) respectively, is a simple closed curve. Denoting either of these curves by $C_{12}$, we see that $C_{12}$ intersects every flow line through points of $C_1 \cup C_2$ in a connected set. Principle [P] now applies (to the curves $C_1 \cup C_2$ and $C_{12}$) to show the existence of a singular transversal disc whose boundary is $C_{12}$. This implies that $z_1$ and $z_{11}$ can be connected by a curve $C_3$ in $U_1$ which is such that the flow lines through any point of it intersects the curve $(I_1 \cup I_{11} \cup C_2)$ in a connected set, while $C_3$ itself cuts all flow line it meets, transversally. Applying Principle [P] to $C_{12}$ and $C_{11} \cup C_3$, we see that the latter should be the boundary of a transversal disc and hence $z_1$ and $z_{11}$ are the end points of a C-curve. Hence, by the convention followed in constructing the system S of curves that foliate $U_1$, there will be a curve of S connecting $z_{11}$ and $z_1$. But we have already assumed that $z_1$ is connected to $z_2$ by a curve of the system S, while $z_{11}$ has to be below $z_2$, since $C_1$ is an ascending curve. This contradiction rules out this case.

<u>Case (i)(b)</u> : Here, when x reaches $z_2$, the corresponding point on $C_{11}$ has not yet reached $\partial U_1$. We now allow x to move on $C_2$ from $z_2$ to $z_1$ and look at the corresponding point on $C_{11}$. Let us assume that when x reaches some position in the interior of the curve $C_2$, the corresponding point on $C_{11}$ reaches a point on $\partial U_1$. (Fig. 42).

Figure 42

A slight modification of essentially the same argument used above would be sufficient to demonstrate the impossibility of this case too.

<u>Case (i)(c)</u> :  Here, as x moves around $C_1 \cup C_2$ and completes a full circuit by returning to the point $z_1$, the corresponding point on $C_{11}$ reaches a point $z_{12}$ without $C_{11}$ meeting $\partial U_1$ anywhere. (Fig. 43).

Figure  43

Let $I_1$ be the interval with end points $z_1$ and $z_{12}$ on the flow line through $z_1$. In the sequel, we will designate by $C_{11}$ some part of the curve of the ruling through $z_1$ which has been selected for attention. In particular, in the rest of this proof $C_{11}$ will be the segment of the curve of the ruling through $z_1$ with end points $z_1$ and $z_{12}$.  It is to be observed that every point of the curve $C = C_{11} \cup I_1$ is below some point of the curve $C_1 \cup C_2$ and intersects flow lines through every point of the latter in a connected set. Principle [P] can now be applied to the curves $C$ and $C_1 \cup C_2$ to show that there is a singular transversal disc whose boundary is $C$. But since $C$ is an ascending curve in $\partial S_n$, this implies that $\partial S_n$ is not acyclic. Hence, this case cannot arise.

<u>Case (ii)</u> : In this case, as x moves to the right on $C_1 \cup C_2$, beginning at the point $z_1$, it reaches a position $x = x_1$, such that for points immediately to the right of $x_1$, the flow line through x will fail to meet $C_{11}$ while it does so at $x = x_1$. This is the situation described in Lemma 3.8 and by that Lemma, there is a vertex p on the interval $[x_1, y_1]$, where $y_1$ is the point at which the flow line through $x_1$ meets $C_{11}$.

Figure   44

We now consider a point y which moves on $C_{11}$ back towards $z_1$, beginning at the point y = $y_1$. Let the flow line through y meet the curve of the ruling $C_{12}$ through p above it, at a variable point we will denote by $\alpha$. As y moves towards $z_1$, it may happen that it encounters a left obstruction at a point y = $y_2$. Let the corresponding point on $C_{12}$ be $\alpha_2$. According to Lemma 3.8, there is a vertex $p_1$ on the interval of the flow line through $y_2$ with end points $\alpha_2$ and $y_2$ (Fig. 45). There are now two possibilities to be considered.

Case (ii)(a) :  $p_1$ is a left vertex. (Fig. 45). It may be possible that $p_1$ may coincide with $y_2$.

Figure   45

Let $I_1$ be the closed sub-interval of the flow line with end points $\alpha_2$ and $y_2$. Let $\gamma_1$ be the sub-interval of the flow line through $y_2$ with end points $\alpha_2$ and $p_1$. Since $p_1$ is a left vertex, there is a nodal line through $p_1$ and lying to the left of it, which we will denote by $C_{13}$.

As y moves on $C_{11}$ to the left of $y_2$, we can continue to consider the intersection with $C_{13}$ of the flow line through y. (By note 3.8, there may be more than one exterior node in $I_1$. In this case, $p_1$ will be the exterior node closest to $y_2$).

Case (ii)(b) :  $p_1$ is a right vertex. This can happen, for instance,  if $p_1$ is a node on the boundary of an unfoliated region $U_k$ of $\partial T$, as shown in Fig.46.

Figure   46

In this case, if we consider the interval $I_1$ of the flow line with end points $\alpha_2$, $y_2$, we find that one end point $\alpha_2$ is on $\partial U_k$ and the other end point $y_2$ is outside $U_k$. (At worst, $y_2$ can lie on $\partial U_k$ – in which case $y_2$ will lie on a nodal line. This will not affect the rest of the argument. But it has to be observed $C_{11}$ cannot intersect the interior of $U_k$, since this region can contain no curves of the ruling). Hence the interval $I_1$ <u>will enter the interior of $U_k$ at a vertex $p_1$</u>. It will therefore intersect the boundary of $U_k$ again at some point $\alpha_{12}$ (Fig. 46), at <u>a nodal line which we will denote by $C_{13}$.</u> Hence, it will be now possible to continue the procedure of considering the point of intersection of the flow line through y with the curve of the ruling through $\alpha_{12}$, as, y moves to the left of $y_2$ along $C_{11}$.

The sub-interval of the flow line with end points $\alpha_2$ and $\alpha_{12}$ will be denoted by $\gamma_1$ (the same symbol as in Case (ii)(a) is intentionally employed).

In either of the above cases, it is thus seen to be possible to continue the process of considering the intersection of the flow line through y, as y moves along $C_{11}$ to the left with, first, the curve $C_{12}$ and then with $C_{13}$. In fact, the process may be continued in stages along curves $C_{12}$, $C_{13}$, $C_{14}$,…….. as y encounters left obstructions at $y_2$, $y_3$, $y_4$,…… . Further, it will be observed that the curves $C_{12}$, $C_{13}$, $C_{14}$,….. are <u>nodal lines</u>. At each of these points $y_i$, i = 1,2,3….. let the sub-intervals of the curves of the ruling and the flow lines that meet at the points $\alpha_{1i}$ above $y_i$ be denoted by $\beta_1$, $\gamma_1$, $\beta_2$, $\gamma_2$, $\beta_3$, $\gamma_3$ …….. respectively, so that the union of these segments is a path which describes the trajectory of the point which begins its journey at P and is on the flow line through y and above it as y as y moves to the left along $C_{11}$. There are now three possibilities to be considered. These correspond, by analogy, to the earlier Cases (i)(a), (i)(b) and (i)(c).

**Possibility (a)** : As y moves on $C_{11}$, for some position of y = $z_0$ in $C_{11}$, let the flow line through $z_0$ meet the nodal line $C_{1k}$ at a point q on $\partial U_1$. ($z_0$ may coincide with $z_1$). In other words, y will reach $z_0$ after encountering left obstructions at $y_2, y_3, \ldots\ldots, y_{k-1}$.

Figure 47

This corresponds to the Case (i)(a). It will be noted that since $C_{1k}$ is a nodal line, q itself becomes a node. (In the extreme case, q, $z_0$ and $z_1$ may coincide ). Let $I_0$ be the sub-interval of the flow line through $x_1$, with end points $x_1$ and p. Then, the curve $\alpha_0 = I_0 \cup \beta_1 \cup \gamma_1 \cup \beta_2 \cup \gamma_2 \cup \ldots.\cup \beta_{k-1}$ is a curve connecting q and $x_1$ and descends from $x_1$ to q − i.e, it ascends from q to $x_1$. (Fig. 48).

Figure 48

Let $\alpha_1$ be the closed segment of the curve $C_1$ between $x_1$ and $z_2$. We look at the curve $\lambda = \alpha_0 \cup \alpha_1 \cup C_2 \cup C_3 \cup I$ where $C_3$ is the sub-interval of the curve of the ruling with end points $z_1$ and $z_0$ (when $z_0$ and $z_1$ coincide $C_3$ may become degenerate) and I is the segment of the flow line through q with end points $z_0$ and q. By applying Principle P to $\lambda$ and $C_1 \cup C_2$, we see that $\lambda$ is the boundary of a transversal disc. Since q is a node and $\lambda$ is obviously an ascending curve in the part of it that lies outside $U_1$, it is seen to be a singular curve. By Lemma 3.4, if [a, b] is the choice interval of q, a will be below $z_2$ or coincide with it. Hence, there will be a curve of the system S in $U_1$ connecting q to <u>a point below $z_2$,</u> and so this possibility cannot arise.

**Possibility (b)** :  Here, by analogy with Case (i)(b), we assume that when y reaches $z_1$, the corresponding point on the curve above it has not yet reached $\partial U_1$.

In this case, we allow y to move from $z_1$ to $z_2$ along $C_2$, using the procedure outlined in Possibility (a), <u>wherever y encounters a left obstruction.</u>  At some intermediate position $z_0$ between $z_1$ and $z_2$, we assume that the flow line through $y = z_0$ intersects the curve $C_{1m}$ at a point $q_1$ at which $C_{1m}$ first intersects $\partial U_1$ (Fig. 49).

Figure 49

 The following observations can be made about $q_1$ :

(i) $q_1$ is a node since $C_{1m}$ is a nodal line.

(ii) Let $\alpha_0 = I_0 \cup \beta_1 \cup \gamma_1 \cup \beta_2 \cup \gamma_2 \ldots \ldots \cup \beta_{m-1}$, where $I_0$ is the interval of the flow line with end points $x_1$ and p, and the $\beta_i$'s and the $\gamma_i$'s are defined as in the earlier case. Then $\alpha_0$ is a curve that ascends in going from $q_1$ to $x_1$.

As in the earlier case, let $\alpha_1$ represent the sub-interval of $C_1$ with end points $x_1$ and $z_2$.We will denote by $\alpha_2$ the sub-interval of $C_2$ with end points $z_2$ and $z_0$ and by I the sub-interval of the flow line through $z_0$ with end points $z_0$ and $q_1$. By analogy with the earlier case, we can consider the curve $\lambda = \alpha_0 \cup \alpha_1 \cup \alpha_2 \cup I$. This is a closed ascending curve which can be considered to begin and end at $q_1$. Applying Principle P to the curves $\lambda$ and $C_1 \cup C_2$, we see that the part of $\lambda$ outside $U_1$ is an ascending singular curve. Exactly as in Possibility (a), this can be seen to lead to a contradiction, thereby ruling out this case.

Possibility (c) :  Here, as in Case (i)(c), as y moves on $C_2$ from $z_1$ to $z_2$, let the flow line through y continue to meet the curve above it till y reaches $z_2$, without the curve above it reaching $\partial U_1$ (Fig. 50). $q_1$ will now denote the point at which the flow line through $z_2$ meets the curve $C_{1r}$ above it. Again, $\alpha_0$ will denote the curve $I_0 \cup \beta_1 \cup \gamma_1 \cup \beta_2 \cup \gamma_2 \cup \ldots \cup \beta_{r-1}$. $\alpha_1$ will be the sub-interval of $C_1$ with end points $x_1$ and $z_2$ and $I$ the sub-interval of the flow line through $z_2$ with end points $z_2$ and $q_1$ (Fig. 50). Let $\lambda = \alpha_0 \cup \alpha_1 \cup I$. The rest of the argument proceeds as before by applying Principle P to the curves $\lambda$ and $(C_1 \cup C_2)$. The Principle implies that $\lambda$ is a singular curve and since $\lambda$ is in $\partial S_n$, we reach the contradiction that $S_n$ is not acyclic. This disposes off this possibility.

Figure  50

This establishes the absence of singular curves on the rulings of $\partial S_n \cup U_1$. Hence, the system of curves of the ruling on $\partial S_n \cup U_1$ is acyclic.

This procedure is now carried out with each unfoliated region $U_2$, $U_3$, ........., $U_k$ of $\partial T$ in turn, to reach a foliation on $\partial T$ without any ascending singular path. Consequently, the resulting foliation on $\partial S_n \cup \partial T$ – and hence on $\partial T$ – is acyclic. In particular, this implies that each curve of this foliation of $\partial T$ is a circle. T is now foliated by discs whose boundaries are precisely these circles, such that the union of the leaves foliating $S_n$ and these discs gives a smooth foliation for $S_{n+1} = S_n \cup T$.

By foliating each tube successively by the procedure described in theorem 3.1, we arrive at the following theorem.

**<u>Theorem 3.2</u>** : Beginning with the acyclic ruling $\Omega$ on the tube space $S_0$, the progressive application of the foliation procedure in the presence of the Constraint Condition gives a smooth co-dimension one foliation on $S^3$ which is transverse to the flow v.

By Novikov's theorem ([2]), there is a compact leaf F of this foliation. We now make the following assumptions regarding F :

(i) $(S^3 - F)$ consists of two components $C_1$ and $C_2$.

(ii) The space $S_0$ with which we began the Foliation process is in $C_1$.

(iii) The vector field v/F gives a field of "outward" pointing vectors on F, with respect to its interior $C_1$. This follows from the fact that v is transversal to F.

It is now evident that at least one flow line originating at a point of $S_0$ will intersect F before returning to $S_0$, since otherwise all flow lines originating at $S_0$ will be contained in $C_1$. This will imply that $S^3$ itself is contained in $C_1$. But this flow line will leave the component $C_1$, but by condition (iii) above, cannot return to $S_0$, violating the fundamental property of v, that every flow line originating from a point of $S_0$ will return to it.

A similar argument serves to establish that $S_0$ cannot be in $C_2$ or F. This implies the non-existence of a minimal flow on $S^3$.

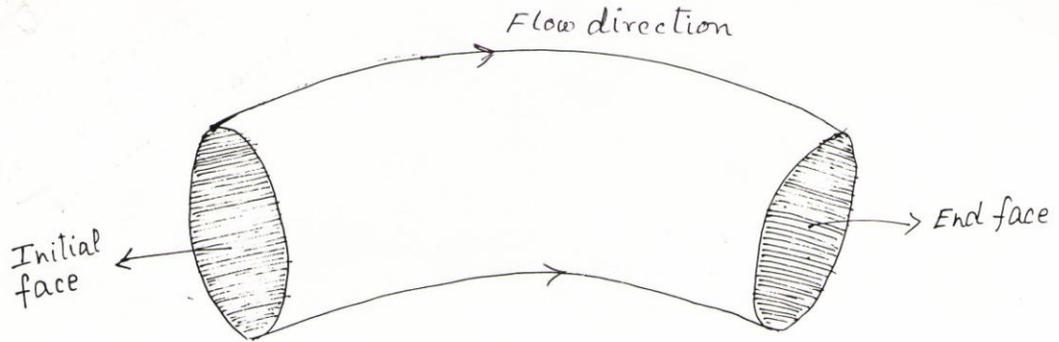

Flow direction

Initial
face

End face

FIGURE 1

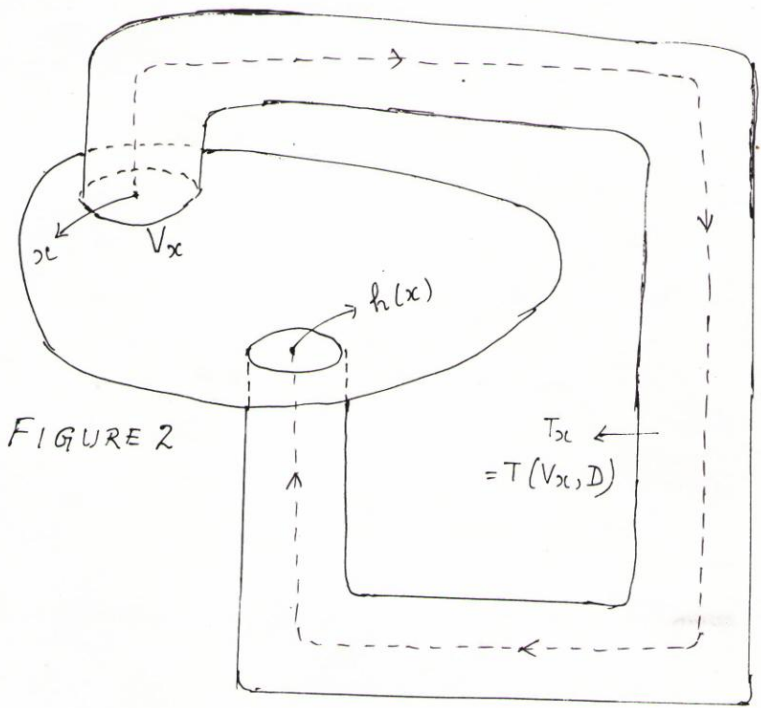

$x$

$V_x$

$h(x)$

FIGURE 2

$T_x$
$= T(V_x, D)$

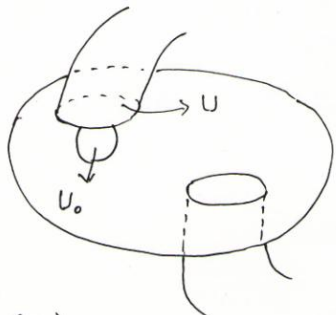

$U$

$U_0$

(a)

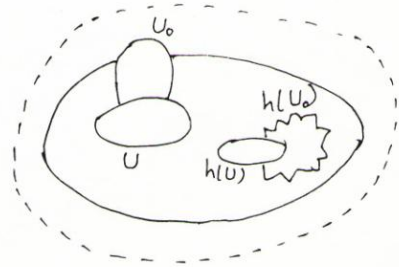

$U_0$

$h(U_0)$

$U$

$h(U)$

FIGURE 3

(b)

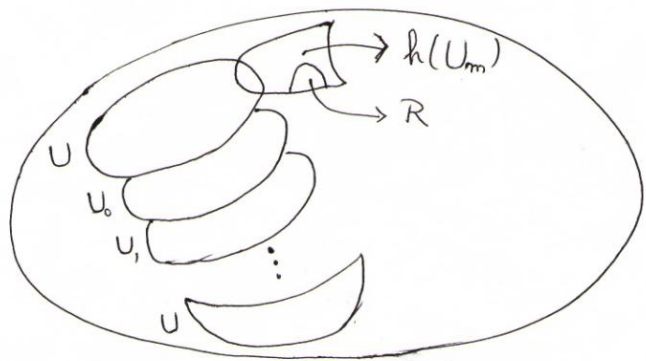

$$F_{IGURE} \, 4$$

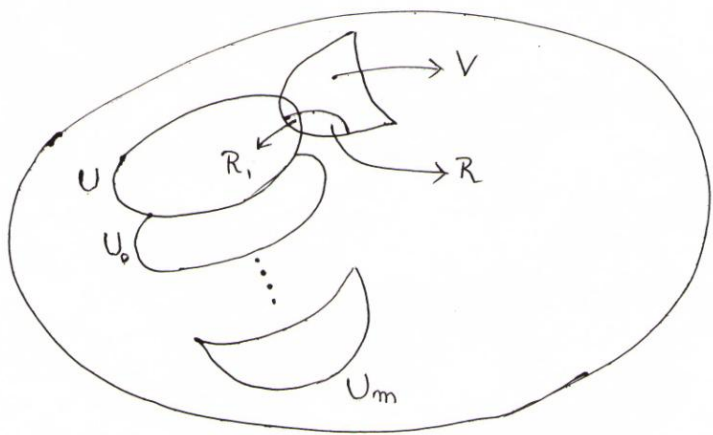

$$F_{IGURE} \, 5.$$

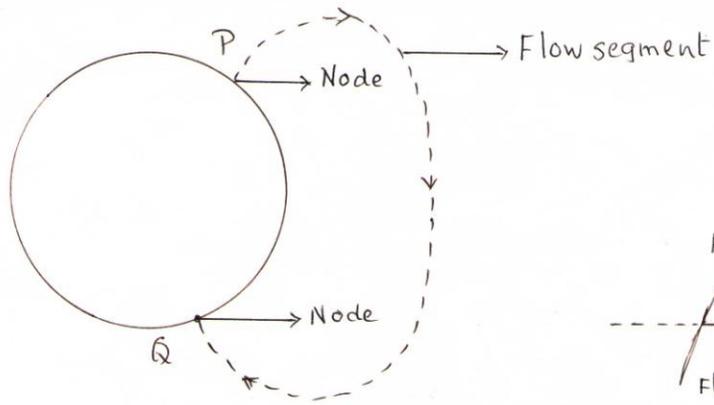

P

→ Node

→ Flow segment

Q

→ Node

(a)

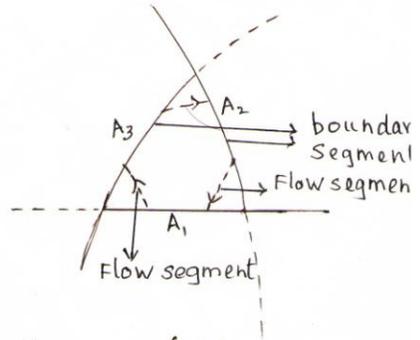

A₃

A₂

→ boundary Segment

→ Flow segment

A₁

Flow segment

FIGURE 6

(b)

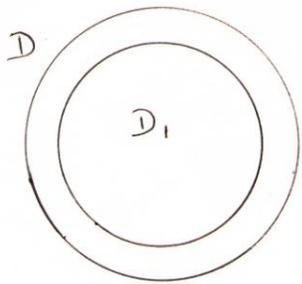

D

D₁

FIGURE 7

(a)

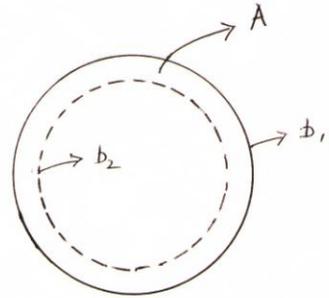

A

b₂

b₁

(b)

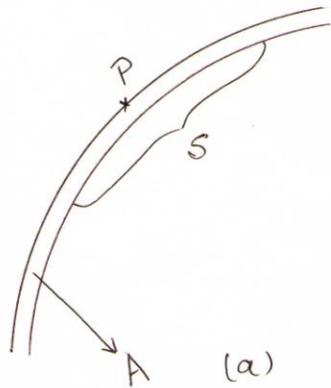

P

S

A

(a)

FIGURE 8

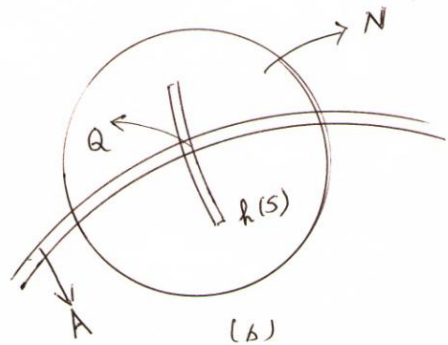

N

Q

l(S)

A

(b)

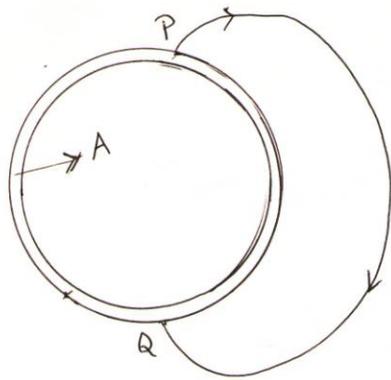

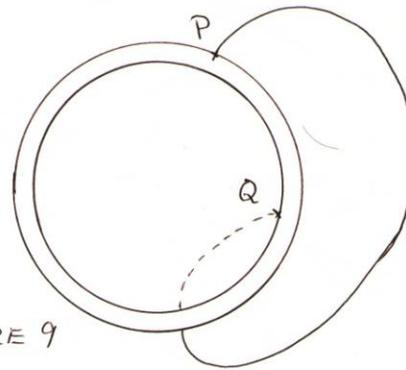

FIGURE 9

(a) Nodal segments
of the first class

(b) Nodal segments of the
second class.

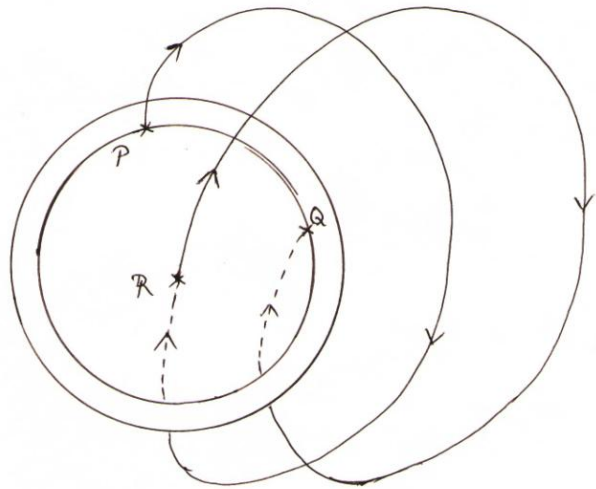

FIGURE 10

An irregular singular A-segment which intersec
the interior of D, at the point R.

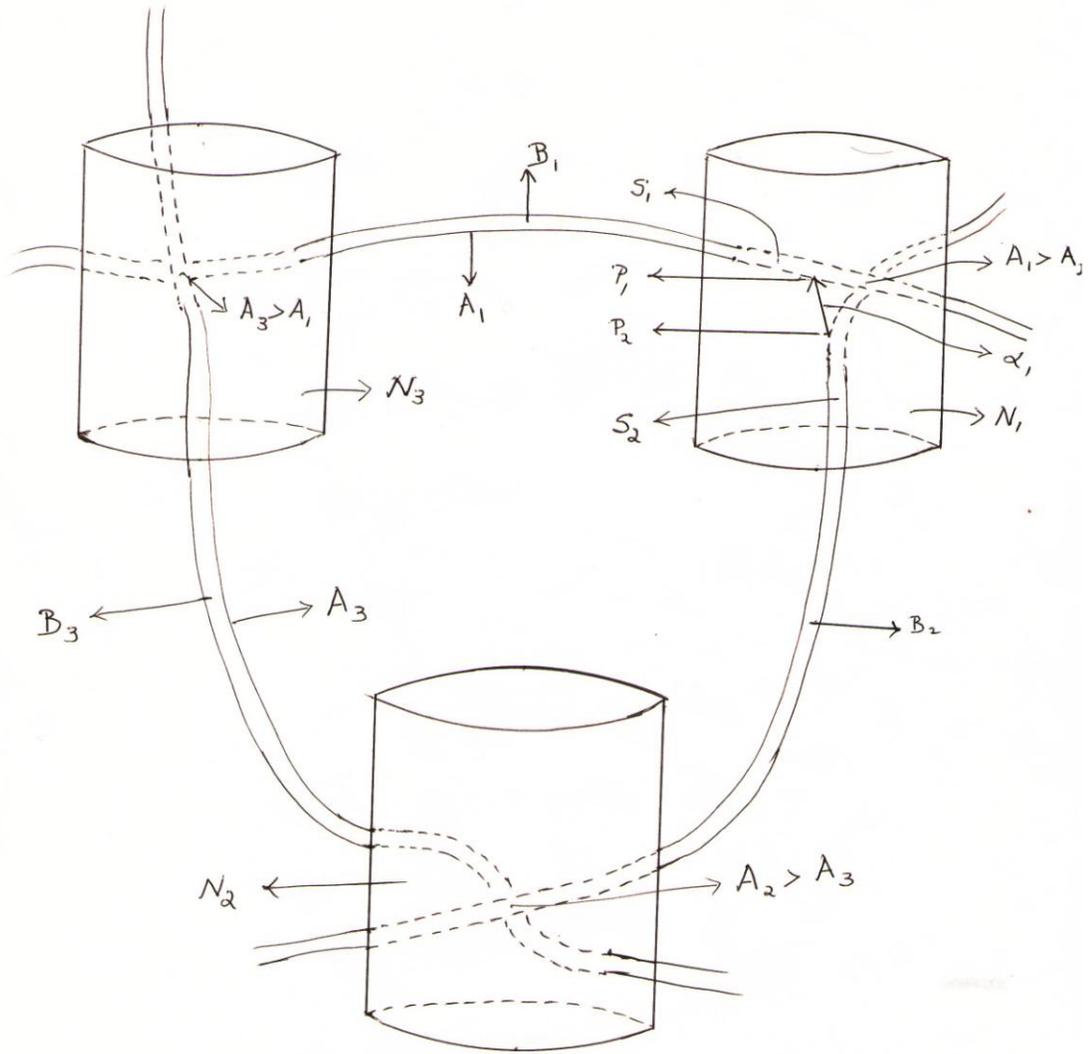

FIGURE 11

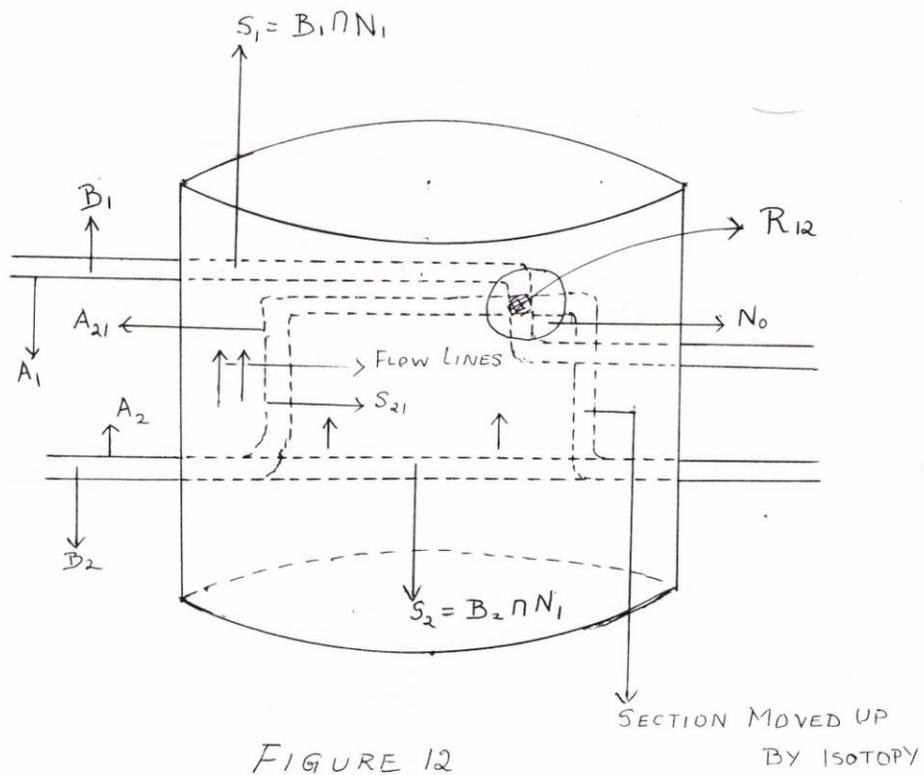

$S_1 = B_1 \cap N_1$

$B_1$

$R_{12}$

$A_{21}$

$N_0$

$A_1$

Flow Lines

$A_2$

$S_{21}$

$D_2$

$S_2 = B_2 \cap N_1$

Section Moved up
By Isotopy

*FIGURE 12*

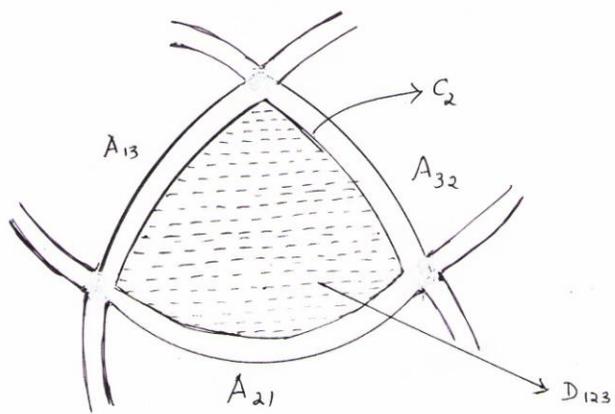

$C_2$

$A_{13}$

$A_{32}$

$A_{21}$

$D_{123}$

*FIG. 12.1*

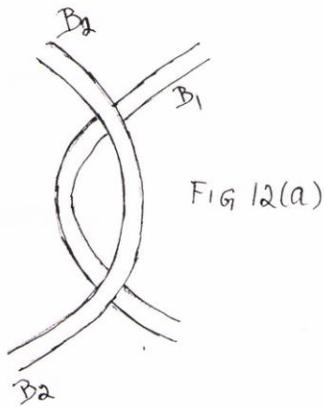

FIG 12(a)

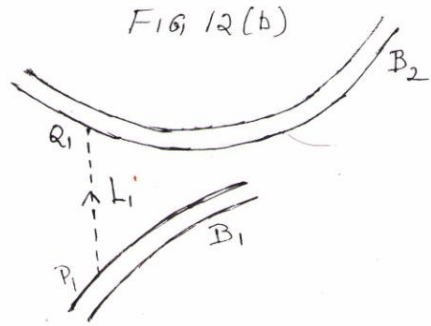

FIG 12(b)

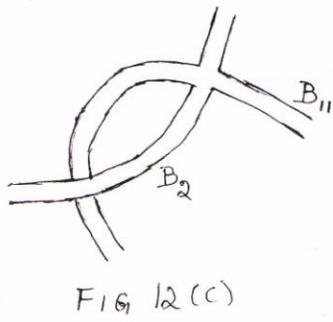

FIG 12(c)

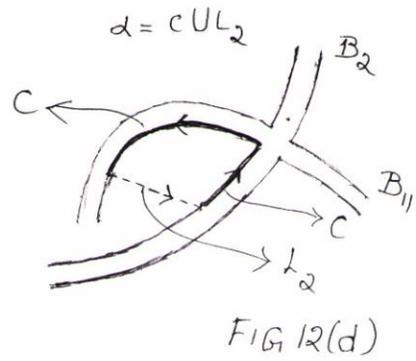

$\alpha = C \cup L_2$

FIG 12(d)

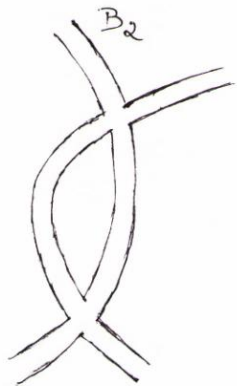

FIG 12(e)

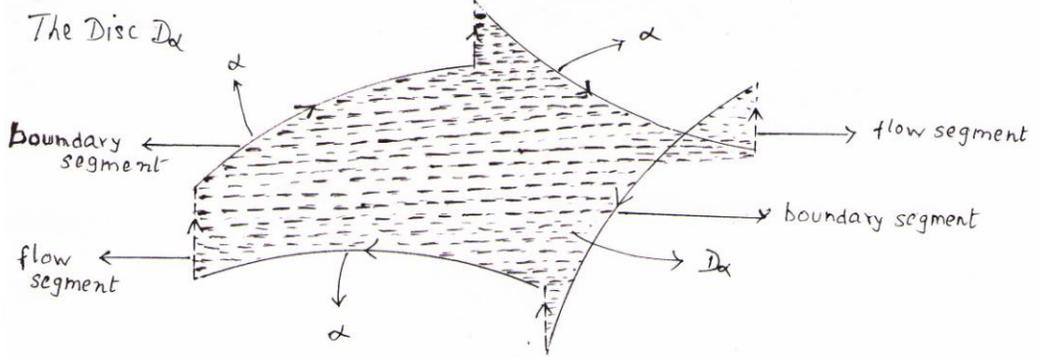

FIG 13(a)

The Disc $D_\alpha$

α  α

Boundary segment ← → flow segment

→ boundary segment

flow segment ← $D_\alpha$

α

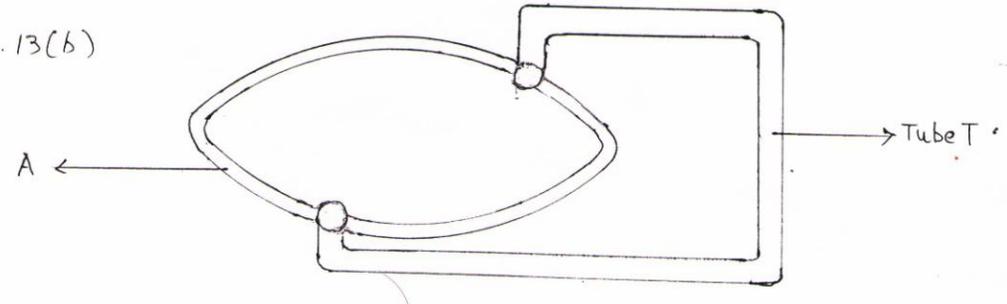

FIG. 13(b)

A ← → Tube T

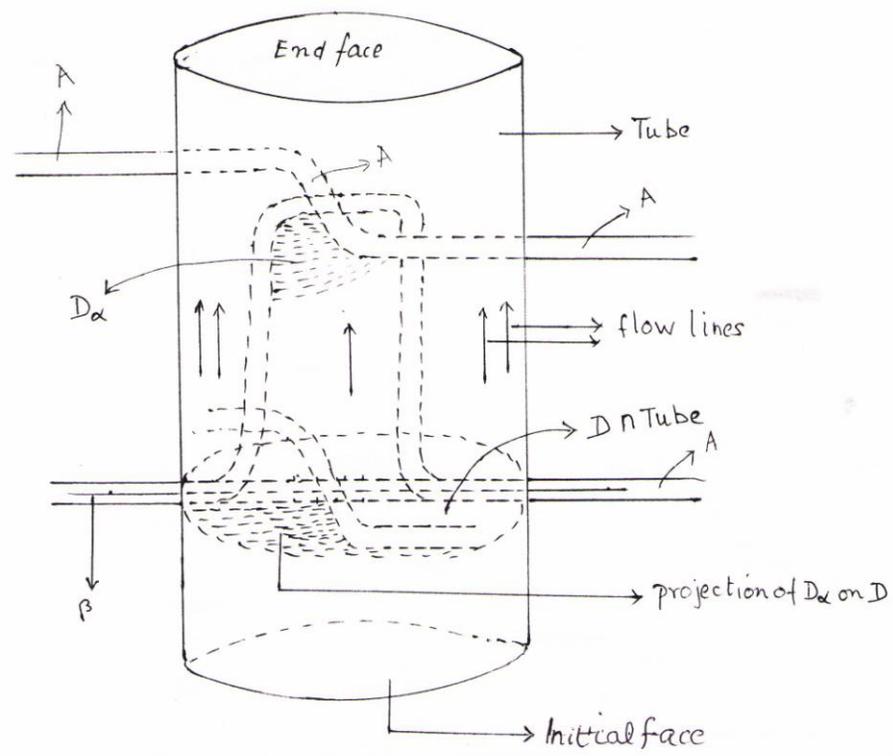

End face

A

→ Tube

λ

A

$D_\alpha$

↑↑  ↑  ↑↑ → flow lines

→ D ∩ Tube

A

→ projection of $D_\alpha$ on D

β

→ Initial face

FIG 14.

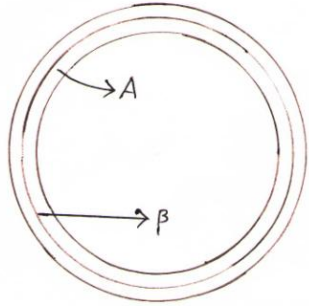

FIG. 15

A

β

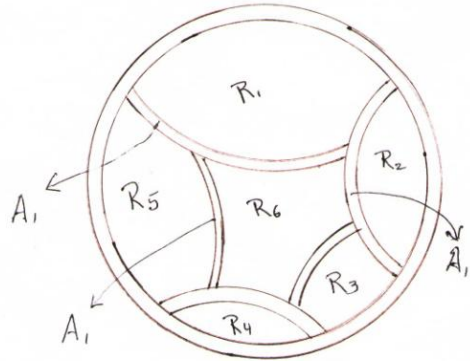

FIGURE 16

$R_1$

$R_2$

$A_1$

$R_5$

$R_6$

$A_1$

$R_3$

$A_1$

$R_4$

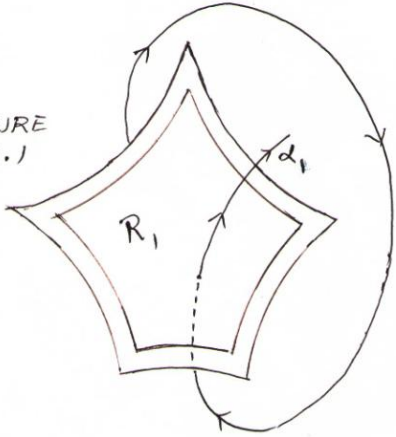

FIGURE 16.1

$\alpha_1$

$R_1$

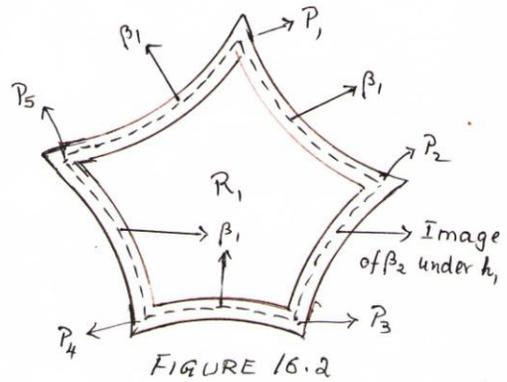

$\beta_1$

$P_1$

$P_5$

$\beta_1$

$P_2$

$R_1$

$\beta_1$

Image of $\beta_2$ under $h_1$

$P_4$

$P_3$

FIGURE 16.2

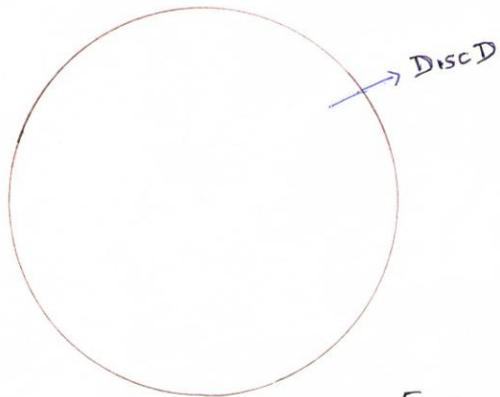

Disc D

(a)

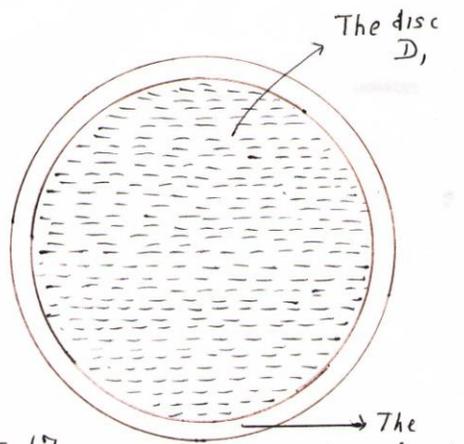

The disc $D_1$

The annulus A

FIGURE 17

(b)

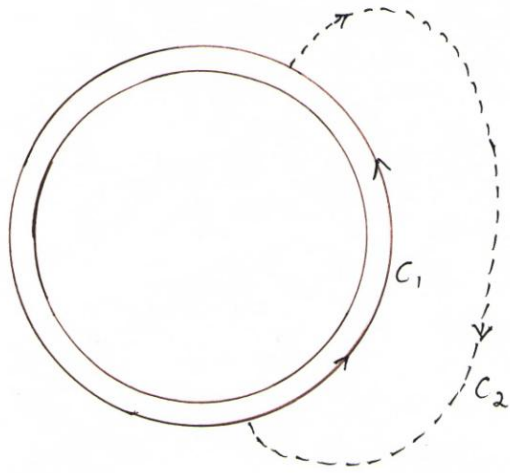

(c) $C_{12} = C_1 \cup C_2$ is a singular cycle of $D$ and $A$.

## FIGURE 17.

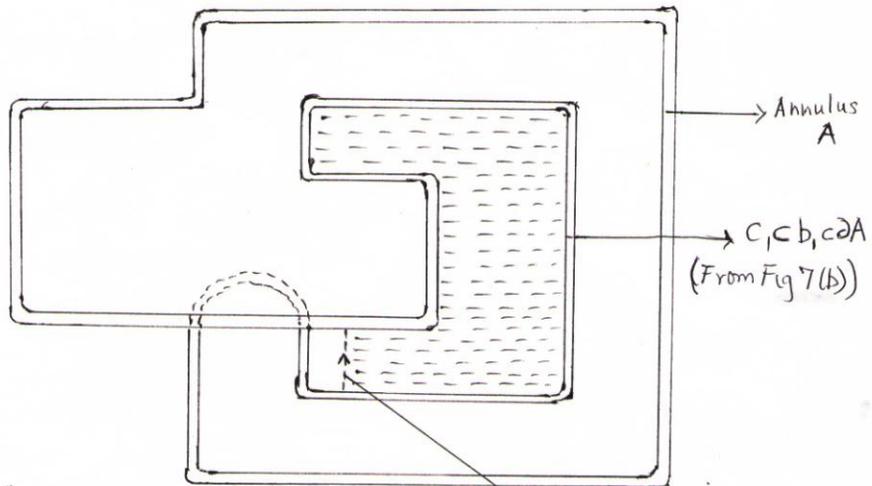

→ Annulus
   A

→ $C_1 \subset b_1 \subset \partial A$
   (From Fig 7(b))

(d)

$C_2$ (Flow segment)

Shaded region shows a singular
disc spanning $C_{12}$.

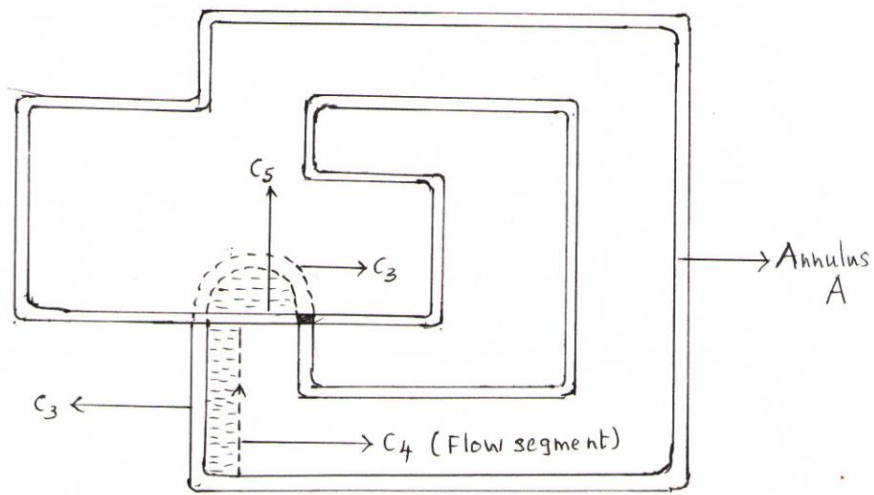

(e) $C_2$ is collapsed to a point. This gives rise to a new singular cycle $C_{22} = C_3 \cup C_4 \cup C_5$. Shaded part shows a singular disc spanning $C_{22}$.

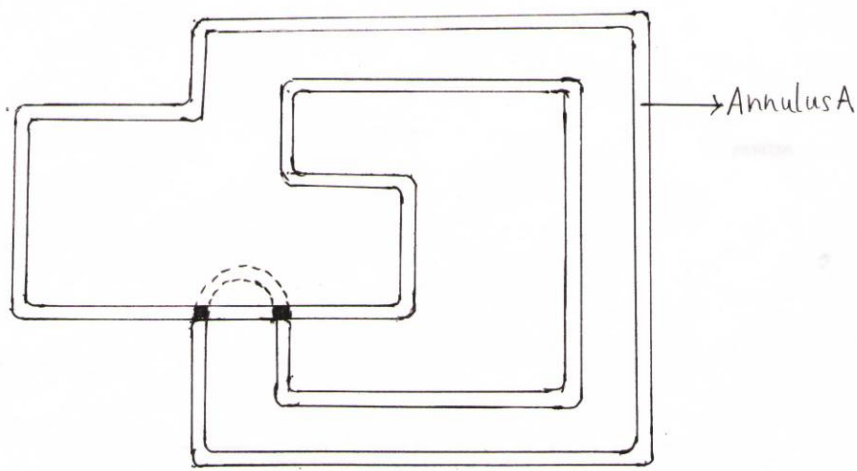

(f) In order to kill $C_{22}$, $C_4$ is collapsed to a point

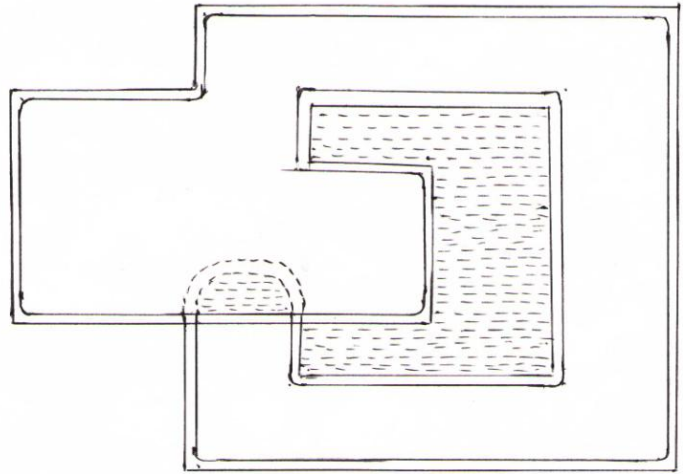

(g)

In the last step, the holes are sealed with discs. The shaded regions show the discs.

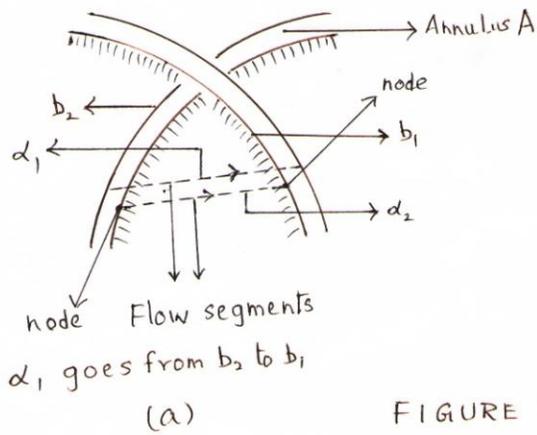

node    Flow segments
$\alpha_1$ goes from $b_2$ to $b_1$

(a)

FIGURE 18

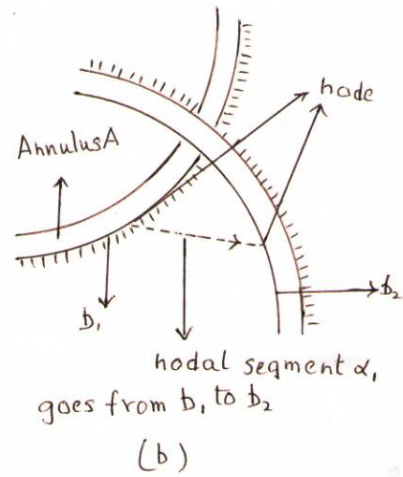

hodal segment $\alpha_1$
goes from $b_1$ to $b_2$

(b)

Shaded curve represents outer boundary $b_1$ of A

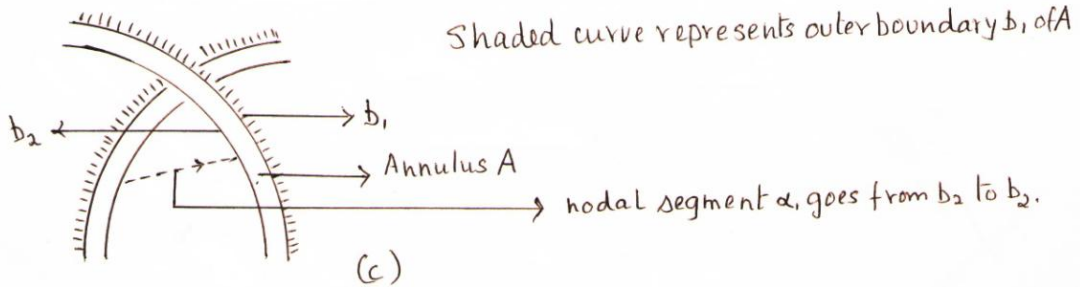

(c)

nodal segment $\alpha$, goes from $b_2$ to $b_2$.

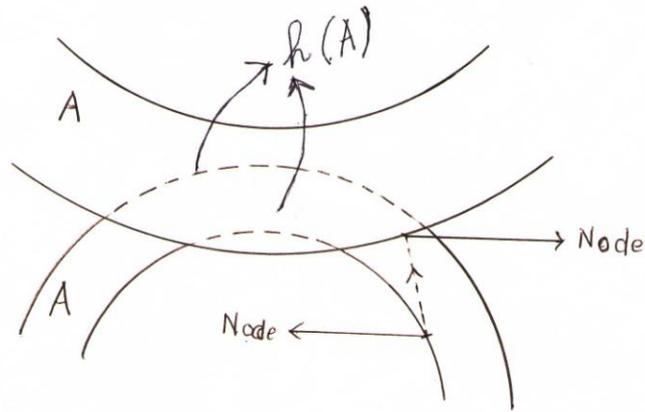

A

h (A)

A

Node

Node

FIGURE 18(d)

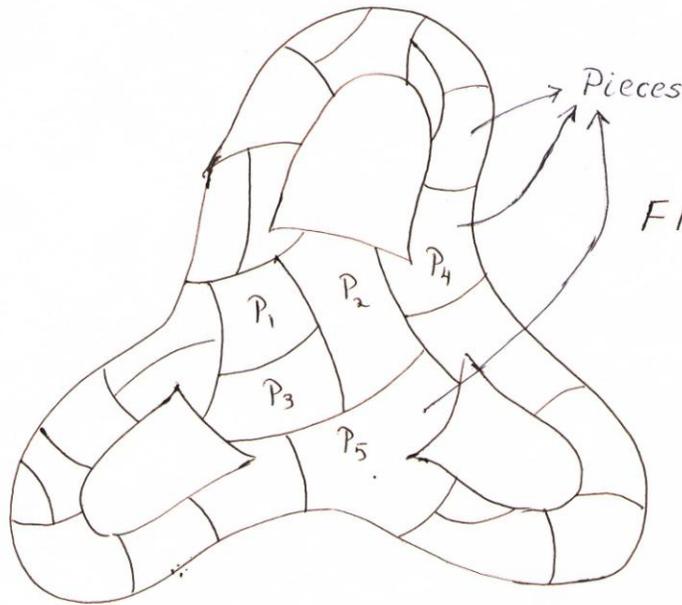

Pieces

FIGURE 19

$P_1$

$P_2$

$P_4$

$P_3$

$P_5$

The space X partitioned into pieces $P_i$ by the system of arcs Y

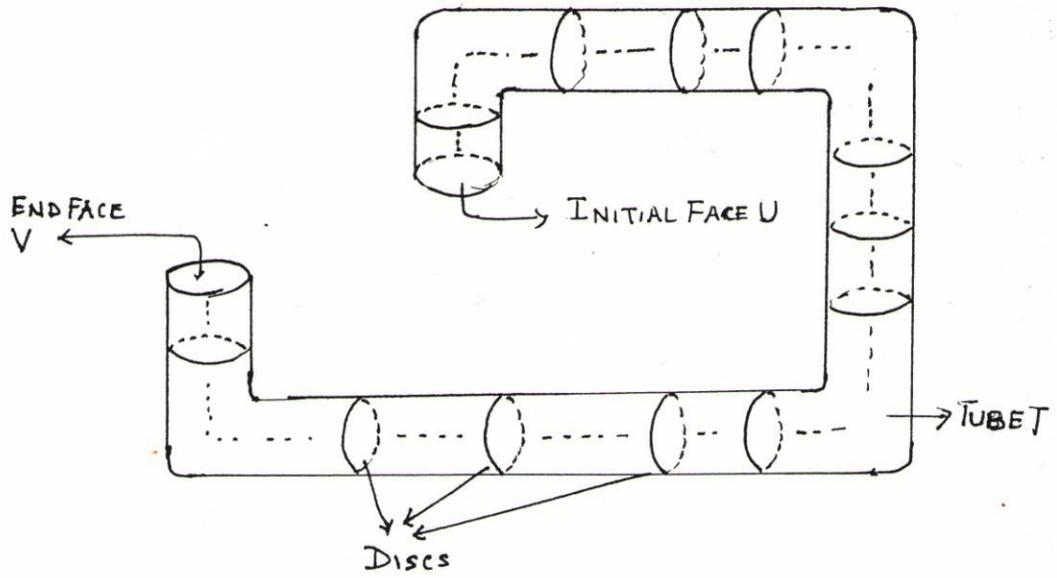

END FACE
V ←

→ INITIAL FACE U

→ TUBE T

Discs

FIGURE 20.

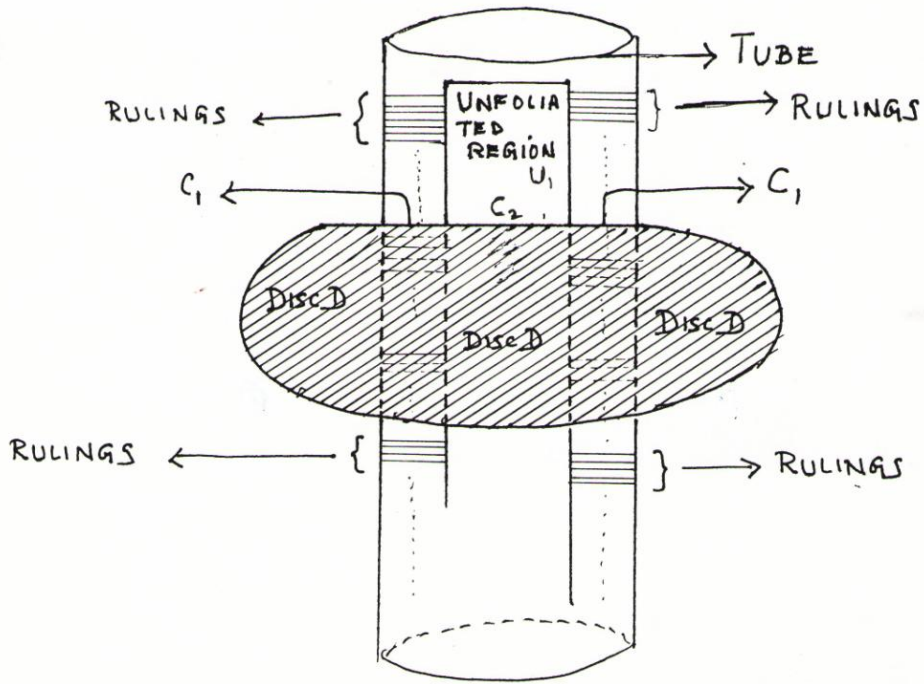

→ TUBE

RULINGS ←    { UNFOLIATED REGION $U_1$ }   → RULINGS

$C_1$ ←    $C_2$    → $C_1$

Disc D    Disc D    Disc D

RULINGS ←    {    } → RULINGS

FIGURE 21.

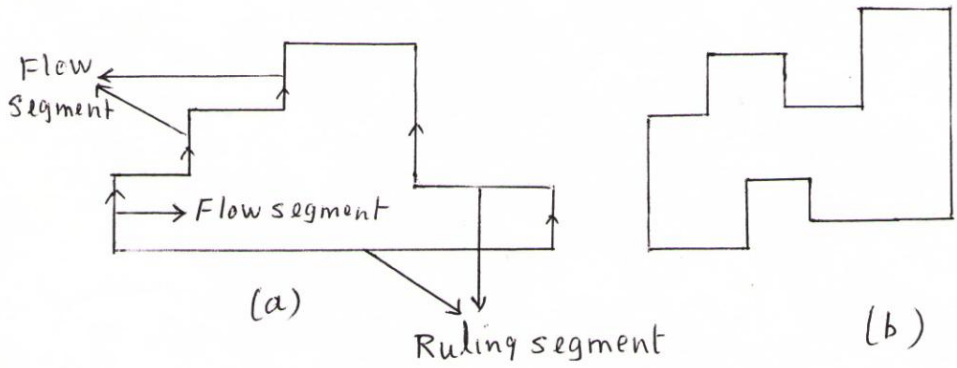

Flow Segment

Flow segment

Ruling segment

(a)

(b)

FIGURE 21.1

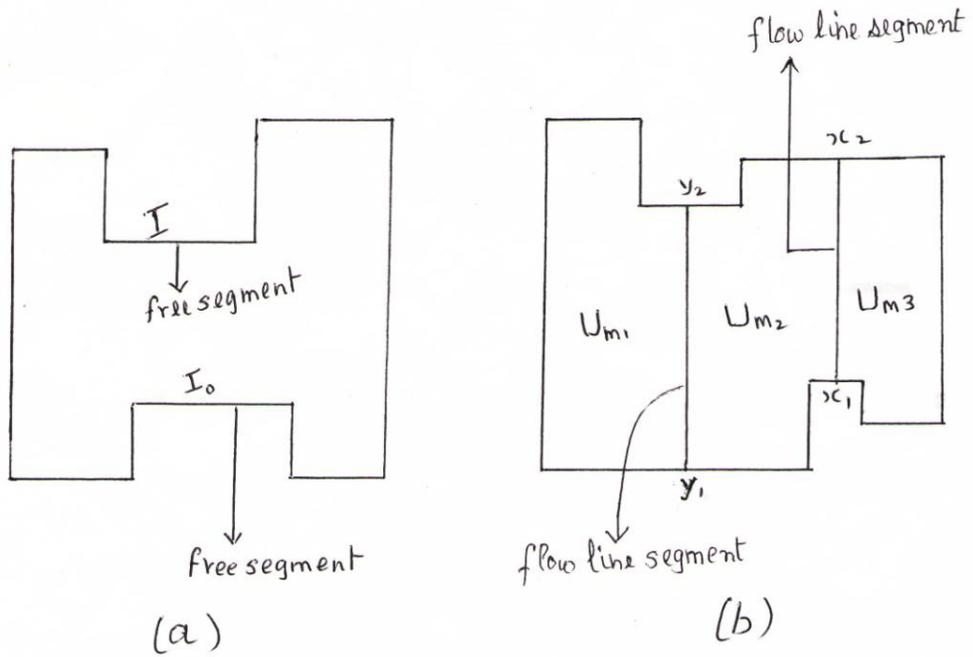

flow line segment

$I$

free segment

$I_0$

free segment

(a)

$y_2$

$x_2$

$U_{m_1}$

$U_{m_2}$

$U_{m_3}$

$x_1$

$y_1$

flow line segment

(b)

FIGURE 21.2

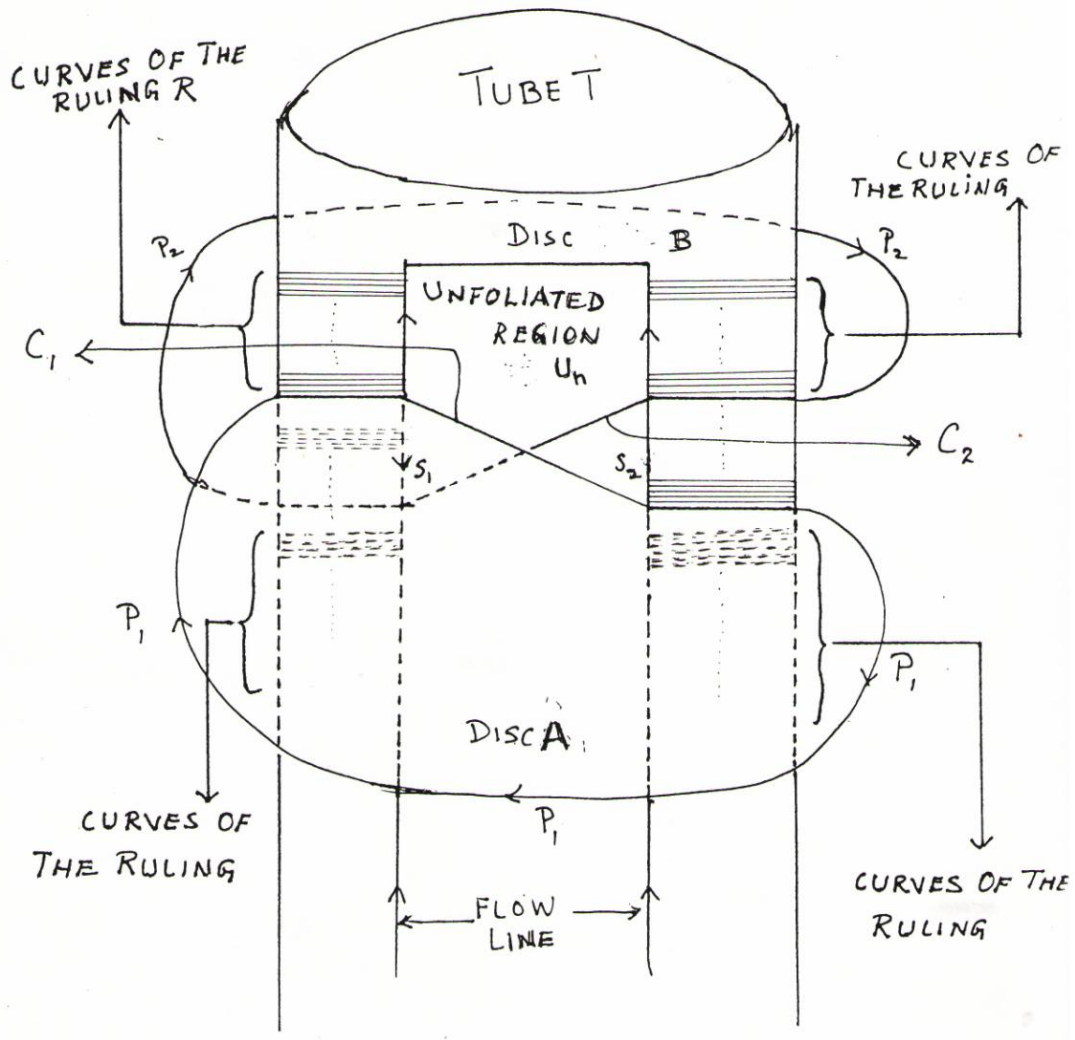

CURVES OF THE RULING R

TUBE T

CURVES OF THE RULING

$P_2$

Disc B

UNFOLIATED REGION $U_n$

$C_1$

$P_2$

$S_1$

$S_2$

$C_2$

$P_1$

$P_1$

CURVES OF THE RULING

DISC A

$P_1$

FLOW LINE

CURVES OF THE RULING

FIGURE 22

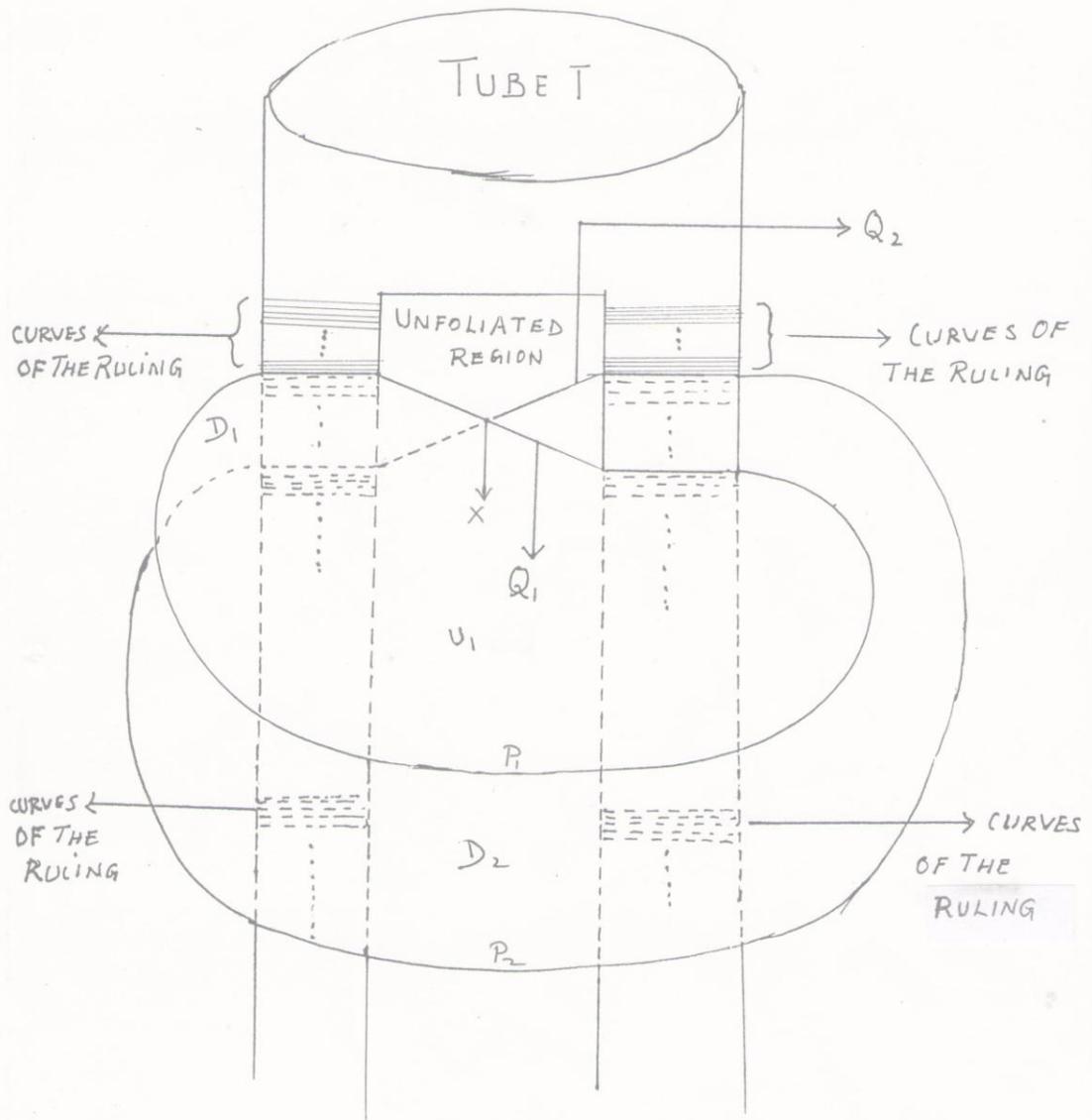

FIGURE 23

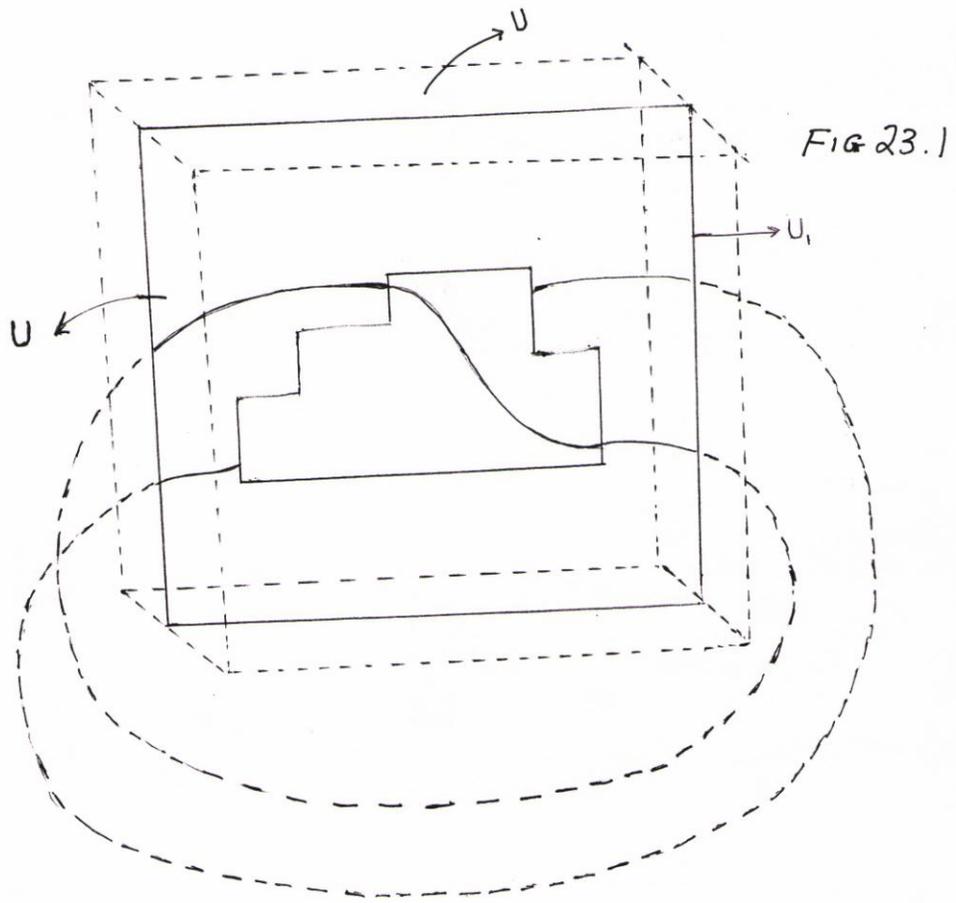

FIG 23.1

U

U

U₁

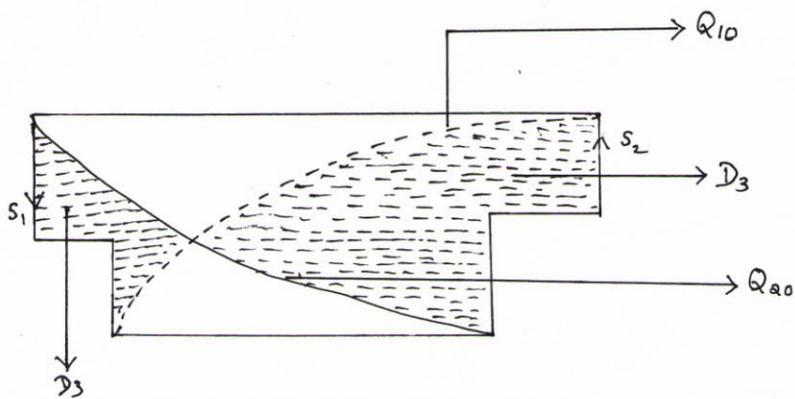

Q₁₀

S₂

D₃

Q₂₀

S₁

D₃

FIG.23.2

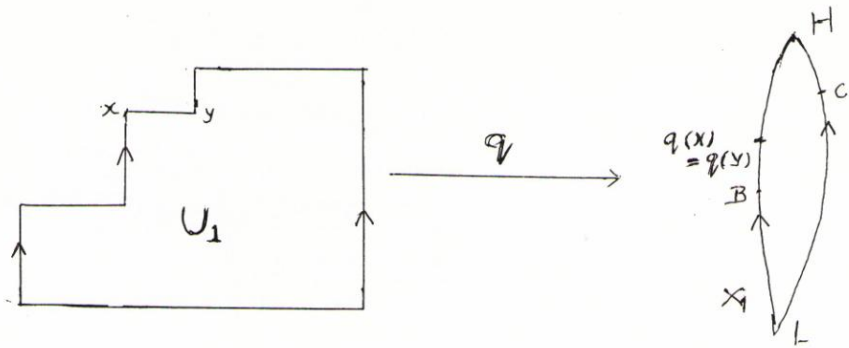

FIGURE 24.

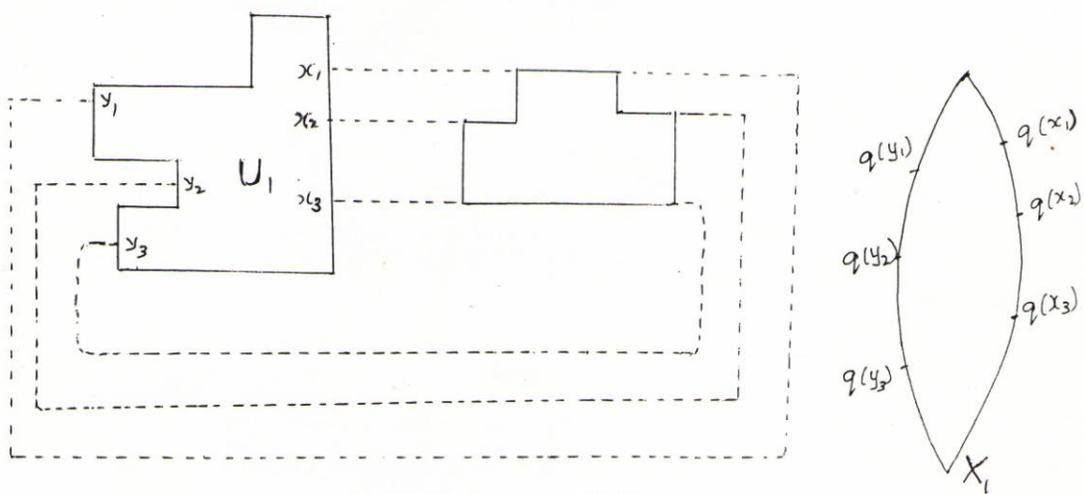

FIGURE 25.

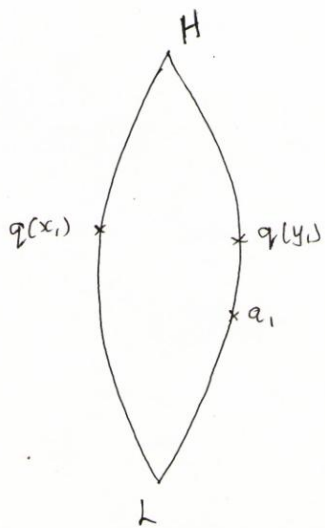

FIGURE 26.

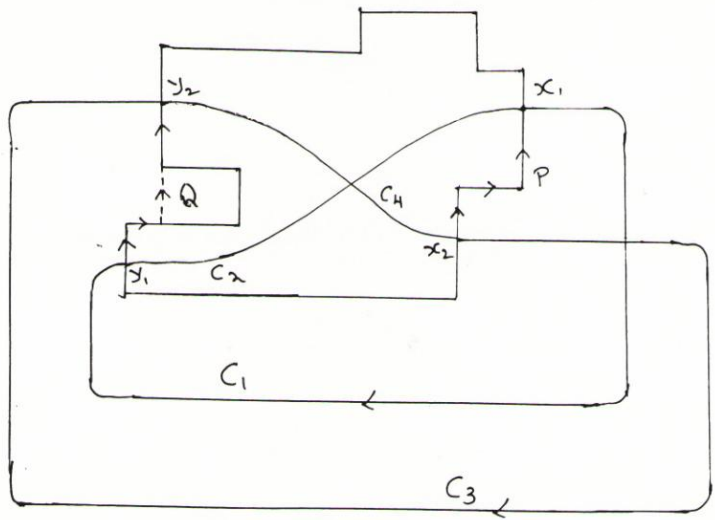

FIGURE 27

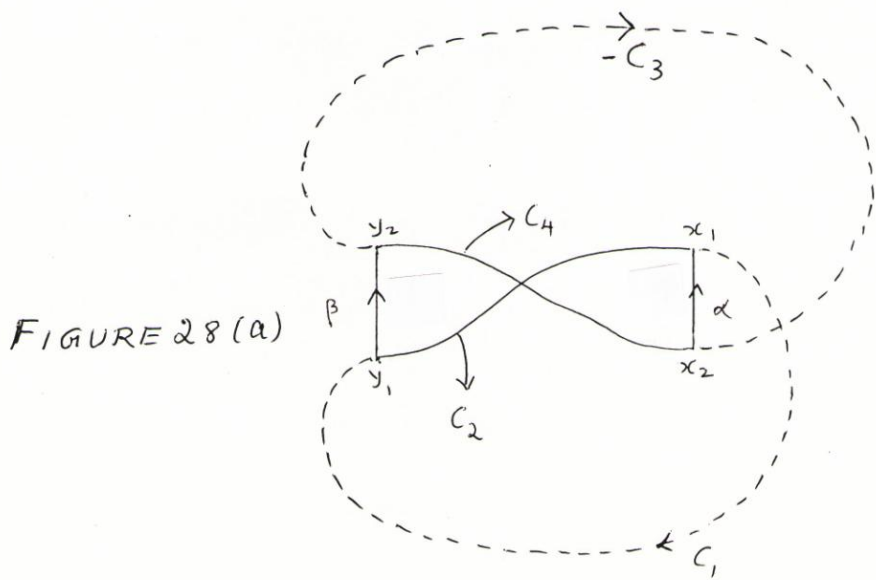

FIGURE 28 (a)

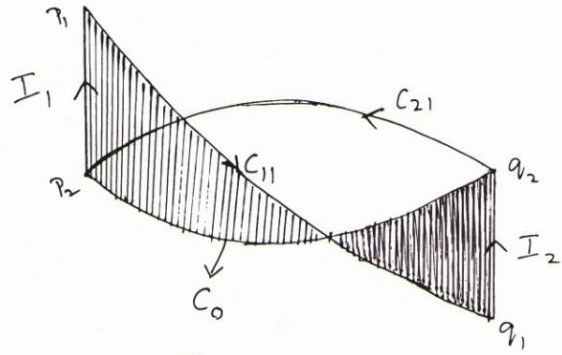

FIGURE 28 (b)

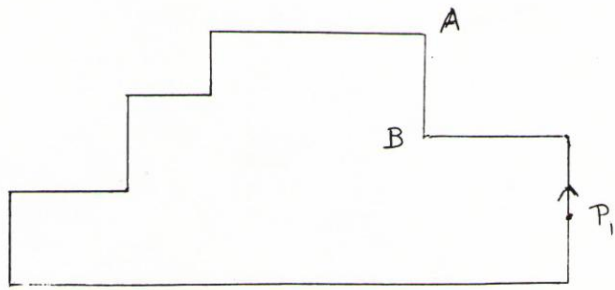

(a)

FIGURE 29

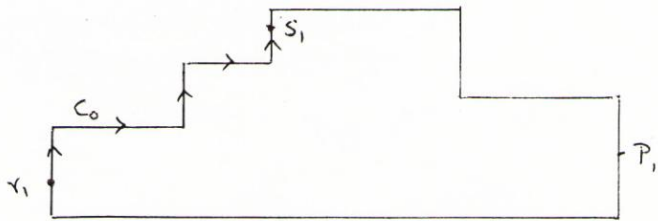

(b)

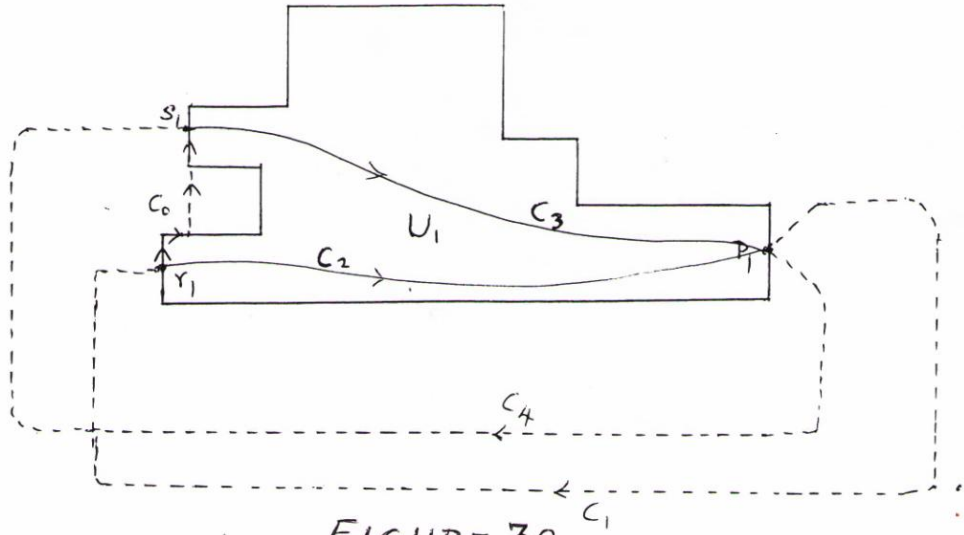

FIGURE 30

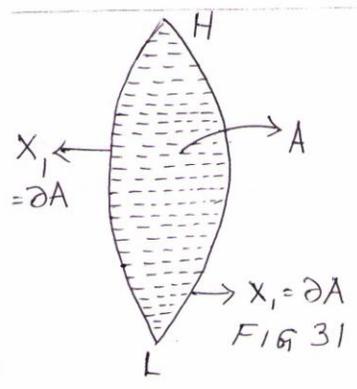

FIG 31

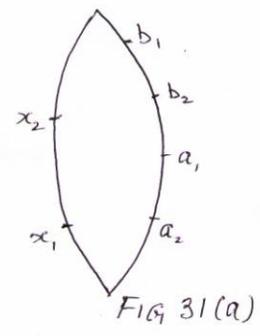

FIG 31 (a)

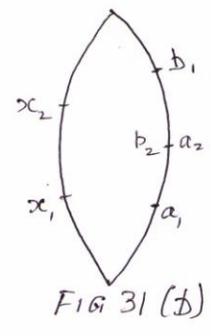

FIG 31 (b)

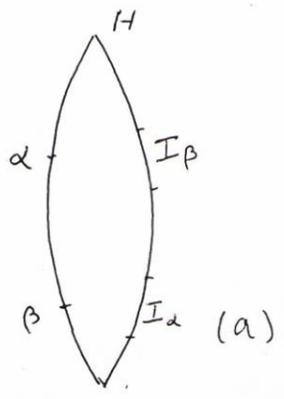

(a)

FIGURE 32

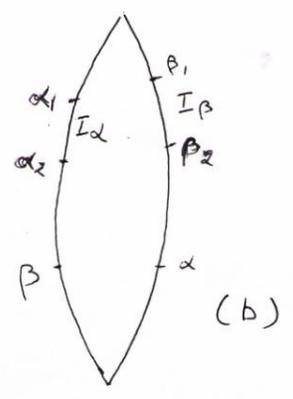

(b)

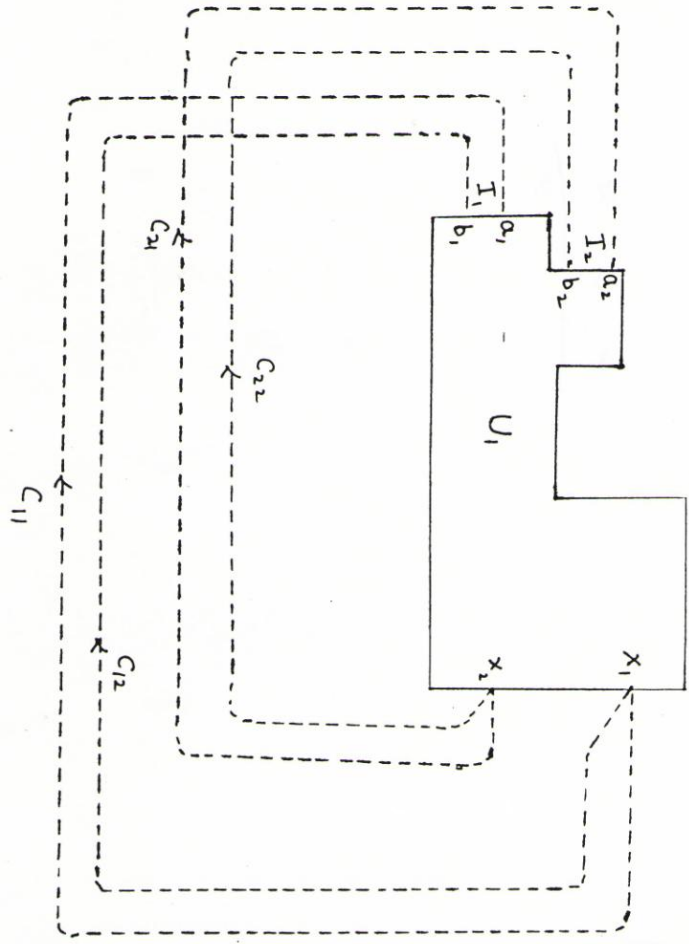

FIGURE 33 (a)

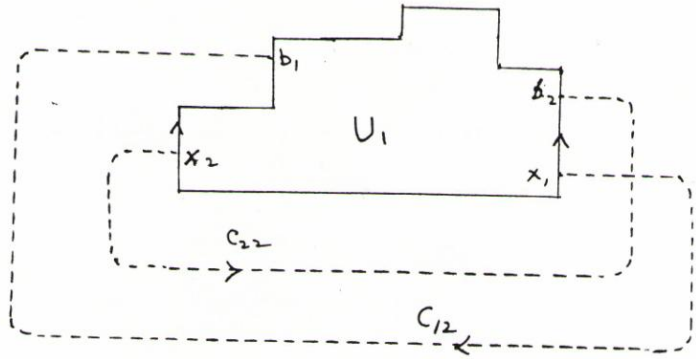

$$b_1 \qquad b_2$$
$$x_2 \qquad U_1 \qquad x_1$$
$$c_{22}$$
$$C_{12}$$

FIGURE 33 (b)

FIGURE 34

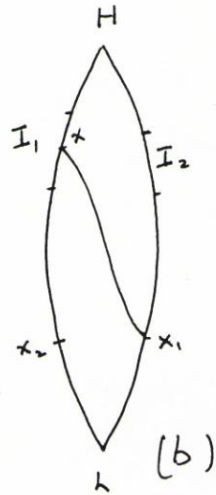

(b)

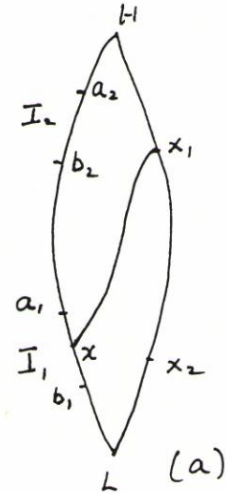

(a)

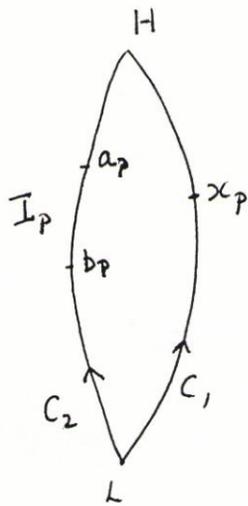

FIGURE 35.

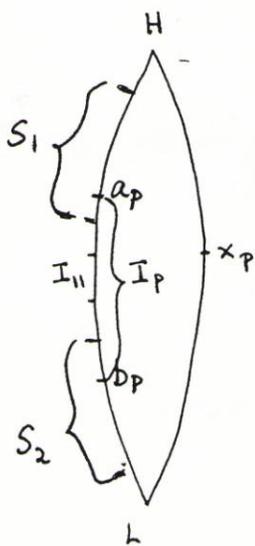

FIGURE 36

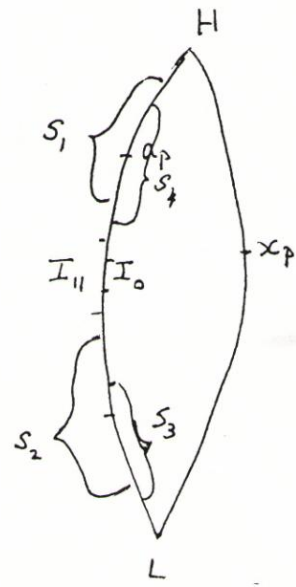

FIGURE 37.

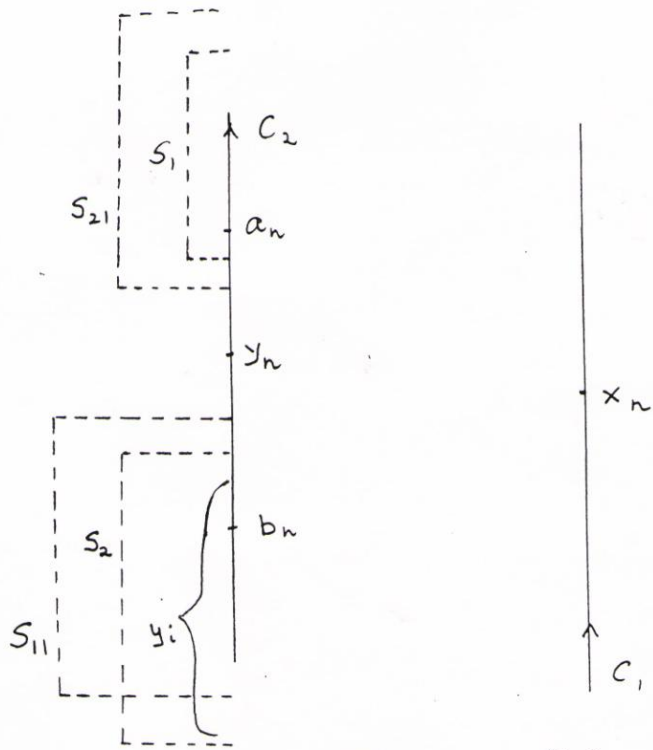

FIGURE 38

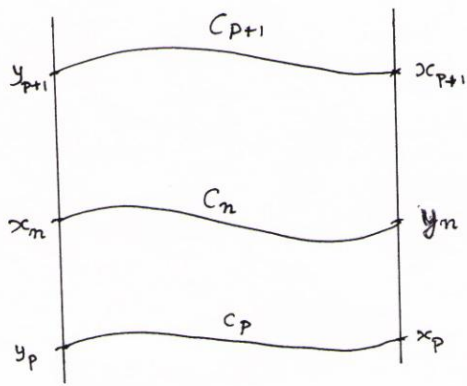

FIG 39.

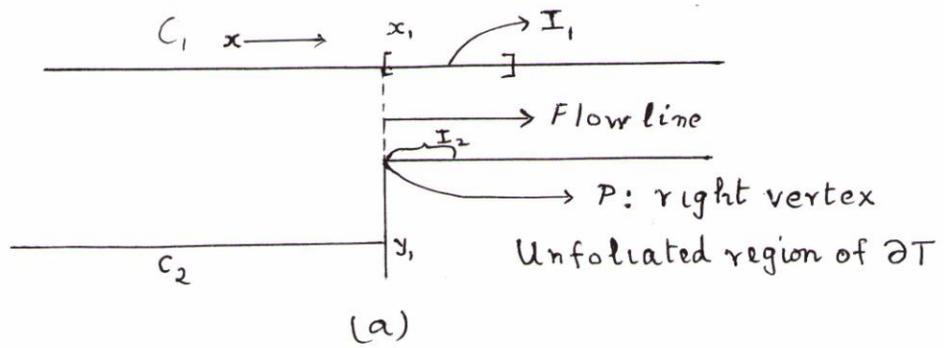

(a)

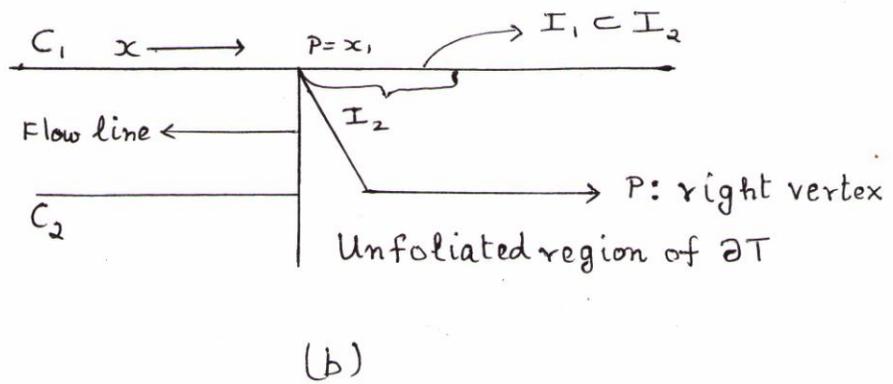

(b)

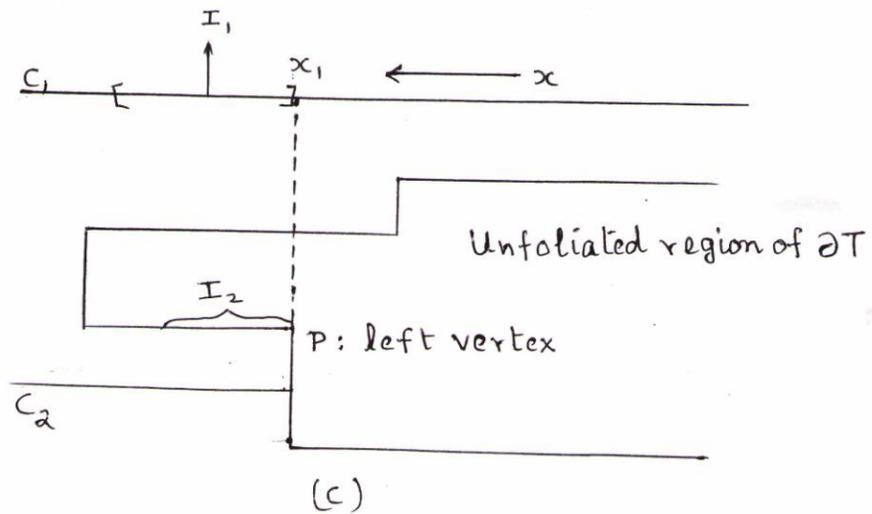

(c)

FIGURE 40 (a), (b), (c).

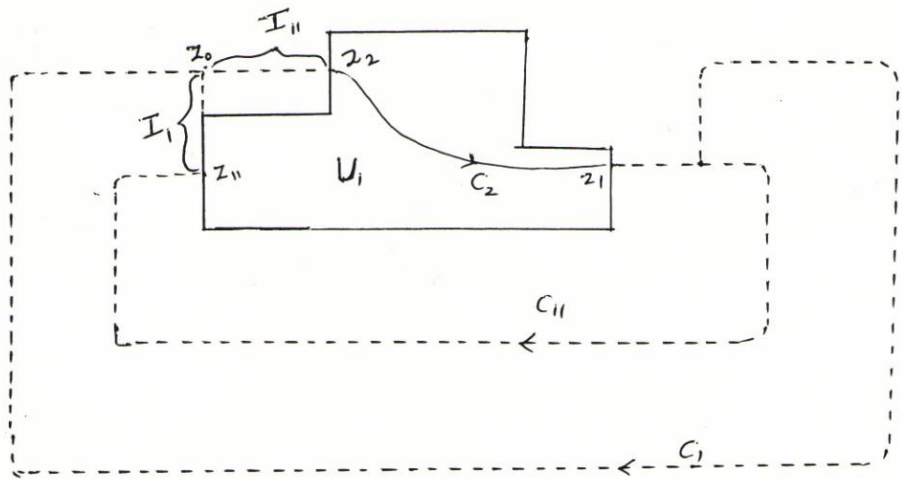

$$FIGURE \; 41 \, (a)$$

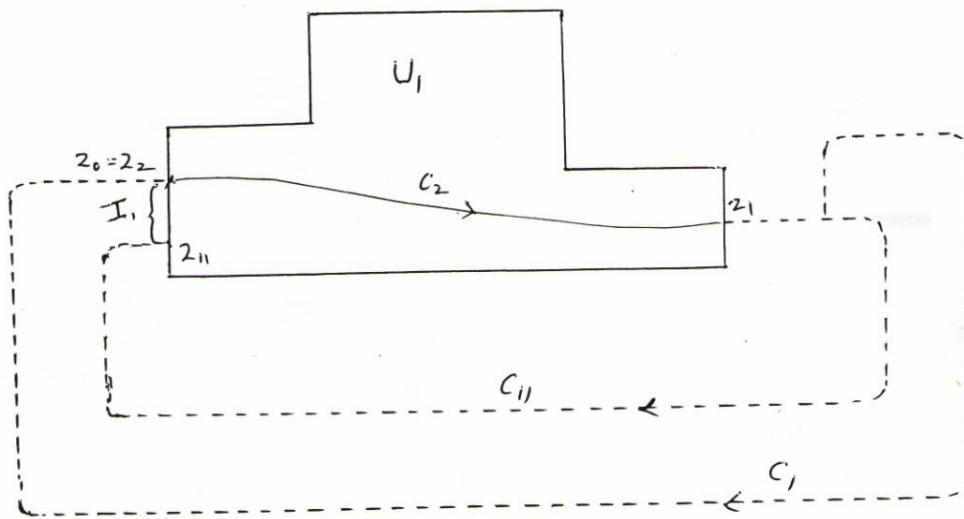

$$FIGURE \; 41 \, (b)$$

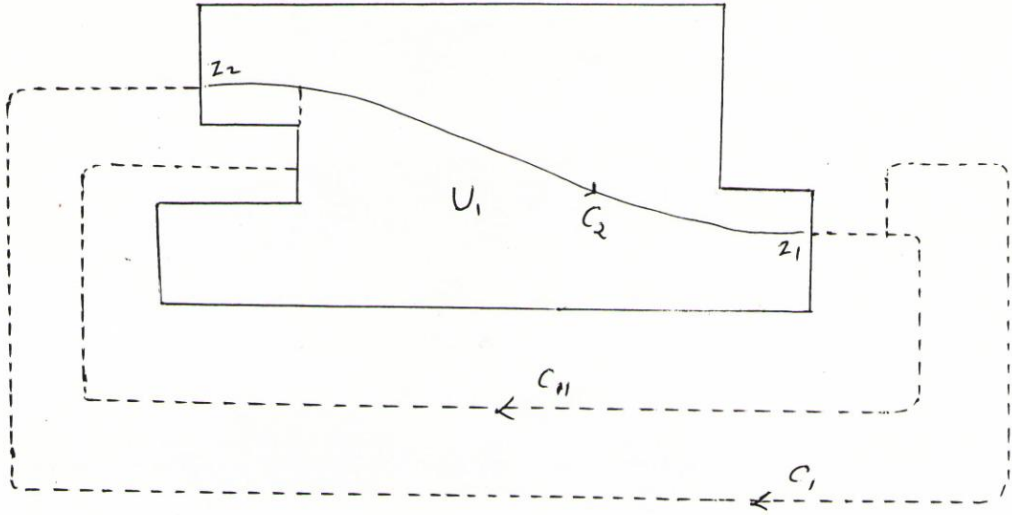

FIGURE 42

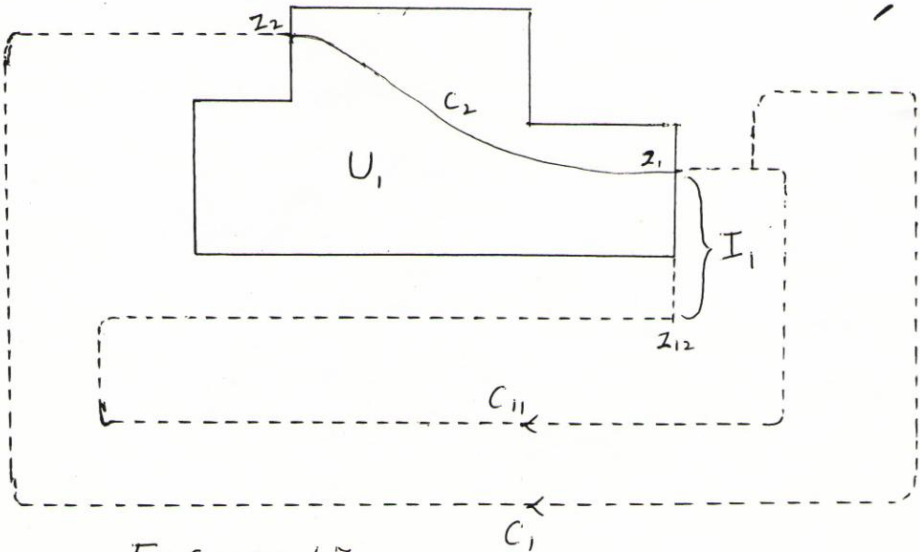

FIGURE 43

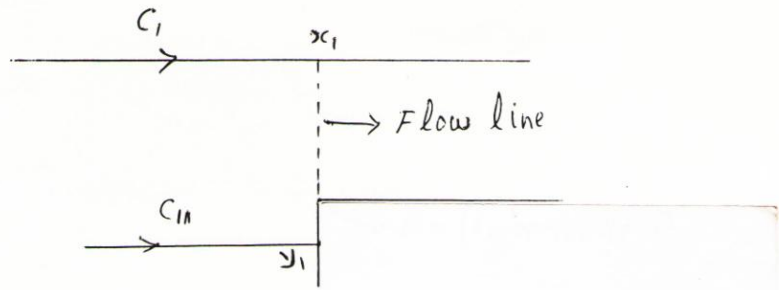

FIGURE 44

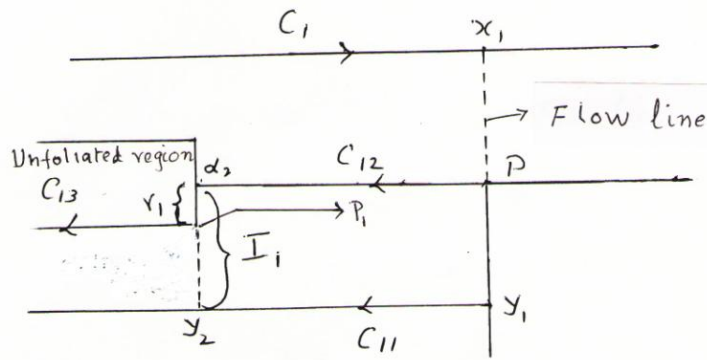

FIGURE 45

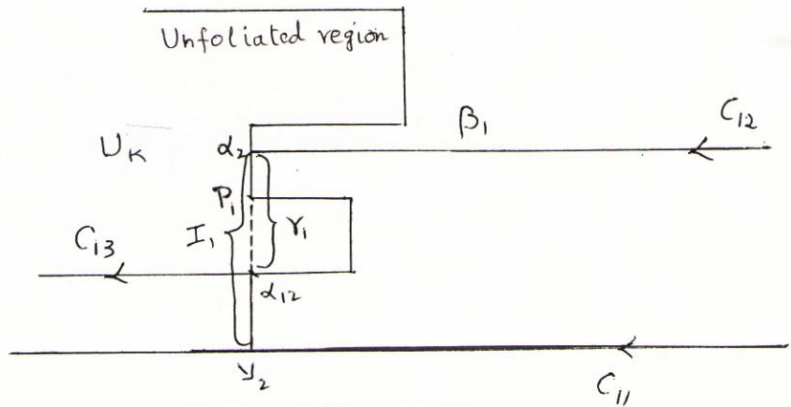

FIGURE 46

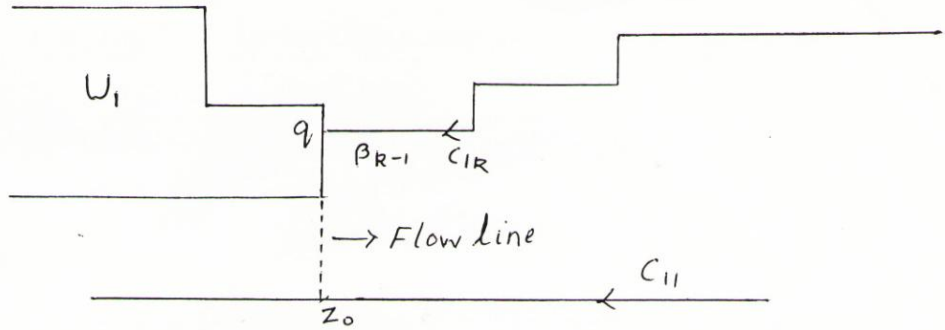

FIGURE 47

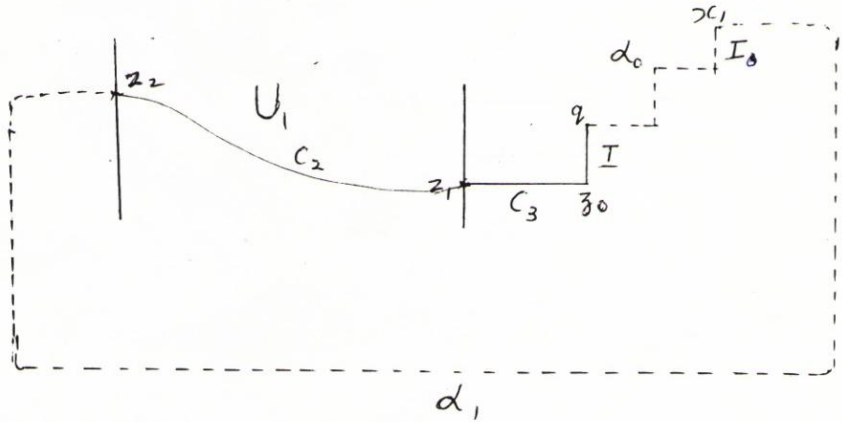

FIGURE 48

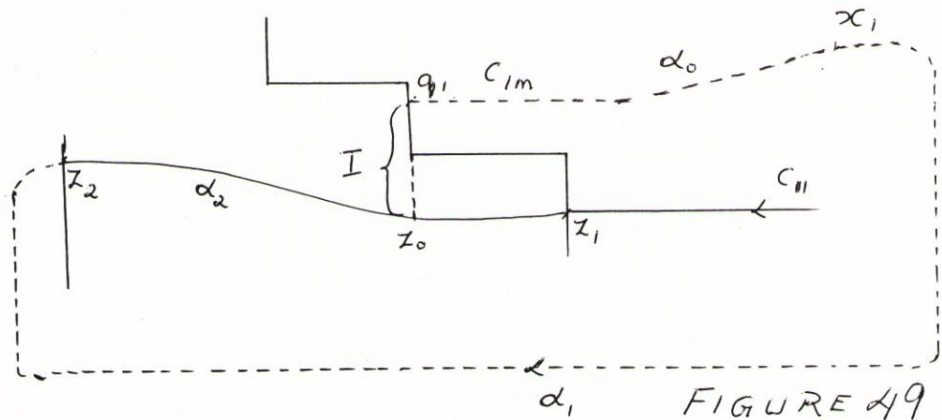

FIGURE 49

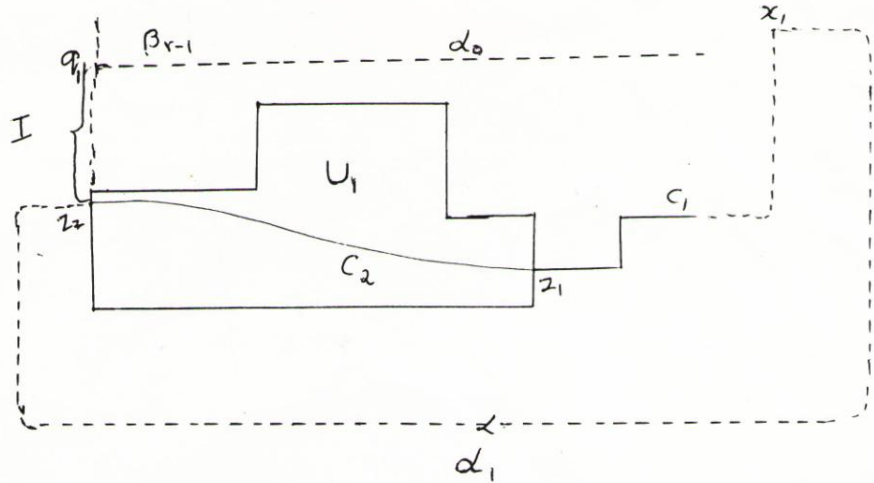

FIGURE 50